\documentclass[12pt,letterpaper]{article}

\RequirePackage{amsthm,amsmath}
\RequirePackage[numbers]{natbib}
\usepackage{a4wide,epsfig,latexsym,amsfonts, amssymb, enumerate,amsmath}
\usepackage{hyperref,relsize}
\usepackage[english]{babel}
\newtheorem{defn}{Definition}[section]
\newtheorem{lem}[defn]{Lemma}
\newtheorem{thm}[defn]{Theorem}
\newtheorem{cor}[defn]{Corollary}
\newtheorem{assu}[defn]{Assumption}
\newtheorem{exa}[defn]{Example}

\newtheorem{nota}[defn]{Notation}

\newtheorem{defi}[defn]{Definition}
\newtheorem{prop}[defn]{Proposition}
\let \epsilon=\varepsilon

\begin{document}

\title{Mimicking counterfactual outcomes to estimate causal effects}
%\runtitle{Mimicking counterfactual outcomes}

\author{Judith J.~Lok\\
  Department of Biostatistics, Harvard School of Public Health\\
  655 Huntington Avenue, Building 2, Room 409\\
Boston, Massachusetts 02115, US\\
  jlok@hsph.harvard.edu}

\maketitle

%Annals: <=35 characters, is 33 now

%Check style for numbering convention.

\begin{abstract} $\;$
  In observational studies, treatment may be adapted to covariates at
  several times without a fixed protocol, in continuous time.
  Treatment influences covariates, which influence treatment, which
  influences covariates, and so on. Then even time-dependent
  Cox-models cannot be used to estimate the net treatment effect.
  Structural nested models have been applied in this setting.
  Structural nested models are based on counterfactuals: the outcome a
  person would have had had treatment been withheld after a certain
  time. Previous work on continuous-time structural nested models
  assumes that counterfactuals depend deterministically on observed
  data, while conjecturing that this assumption can be relaxed. This
  article proves that one can mimic counterfactuals by constructing
  random variables, solutions to a differential equation, that have
  the same distribution as the counterfactuals, even given past
  observed data. These ``mimicking'' variables can be used to estimate
  the parameters of structural nested models without assuming the
  treatment effect to be deterministic.
%Annals: <= 150 words. 149 now. 150 Nov 10 2013. 150 Nov 2014.
% of such a time varying treatment
%, an assumption which has frequently been attacked and
%  has made people avoid counterfactuals altogether.
% These mimicking variables are the solution
%  to a differential equation with a final condition, which has one or
%  more parameters describing treatment effect. In this article we prove
%  this open problem.
%We hope that this will contribute to the discussion about causal reasoning.
\end{abstract}


\section{Introduction}
%Causal inference, counterfactuals and mimicking
%Counterfactuals, structural nested models, and mimicking
\label{intro}

%Except for in randomized experiments, statisticians usually restrict
%themselves to the modelling of association, and warn against using
%statistical models for causal inference. By lack of alternatives and
%since causal questions are often of great interest, models for
%association are, in spite of these warnings, often interpreted in a
%causal way. This leads \cite{Spir} to ask:
%``Why is so much of statistical application and so little of
%statistical theory concerned with causal inference?''

%A controversy among researchers in causal inference is whether or not
%to think in terms

Observational studies are no replacement for randomized clinical
trials, but they can be used, for example, where randomization is
unethical or to generate hypotheses for subsequent clinical trials. In an
observational study, treatment may be adapted to patient characteristics
which predict the outcome of interest.
This is called confounding by indication. If the
confounding by indication only takes place at baseline, one can
condition on initial person characteristics in order to get
meaningful estimates of the treatment effect. However, if the
confounding by indication also takes place after baseline, variables used for treatment decisions may be influenced by past
treatment. Thus they may themselves be indications of the treatment
effect, and in that case simply conditioning on them can lead to false
conclusions.

With such time-dependent confounding by indication, even the time-dependent
Cox model does not estimate the net effect of treatment (see e.g.~\cite{R87a}, \cite{Rpcp}, or~\cite{Rlat}).
 %\cite{R87a} Robins (1995), \cite{Rpcp} or Robins (1997)). \cite{Rlat}).
%, or \cite{Lok}.
With a time-dependent Cox model, the rate of events given past treatment and covariate history can be estimated, but the true parameter(s) on treatment may not reflect the treatment effect. A consistent estimator of the effect of the treatment on the outcome of interest has to take into account the effect of treatment on intermediate covariates. This is easily understood when considering a treatment which affects the outcome only because it affects an intermediate variable $L$. In that situation, if $L$ and treatment are both included in the time-dependent Cox model for the event of interest, the true parameter(s) on treatment in this Cox model equal $0$. However, treatment could be beneficial due to its effect on $L$. On the other hand, not including $L$ may also result in an inconsistent estimator, if $L$ predicts future treatment. This follows from the same reasoning as why, in case of non-randomized point treatment, one needs to adjust for predictors of both the treatment and the outcome to consistently estimate the treatment effect: if one does not adjust for $L$, and if persons with $L$ indicating a bad prognosis are more likely to be treated, the treatment may seem to adversely affect the outcome, even if it has no effect on anyone. To conclude, with time-dependent confounding by indication, one needs to take confounders into account, but adding the confounders to an outcome model is not enough.

If all confounders are measured (see Assumption~\ref{nuc} below),
structural nested models, proposed in \cite{R89,R92,Enc}, %\cite{R89,R92,Enc},
and marginal structural models, proposed in \cite{MSM1,MSM2},
%I checked that R89 has SNFTMs. R86 and R87 do not have them it seems to me.
can be used to consistently estimate treatment effects in the presence of
time-dependent confounding by indication. Structural nested models and marginal structural models make a distinction between the effect of the treatment and the reason why the treatment was given, by separately modeling the treatment decisions and the treatment effect.
\cite{versus}
compares structural nested models and marginal structural models. The current article focuses on structural nested models.

Structural nested models model relations between counterfactual
outcomes. We allow for general treatment regimes. Consider a single person, who received a
particular treatment regime with outcome $Y$. For example, the particular treatment regime could
be as follows: first, no treatment, then, after a certain time, initiation of treatment, then, the dosage
changed some time thereafter, then treatment stopped, initiated again, etcetera. Had the person's
treatment been stopped (prematurely) at time $t$ and not been re-initiated thereafter (or, had treatment changed to a
``baseline'' treatment regime from time $t$ onwards), the person's outcome, $Y^{(t)}$, might have been
different. Since $Y^{(t)}$ is generally not observable, it is a
counterfactual outcome. In a discrete-time setting, \cite{RJ} show that existence of counterfactuals places
  no restrictions on the distribution of the observed variables. No comparable proof
  exists for the continuous-time case.

%The principal objection to the use of counterfactuals rests on the
%impossibility to determine exactly what is ``the closest possible
%world to this one in which ...'' (the treatment would have been
%different). The rebuttal is that carrying out the thought experiment
%``what would have happened if'' to make the prediction ``what would
%happen if'' is exactly what scientific research is all about. Of
%course, doing a good thought experiment requires a clear understanding
%of the context.

%\cite{Pea} and \cite{Spir} study causality using graphs that can also
%be expressed in terms of counterfactuals.
%, and I am unable to avoid seeing counterfactuals behind those graphs.
%I consider counterfactuals a very natural and useful framework to
%think about causal questions, and have not made attempts to avoid
%them.

An important controversy in the causal literature is that
counterfactuals are often assumed to depend deterministically on the
observed data: given the model and the parameter values, all
counterfactual outcomes $Y^{(t)}$ for each person can simply be
calculated from the observed data. \cite{Enc} calls this local rank
preservation (in most cases, this implies global rank preservation),
when the counterfactual outcomes are solutions to the differential
equation (\ref{Xdef}) in Section~\ref{cmim} below. Treatment is then
said \emph{not} to affect the outcome of interest if the outcome for
any particular person would have been exactly the same regardless
which treatment was given. The assumption of deterministic dependence
is related to the assumption of constant effect in \cite{Hol}:
%Section 4.4
that is, the difference between counterfactual outcomes belonging to
different treatments is a constant identical for all persons.

The assumption of deterministic dependence/ (local) rank preservation
has frequently been attacked.
%, and has made some people avoid counterfactuals altogether
This assumption does not hold if, for example, two persons with the
same observed data (e.g., both receiving some prophylactic drug) could
have had a different outcome had they not been treated starting from
some time $t$ (e.g., one might have contacted a virus and the other
might not). In addition, deterministic dependence can never be tested,
with only one outcome observed for each person. For these two
reasons, the assumption of (local) rank preservation should be avoided
if at all possible.

%Structural Nested Models use the observed data to define a
%differential equation with a final condition. For each person we
%construct the solution to that differential equation with final
%condition, called $X(t)$ (where $t$ stands for time). Then we show
%that this $X(t)$ has the same distribution as the counterfactual
%outcome with treatment as given until time $t$ and no treatment after
%$t$, $Y^{(t)}$, given all covariate- and treatment information known
%at time $t$. This remarkable result is Theorem~\ref{thm1} in the
%current article. It is the cornerstone for use of Structural Nested
%Models without assuming constant treatment effect.

%To estimate the treatment effect we need that there is no unmeasured
%confounding: all prognostic factors for the outcome of interest which
%influence treatment decisions have to be observed. For the main result
%of this article, $X(t)$ has the same distribution as $Y^{(t)}$ given
%past treatment and covariates, we do not need to assume no unmeasured
%confounding. However, without no unmeasured confounding the result of
%this article cannot be used to estimate treatment effect or test
%whether treatment affects the outcome, and is thus of no interest.

%Structural nested models express treatment effect in terms of the
%difference (in the distribution, as we show here) of various
%counterfactual outcomes. This article proves a long standing conjecture
%from \cite{Enc} about structural nested models in continuous time.

In discrete time, when treatment and covariates change at fixed times
which are the same for all persons, the theory of structural nested
models is well developed. \cite{SNart} prove that it is not necessary
to assume a deterministic treatment effect.  In order to do so, they
show that a certain ``blipped down'' outcome $X(t)$ mimics the
outcome $Y^{(t)}$ had treatment been withheld from time $t$ onwards, in the
sense that $X(t)$ has the same distribution as $Y^{(t)}$ given past
treatment and covariate history. They also indicate why the resulting
estimators for treatment effect are consistent and asymptotically
normal.

However, in reality covariates and treatment often change in continuous
  time. Moreover, in discrete time the interpretation of the treatment
  effect (shift- or blip function) depends on the time scale
  chosen. In continuous time, the treatment effect (infinitesimal
  shift function) can often be interpreted as speed or rate. For
  these reasons, \cite{R92, Aids, smoke, Keiding, Kei} have applied
  continuous-time structural nested distribution models. However,
  because of a lack of theory for these models, the applications have
  relied on the assumption of (local) rank preservation. The models
  fitted in \cite{Aids, Keiding, Kei} are described in Examples~\ref{DAexa2}
and~\ref{DAexa}.
Section~\ref{nuc} or, in greater detail, \cite{ASarx} describes
  how to use the results in the current article in order to show that assuming
  (local) rank preservation is not necessary to estimate treatment effects with
  structural nested models (an example can be found in Section~\ref{sim}).
  Therefore, the main contribution of the current article is showing that the methods in
  \cite{R92, Aids, smoke, Keiding, Kei,Enc} are robust to failures of the assumption of (local) rank
  preservation.

Structural nested models in continuous time are meant to estimate the
effect of a continuous treatment, for which the effect of a small
duration is small. \cite{Enc} conjectures that the appealing large
sample properties of discrete-time structural nested models extend to
continuous time; however, his proof requires the assumption of (local)
rank preservation. He conjectures that (i) also without (local) rank
preservation, a certain ``blipped down'' outcome $X(t)$ has the same
distribution as $Y^{(t)}$ given past treatment and covariate history,
(ii) the resulting estimators
%solve unbiased estimating equations and
are consistent and asymptotically normal, and (iii) for certain
models, estimators and confidence intervals can be calculated with
standard software, used in a non-standard way. This article proves
conjecture (i), which we call mimicking counterfactual outcomes, and
explains why such a subtle result is true.
% unverified statistical claims \cite{ASarx, plss} proves conjectures (2) and (3), using the result
\cite{ASarx} proves conjecture (ii), using conjecture (i). \cite{plss}
proves conjecture (iii), using a partial likelihood approach and
conjectures (i) and (ii). Thus, the current article fills the final link in
this methodology to estimate treatment effects of time-varying
treatments in longitudinal observational studies without relying on (local) rank preservation.
This methodology can be applied to longitudinal observational data, to study the effects
of interventions affecting, for example, economic and health outcomes.
%answer many important public health issues.
%provides an essential building block in the theory of structural
%nested models in continuous time.
%is essential for

This article is organized as follows. Section~\ref{sets} introduces
the setting and notation of this article. Section~\ref{model}
introduces the model for treatment effect, and shows some examples.
%The model has one or more parameters describing treatment effect.
Section~\ref{cmim} defines the mimicking variables $X(t)$ as the
solution to a differential equation with a final condition.
Section~\ref{cmim} also states the main result of this article: $X(t)$
mimics $Y^{(t)}$ in the sense that it has the same distribution as
$Y^{(t)}$, even given past treatment- and covariate history.
Section~\ref{nuc} formalizes the assumption of no unmeasured
confounding, which as shown there is needed to use the result of the
current article to estimate the treatment effect. Section~\ref{nuc}
also indicates how, using the mimicking result, tests and estimators
can be developed without assuming (local) rank preservation.
Section~\ref{XYin} outlines the proof of the main result of this
article: $X(t)$ mimics $Y^{(t)}$ in the sense that it has the same
distribution as $Y^{(t)}$, even given past treatment- and covariate
history. Section~\ref{XY} proves the main result of this article for
non-survival outcomes $Y$. Section~\ref{XT} proves the main result for
survival time outcomes $Y$. Section~\ref{sim} describes a simulation study. Section~\ref{ext} concludes this article with a
discussion.

\section{Setting and notation}
\label{sets}

The setting to which continuous-time structural nested models apply
is as follows. The outcome of interest is a continuous real-valued
variable $Y$. For example, $Y$ is a person's survival time, time to clinical
AIDS, the number of white blood cells, or the CD4 count. Our objective is to estimate the effect of treatment
on $Y$. In this article, we consider a fixed time interval $t\in[0,\tau]$ with finite $\tau$, where $t=0$ is the time at which follow-up of interest starts (for example, $0$ could be the time of enrollment in a study, or a baseline time). During the time interval $[0,\tau]$,
treatment and person characteristics are observed for each person. $Y$
is measured at or after time $\tau$, or, in the case of a survival
time outcome, $Y$ could be measured before time $\tau$ if the person dies
before time $\tau$. We assume that treatment starts at or
after time $0$. We suppose that after time $\tau$, treatment is stopped
or switched to some kind of baseline treatment regime. Most of this article assumes
that there is no censoring, and $Y$ is observed for every person in
the study. Section~\ref{cens} incorporates right
censoring.

The covariate process describes the course of the disease of a
person, e.g.\ the course of the blood pressure and the white blood
cell count. The covariates which \emph{must} be included are those
which both (i)~influence a doctor's treatment decisions \emph{and}
(ii)~predict a person's prognosis with respect to the
outcome of interest.  If such covariates are not observed the
assumption of no unmeasured confounding, see Section~\ref{nuc}, will
not hold.

Denote the probability space by $\left(\Omega,{\cal F},P\right)$.  For
the moment consider a single person. Write $Z(t)$ for the covariate-
\emph{and} treatment values at time $t$. This article assumes that
$Z(t)$ takes values in $\mathbb{R}^m$, and that
$Z(t):\Omega\rightarrow \mathbb{R}^m$ is measurable for each
$t\in\left[0,\tau\right]$. Moreover, we assume that $Z$, seen as a
function on $\left[0,\tau\right]$, is continuous from the right with
limits from the left (cadlag), and that with probability one this
function, or ``sample path'', has only finitely many jumps. We also
assume that the probability that the covariate- and treatment process
$Z$ jumps at time $t$ equals $0$ for every fixed time $t$ (except
possibly for finitely many fixed times $t$, which could have point masses). For
example, the hazard of the jumps of the treatment process could be
continuous for all $t$, and could follow a continuous parametric distribution.
$\overline{Z}_t=\left(Z\left(s\right): 0\leq s\leq t\right)$ denotes
the covariate- and treatment history until time $t$, and
$\overline{{\cal Z}}_t$ is the space of cadlag functions from
$\left[0,t\right]$ to $\mathbb{R}^m$ in which $\overline{Z}_t$ takes
its values. Similarly, $\overline{Z}$ denotes the complete covariate-
and treatment history of the person in the interval
$\left[0,\tau\right]$, and $\overline{{\cal Z}}$ is the space in which
$\overline{Z}$ takes its values. In this article, the $\sigma$-algebra
on $\overline{{\cal Z}}_t$ and $\overline{{\cal Z}}$ is the projection
$\sigma$-algebra; measurability of $Z(s)$ for each $s\leq t$ is then
equivalent to measurability of the random variable $\overline{Z}_t$.

Counterfactual outcomes were already mentioned in the introduction.
$Y^{(t)}$ is the final outcome had treatment been stopped (prematurely) at time $t$
and not been re-initiated thereafter (or changed to
some kind of baseline treatment regime $\overline{0}$ from time $t$ onwards).
This article supposes that all counterfactual outcomes $Y^{(t)}$, for
$t\in\left[0,\tau\right]$ and for each person, are random variables on
the probability space $\left(\Omega,{\cal F},P\right)$. We assume
that observations and counterfactual outcomes of different persons are
independent and identically distributed, and are a random sample from
a larger infinite population of interest. For notational convenience,
we suppress the subscript $i$ for person.

\section{Model for treatment effect}
\label{model}

Structural nested models in continuous time model distributional
relations between $Y^{(t)}$ and $Y^{(t+h)}$, for $h> 0$ small, through a
so-called infinitesimal shift-function $D$. Write $F$ for the
cumulative distribution function and $F^{-1}:\left(0,1\right)\mapsto
{\mathbb R}$ for its generalized inverse
$F^{-1}\left(p\right)=\inf\left\{x:F\left(x\right)\geq p\right\}$.
Then the infinitesimal shift-function $D$ is defined as
\begin{equation} \label{defd}
D\left(y,t;\overline{Z}_t\right)=\left.\frac{\partial}{\partial
  h}\right|_{h=0} \left(F_{Y^{(t+h)}| \overline{Z}_t}^{-1}\circ
F_{Y^{(t)}|\overline{Z}_t}\right)\left(y\right),\end{equation} the
right hand derivative of the quantile-quantile transform which moves
quantiles of the distribution of $Y^{(t)}$ to quantiles of the
distribution of $Y^{(t+h)}$ ($h\geq 0$), given the covariate- and
treatment history until time $t$, $\overline{Z}_t$. In order to define
$D$, no assumptions are necessary about the joint distribution of the
counterfactuals $Y^{(t)}$.

\begin{exa} Survival of AIDS patients.
\label{DAexa2} \emph{\cite{Aids} describe an AIDS clinical trial to
  study the effect of AZT treatment on survival in HIV-positive
  patients. Time $0$ was the time of enrollment in the study. Embedded
  within this trial was an uncontrolled observational
  study of the effect of prophylaxis therapy for PCP on survival. PCP,
  Pneumocystis Carinii Pneumonia, is an opportunistic infection that
  affects HIV-positive patients. \cite{Aids} use continuous-time structural nested models to study
  the effect of PCP prophylaxis therapy on survival of HIV-positive patients. Thus, the outcome of
  interest, $Y$, is the survival time, and the treatment under study is
  prophylaxis for PCP. Although \cite{Aids} estimate the effect of changes in the time the treatment is discontinued, we will consider estimating the effect of changes in the initiation time of the treatment. This conforms better to the clinical practice in HIV/AIDS, where PCP prophylaxis is rarely discontinued, and to the assumption in \cite{Aids} that once PCP prophylaxis is started, it is never stopped. In this example, $Y^{(t)}$ is the counterfactual outcome had PCP prophylaxis treatment been as given in reality until time $t$ and initiated or continued thereafter. We thus define the baseline treatment regime ``$\overline{0}$'' as ``continuously treat with PCP prophylaxis''. In the context of this example, the local rank preservation assumption of \cite{Aids} can be expressed as:
\begin{equation}\label{LRP1}
Y^{(t)} -t = \int_{t}^Y e^{\psi 1_{{\rm no}\;{\rm prophylaxis}\;{\rm at}\; s}} ds.
\end{equation}
Assumption~(\ref{LRP1}) is very strong, because it requires that given the model
parameter $\psi$ and the observed outcome $Y$, all counterfactual outcomes
$Y^{(t)}$ can be calculated from the observed data. The current
article proves that it suffices to assume that
\begin{equation}D_\psi\left(y,t;\overline{Z}_t\right)=
\bigl(1-e^{\psi}\bigr)\; 1_{\left\{{\rm no}\;{\rm prophylaxis}\;{\rm
at}\;t\right\}}.\label{Dexa1}
\end{equation}
We show that under Assumption~(\ref{Dexa1}),
\begin{equation} \label{DAsimAIDS}
Y^{(t)} -t \sim \int_{t}^Y e^{\psi 1_{{\rm no}\;{\rm prophylaxis}\;{\rm at}\; s}} ds
\end{equation}
conditional on $\overline{Z}_t$ and $Y>t$, where $\sim$ means ``has the same distribution as''. Given $\overline{Z}_t$, both $Y^{(t)} -t$ and $\int_{t}^Y e^{\psi
  1_{{\rm no}\;{\rm prophylaxis}\;{\rm at}\; s}} ds$ are random variables,
depending on $Y^{(t)}$ and $Y$, respectively. Assumption~(\ref{Dexa1}) does not impose that $Y^{(t)}-t$ is equal to $\int_{t}^Y e^{\psi 1_{{\rm no}\;{\rm prophylaxis}\;{\rm at}\; s}} ds$, but only that the distribution of these two random variables is the same conditional on $\overline{Z}_t$ and $Y>t$. Thus, under equation~(\ref{Dexa1}), patients who have the exact same observed history over $[0,\tau]$, $\overline{Z}_\tau$ and $Y$, do not necessarily have the same counterfactual outcomes $Y^{(t)}$. This is a substantial relaxation of the assumptions previously adopted in the literature on continuous-time structural nested models. Relaxing assumption~(\ref{LRP1}) is empirically relevant because in clinical practice $Y^{(t)}$ may differ between two patients with the exact same observed history. Suppose for example that two patients with the exact same observed history were both on PCP prophylaxis. If one of the patients got in contact with pneumococcal bacteria (and therefore might have caught PCP without the preventive treatment, PCP prophylaxis), and the other did not get in contact with pneumococcal bacteria (and therefore might not have caught PCP, even without PCP prophylaxis), the outcomes for the two patients without PCP prophylaxis could be different.}

\emph{In equation~(\ref{DAsimAIDS}), the part of the residual survival time, $Y-t$, that is untreated gets multiplied by $e^{\psi}$ to attain the same distribution as $Y^{(t)}-t$ (the residual survival time under ``continuous treatment from $t$ onwards''), conditional on $\overline{Z}_t$ and $Y>t$. Therefore, analogous to accelerated failure
time models (see e.g.~Cox and Oaks, 1984), the multiplication factor $e^{\psi}$ can be interpreted
in a distributional way.}

\emph{Our results do not depend on adopting the particular specification of the infinitesimal shift-function $D$ of equation~(\ref{Dexa1}).
For example, they also apply to an alternative specification of $D$ from \cite{Aids}.
In this alternative specification, the effect of the PCP prophylaxis can depend on the AZT
  treatment the patient received and whether or not the patient had a
  history of PCP prior to the start of PCP prophylaxis. Because the data in \cite{Aids} were from a
  clinical trial for AZT treatment, AZT treatment is described by
  a single variable $R$ indicating the treatment arm the patient was
  randomized to ($R$ equals $1$ or $2$).  Let $P\left(t\right)$ be equal to
  $1$ if the patient had PCP before or at time $t$ \emph{and} before
  prophylaxis treatment started; otherwise $P\left(t\right)$ is equal to
  $0$.  The model described in \cite{Aids}, but adapted to our choice of baseline treatment regime ($\overline{0}$ is continuous treatment with PCP prophylaxis), is
%\begin{equation*}
%Y^{(t)} -t = \int_{t}^Y e^{1_{\left\{{\rm no}\;{\rm prophylaxis}\;{\rm at}\;s\right\}}\left(\psi_1+\psi_2 P\left(s\right) +\psi_3 R\right)} ds.
%\end{equation*}
%For this model, \cite{Aids} stated that
\begin{equation}D_{\psi_1,\psi_2,\psi_3}\left(y,t;\overline{Z}_t\right)
=\bigl(1-e^{\psi_1+\psi_2 P\left(t\right) +\psi_3 R}\bigr)\;
1_{\left\{{\rm no}\;{\rm prophylaxis}\;{\rm at}\;t\right\}}.\label{Dexa2}\end{equation}
This article shows that if equation~(\ref{Dexa2}) holds, then
\begin{equation*}% \label{DAsim2}
Y^{(t)}-t \sim \int_{t}^Y e^{1_{\left\{{\rm no}\;{\rm prophylaxis}\;{\rm
at}\;s\right\}}\left(\psi_1+\psi_2 P\left(s\right) +\psi_3 R\right)} ds
\;\;\;\;\;\;\;given\;\overline{Z}_t,\end{equation*}
for $t<Y$.}
\end{exa}

\begin{exa} \label{DAexa} Effect of Graft versus Host Disease
(GvHD) on time to leukemic relapse. \emph{\cite{Kei} and
    \cite{Keiding} use continuous-time structural nested models to
    study the effect of GvHD on time to leukemic relapse in patients
    who had Bone Marrow Transplantation (BMT). Infection with
    Cytomegalovirus (CMV) is a time-dependent confounder: an
    independent prognostic factor for relapse that both 1. predicts
    the subsequent development of the exposure GvHD and 2. is
    predicted by past exposure GvHD. Write $Y$ for the time until
    leukemic relapse. Assume that $Y$ is observed for every
    patient. In \cite{Kei} and \cite{Keiding}, $Y^{(t)}$ is the outcome
    had the patient been exposed (or not) to GvHD as in reality until time $t$,
    and not exposed afterwards. Based on biological knowledge, \cite{Kei} and
    \cite{Keiding} assume that
\begin{equation}\label{Dexa3} D_\psi\left(y,t;\overline{Z}_t\right)=
\bigl(1-e^{\psi}\bigr)\; 1_{\left\{{\rm GvHD}\;{\rm
at}\;t\right\}}.
\end{equation}
This article
shows that then, for $t<Y$, preventing GvHD from $t$ onwards leads to
\begin{equation} \label{DAsim}
Y^{(t)} -t \sim \int_{t}^Y e^{\psi 1_{\left\{{\rm GvHD}\;{\rm
at}\;s\right\}}} ds
\;\;\;\;{\rm given}\;\overline{Z}_t.
\end{equation}
\cite{Kei} and \cite{Keiding} assume that (\ref{DAsim}) is true even
with $\sim$ replaced by $=$ (although only for $t=0$), hoping that could
be relaxed. This article shows that indeed (\ref{Dexa3}) is sufficient to
estimate the effect of GvHD.}
\end{exa}

\begin{exa} (Incorporating a-priori biological knowledge, following \cite{Enc}). \emph{Again consider survival as the outcome of
interest.  Suppose that it is known that treatment received at time
$t$ only affects survival for patients who would die by time $t+5$ if
they would receive no further treatment. An example would be a setting
in which failure is death from an infectious disease, the treatment is
a preventive antibiotic treatment which is of no benefit unless the
person is already infected and, if death occurs, it always does within
five weeks from the time of initial unrecorded subclinical infection.
In that case, the natural restriction on $D$ is that
\begin{equation*}D\left(y,t;\overline{Z}_t\right)=0 \;\;\;\;\; {\rm if}\; y-t>5.
\end{equation*}}
\end{exa}
As can be seen from these examples, the parameters of a continuous-time structural
nested model are often rates. More biostatistical
examples of models for $D$ can be found in
e.g.\ \cite{smoke,Enc,RGPCP,R92,Tilling}.

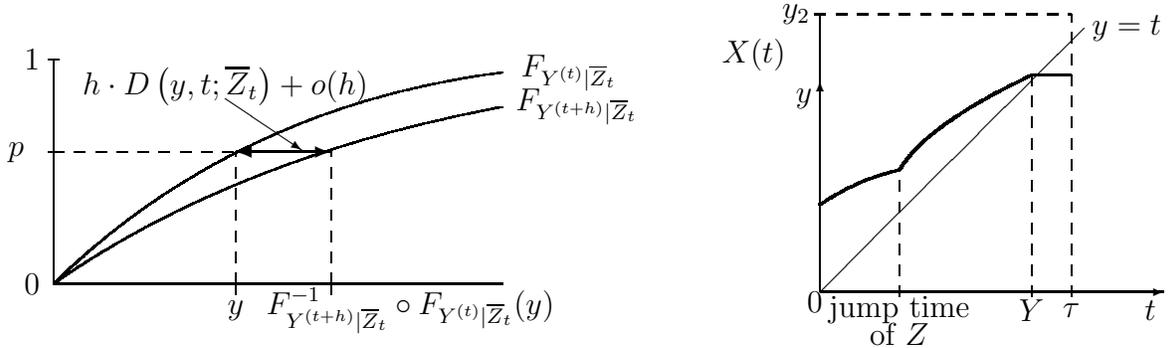
\begin{figure}[htb!]
\begin{picture}(450,112)
\put(30,10){\line(1,0){170}}
\put(30,10){\line(0,1){85}}
\put(22,10){\makebox(0,0){$0$}}
\put(22,95){\makebox(0,0){$1$}}
\put(28,95){\line(1,0){4}}
\qbezier[4000](30,10)(99,58)(200,77)
\qbezier[4000](30,10)(99,77)(200,90)
\put(165,0){\makebox(0,0){$F^{-1}_{Y^{(t+h)}|\overline{Z}_t}\circ F_{Y^{(t)}|\overline{Z}_t}(y)$}}
\put(99,6){\line(0,1){6}}
\put(99,0){\makebox(0,0){$y$}}
\put(135,6){\line(0,1){6}}
\put(228,77){\makebox(0,0){$F_{Y^{(t+h)}|\overline{Z}_t}$}}
\put(225,90){\makebox(0,0){$F_{Y^{(t)}|\overline{Z}_t}$}}
\put(135,10){\line(0,1){4}}
\put(135,18){\line(0,1){4}}
\put(135,26){\line(0,1){4}}
\put(135,34){\line(0,1){4}}
\put(135,42){\line(0,1){4}}
\put(135,50){\line(0,1){4}}
\put(135,58){\line(0,1){2}}
%\put(38,135){\line(1,0){4}}
\put(99,10){\line(0,1){4}}
\put(99,18){\line(0,1){4}}
\put(99,26){\line(0,1){4}}
\put(99,34){\line(0,1){4}}
\put(99,42){\line(0,1){4}}
\put(99,50){\line(0,1){4}}
\put(99,58){\line(0,1){2}}
\put(28,60){\line(1,0){4}}
\put(16,60){\makebox(0,0){$p$}}
\put(95,85){\makebox(0,0){$h\cdot D\left(y,t;\overline{Z}_t\right)+o(h)$}}
%\put(222,135){\line(-4,1){62}}
%\put(222,135){\vector(4,-1){62}}
\put(95,80){\vector(3,-2){29}}
\put(30,60){\line(1,0){4}}
\put(38,60){\line(1,0){4}}
\put(46,60){\line(1,0){4}}
\put(54,60){\line(1,0){4}}
\put(62,60){\line(1,0){4}}
\put(70,60){\line(1,0){4}}
\put(78,60){\line(1,0){4}}
\put(86,60){\line(1,0){4}}
\put(94,60){\line(1,0){4}}
\put(102,60){\line(1,0){4}}
\put(110,60){\line(1,0){4}}
\thicklines
\put(99,60){\vector(1,0){36}}
\put(135,60){\vector(-1,0){36}}
\thinlines
\put(295,97){\makebox(0,0){$X(t)$}}
\put(320,7){\vector(1,0){130}}
\put(320,7){\vector(0,1){80}}
\put(320,7){\line(0,1){105}}
\put(318,3){\makebox(0,0){$0$}}
\put(445,0){\makebox(0,0){$t$}}
\put(311,112){\makebox(0,0){$y_2$}}
\put(314,80){\makebox(0,0){$y$}}
\put(318,112){\line(1,0){4}}
\put(326,112){\line(1,0){4}}
\put(334,112){\line(1,0){4}}
\put(342,112){\line(1,0){4}}
\put(350,112){\line(1,0){4}}
\put(358,112){\line(1,0){4}}
\put(366,112){\line(1,0){4}}
\put(374,112){\line(1,0){4}}
\put(382,112){\line(1,0){4}}
\put(390,112){\line(1,0){4}}
\put(398,112){\line(1,0){4}}
\put(406,112){\line(1,0){4}}
\put(414,112){\line(1,0){1}}
\put(320,7){\line(1,1){100}}
\put(436,107){\makebox(0,0){$y=t$}}
\put(350,5){\line(0,1){4}}
\put(350,13){\line(0,1){4}}
\put(350,21){\line(0,1){4}}
\put(350,29){\line(0,1){4}}
\put(350,37){\line(0,1){4}}
\put(350,45){\line(0,1){4}}
\put(400,5){\line(0,1){4}}
\put(400,13){\line(0,1){4}}
\put(400,21){\line(0,1){4}}
\put(400,29){\line(0,1){4}}
\put(400,37){\line(0,1){4}}
\put(400,45){\line(0,1){4}}
\put(400,53){\line(0,1){4}}
\put(400,61){\line(0,1){4}}
\put(400,69){\line(0,1){4}}
\put(400,77){\line(0,1){4}}
\put(400,85){\line(0,1){4}}
\put(415,5){\line(0,1){4}}
\put(415,13){\line(0,1){4}}
\put(415,21){\line(0,1){4}}
\put(415,29){\line(0,1){4}}
\put(415,37){\line(0,1){4}}
\put(415,45){\line(0,1){4}}
\put(415,53){\line(0,1){4}}
\put(415,61){\line(0,1){4}}
\put(415,69){\line(0,1){4}}
\put(415,77){\line(0,1){4}}
\put(415,85){\line(0,1){4}}
\put(415,93){\line(0,1){4}}
\put(415,101){\line(0,1){4}}
\put(415,109){\line(0,1){3}}
\put(400,0){\makebox(0,0){$Y$}}
\put(415,0){\makebox(0,0){$\tau$}}
\put(400,5){\line(0,1){4}}
\put(350,0){\makebox(0,0){jump time}}
\put(350,-9){\makebox(0,0){of $Z$}}
\put(350,5){\line(0,1){4}}
\thicklines
\put(400,89){\line(1,0){15}}
\qbezier(350,53)(360,70)(400,89)
\qbezier(320,40)(335,50)(350,53)
\thinlines
\end{picture}
\caption{Left: Illustration of the infinitesimal shift-function $D$.\break
  Right: An example of a solution $X(t)$ to the
  differential equation $dX(t)/dt=D\left(X(t),t;\overline{Z}_t\right)$
  with final condition $X\left(\tau\right)=Y$ in case the outcome is
  survival time.}
\label{Dfig}
\end{figure}

$h\cdot D\left(y,t;\overline{Z}_t\right)$ can be interpreted as the
infinitesimal effect on the outcome of the treatment actually given in
the time interval $\left[t,t+h\right)$ (relative to the baseline
  treatment regime). To be more precise, from the definition of $D$, it
  follows that
\begin{equation*}h\cdot D\left(y,t;\overline{Z}_t\right)=
\left(F_{Y^{(t+h)}| \overline{Z}_t}^{-1}\circ
F_{Y^{(t)}|\overline{Z}_t}\right)\left(y\right)-y+o\left(h\right).\end{equation*}
In Figure~\ref{Dfig} (left) this is sketched.

It can be shown that $D\equiv 0$ if and only if treatment does not
affect the outcome of interest, as was conjectured in \cite{Enc}. To
be more precise, \cite{Lok} shows that, for
example, $D\equiv 0$ if and only if for every $h>0$ and $t$,
$Y^{(t+h)}$ has the same distribution as $Y^{(t)}$ given
$\overline{Z}_t$. That is, $D\equiv 0$ if and only if ``at any time
$t$, whatever person characteristics are selected at that time
($\overline{Z}_t$), switching `treatment as given' to `baseline treatment regime' at some fixed time
after $t$ would not change the distribution of the outcome in persons
with these person characteristics''. To prove this one needs the
mimicking result of the current article.

\section{Mimicking counterfactual outcomes}
\label{cmim}

Define $X\left(t\right)$ as the
continuous solution to the differential equation
\begin{equation} \label{Xdef}
dX(t)/dt=D\left(X\left(t\right),t;\overline{Z}_t\right)\end{equation}
with final condition $X\left(\tau\right)=Y$, the observed outcome (see
Figure~\ref{Dfig}, right). Then $X\left(t\right)$ mimics $Y^{(t)}$ in
the sense that it has the same distribution as $Y^{(t)}$, even given
the person's treatment- and covariate history at time $t$, $\overline{Z}_t$.
%This rather surprising result was conjectured in \cite{Enc}
% and is the main result of the current article.
To prove this main result we need the following consistency assumption.
\begin{assu} \label{inst} \emph{(Consistency)}.
$Y^{\left(\tau\right)}$ has the same distribution as $Y$ given
$\overline{Z}_{\tau}$.
\end{assu}
\noindent Notice that because no treatment was given after time $\tau$
and the treatment process is right continuous, there is no difference
in treatment between $Y^{(\tau)}$ and $Y$. Under this consistency
assumption and regularity conditions only, it is proved in
Sections~\ref{XY} and~\ref{XT} that indeed (\ref{Xdef}) has a unique
solution $X$ for every $\omega\in\Omega$, and that this solution
$X(t)$ mimics $Y^{(t)}$ in the sense that it has the same distribution
as $Y^{(t)}$ given $\overline{Z}_t$.
% in the sense that it has the same distribution as $Y^{(t)}$ given
%$\overline{Z}_t$.
%Section~\ref{XY} deals with non-survival outcomes, Section~\ref{XT}
%with survival outcomes. Survival outcomes require a different set of
%assumptions, as will be explained in Section~\ref{XT}.

\begin{exa} \label{DAexai}
  Survival of AIDS patients \emph{ (continuation of Example~\ref{DAexa2}). If equation~(\ref{Dexa1}) holds, then
\begin{equation*}
X(t)=t+ \int_{t}^Y e^{\psi
1_{\left\{{\rm prophylaxis}\;{\rm at}\;s\right\}}} ds
\end{equation*}
for $t<Y$, and $X(t)=Y$ for $t\geq Y$. Alternatively, if
equation~(\ref{Dexa2}) holds, then
\begin{equation*}
X(t)=t+ \int_{t}^Y e^{1_{\left\{{\rm prophylaxis}\;{\rm at}\;s\right\}}\left(\psi_1+\psi_2P(s)+\psi_3 R\right)} ds
\end{equation*}
for $t<Y$, and $X(t)=Y$ for $t\geq Y$.}
\end{exa}

\section{Estimators, tests, and ``no unmeasured confounding''}
\label{nuc}

This section contains a brief summary of \cite{ASarx}, who shows how
the result of the current article leads to testing and estimation.
In addition, Appendix~\ref{simapp} provides an example of estimation
in our simulation study.

The main assumption underlying structural nested models is that all
information the doctors used to make treatment decisions, and which is
predictive of the person's prognosis with respect to the final
outcome, is available for analysis. This assumption of no unmeasured
confounding makes it possible to distinguish between treatment effect
and selection bias; see e.g.\ \cite{Aids}, \cite{Enc}, \cite{SNart} or
\cite{ASarx}.

Assume that the treatment process gives rise to a
counting process $N(t)$. For example, $N(t)$ is the number of
treatment changes until time $t$. The assumption of no unmeasured
confounding is then formalized as
\begin{assu} \emph{(No unmeasured confounding)}. The rate with which $N$ jumps
  given $\overline{Z}_{t-}$ is the same as the rate with which $N(t)$
  jumps given $\overline{Z}_{t-}$ and $\left(Y^{(s)}:s< t\right)$,
\end{assu}
\noindent because given the observed $\overline{Z}_{t-}$, the
(unobserved) prognosis of a person, represented by $Y^{(s)}$ for $s<
t$, should not predict treatment at or after time $t$. If it does, there is no
way to distinguish between the effect of the treatment and the reason
why it is initiated.

Notice that if $X(t)$ mimics $Y^{(t)}$ in the sense that it has the
same distribution as $Y^{(t)}$ given $\overline{Z}_{t-}$, it can be
expected that under no unmeasured confounding, the rate with which
$N(t)$ jumps at time $t$ also does not depend on $X(t)$, given
$\overline{Z}_{t-}$. It can formally be shown that this is indeed true.

First consider how this leads to testing. If treatment does not
affect the outcome of interest, $D\equiv 0$ and thus $X(t)\equiv Y$.
So if treatment does not affect the outcome of interest, changes of
treatment at time $t$ should be independent of $Y$, given
$\overline{Z}_{t-}$. Thus one can test whether treatment affects the
outcome of interest by testing whether, given $\overline{Z}_{t-}$, $Y$
adds to the prediction model for treatment changes.

Also for estimation of the infinitesimal shift-function $D$ we assume
that there is no unmeasured confounding. Suppose that one has a
correctly specified parametric model $D_\psi$ for $D$. Then one can
calculate ``$X_\psi(t)$'', the solution to
\begin{equation} \label{Xpsi}
dX_\psi\left(t\right)/dt=
D_{\psi}\left(X_\psi\left(t\right),t;\overline{Z}_t\right)
\end{equation}
with final condition $X(\tau)=Y$. If $X(t)$ mimics $Y^{(t)}$, then
$X_\psi(t)$ has the same distribution as $Y^{(t)}$ given
$\overline{Z}_t$ for the true $\psi$.  Since $Y^{(t)}$ does not add to
the prediction model for treatment changes given $\overline{Z}_{t-}$,
$\psi$ could then be estimated by picking the $\psi$ for which, given
$\overline{Z}_{t-}$, $X_\psi(t)$ adds the least to the prediction model
for $N$, treatment changes.
%For details and a complete survey we refer to \cite{ASarx}.
This can be proven to lead to the following theorem:
\begin{thm} \label{thec}
  Suppose that the intensity process $\lambda$ is bounded, $Y^{()}$ is
  cadlag, there is no unmeasured confounding and no instantaneous
  treatment effect (with probability $1$, $N()$ and $Y^{()}$ do not
  jump at the same time). Suppose also that for every
  $t\in\left[0,\tau\right]$, $X(t)$ has the same distribution as
  $Y^{(t)}$ given $\overline{Z}_t$.  Then
\begin{equation*}E\int_0^\tau  h_t\left(X\left(t\right),\overline{Z}_{t-}\right)
 \left(dN(t)-\lambda\left(t\right)dt\right)=0\end{equation*}
for each $h_t$ satisfying a regularity restriction. Thus if
$D_\psi$ and $\lambda_\theta$ are correctly specified parametric
models for $D$ and $\lambda$, respectively,
\begin{equation*}P_n\int_0^\tau  h_t\left(X_{\psi}\left(t\right),\overline{Z}_{t-}\right)
 \left(dN(t)-\lambda_{\theta}\left(t\right)dt\right)=0,\end{equation*}
with $P_n$ the empirical measure $P_n \,X=1/n \sum_{i=1}^n X_i$, is an
unbiased estimating equation for $\left(\theta_0,\psi_{0}\right)$, for
each $h_t$ satisfying a regularity restriction. $h_t$ here is allowed to
depend on $\psi$ and $\theta$, as long as it satisfies
the regularity restriction for $\left(\theta_0,\psi_{0}\right)$.
\end{thm}
In fact, these estimating equations are often martingales at the true
parameter. From e.g.~\cite{Vaart}, Theorem~\ref{thec} implies that
the resulting estimators are, under regularity conditions, consistent
and asymptotically normal.

%A subset of these estimators can be calculated with standard software,
%as was conjectured in \cite{Enc} and proved in
%\cite{plss} using a partial likelihood approach.

\section{Outline of the proof}
\label{XYin}

Throughout the proof this article uses fixed versions of
$F_{Y^{(t+h)}|\overline{Z}_t}$ satisfying all regularity conditions of
Section~\ref{result}. Section~\ref{Dsec} shows existence of $D$.
It also derives a different expression for $D$, which is often used in
the rest of the proof. Section~\ref{Dexu} shows existence and
uniqueness of solutions $X(t)$ to the differential equation with $D$, equation~(\ref{Xdef}),
with final condition $X\left(\tau\right)=Y$.

The proof that this $X(t)$ mimics $Y^{(t)}$ is based on
discretization. Section~\ref{det} therefore considers the situation
where the treatment- and covariate process $Z$ can be fully described
by its values at finitely many fixed times
$0<\tau_1<\tau_2<\ldots<\tau_K$ and $\tau$. In fact this is the
discrete-time situation studied in \cite{SNart}, but instead of using
the shift-function $\gamma$ described there as a model this article uses the
infinitesimal shift-function $D$. Proposition~\ref{detth} in
Section~\ref{det} states that in this discrete-time setting with $D$
instead of $\gamma$, $X\left(t\right)$ mimics $Y^{(t)}$ in the sense
that it has the same distribution as $Y^{(t)}$ given the discrete-time
$\overline{Z}_t$, under a regularity condition and Consistency
Assumption~\ref{inst}.  The proof of Proposition~\ref{detth} is
relatively easy, because in this discrete-time setting the continuous
solution to the differential equation can be written down explicitly,
in terms of conditional distribution functions.

Sections~\ref{discsch}--\ref{Concl} consider the situation where the
probability that $Z$ jumps \emph{at} $t$ equals zero for all $t$.  We
prove that also in this case, $X\left(t\right)$ mimics $Y^{(t)}$,
under the conditions of Section~\ref{result}. First,
Section~\ref{discsch} prepares the discretization by constructing a
series $\overline{Z}^{(n)}$, containing more and more information on
the covariate- and treatment history $\overline{Z}$ as $n$
increases. $\overline{Z}^{(n)}$ depends deterministically on
$\overline{Z}$, so that no extra randomness is necessary to construct
$\overline{Z}^{(n)}$. The discretization does not change $Y^{(t)}$;
just the information on the treatment- and covariate process
considered is reduced.
% The discretized process
$\overline{Z}^{(n)}$ is a covariate- and treatment history as
considered in Section~\ref{det}. Therefore, $D^{(n)}$ can be defined as
\begin{equation} D^{(n)}\Bigl(y,t,\overline{Z}_t^{(n)}\Bigr)
=\left.\frac{\partial}{\partial h}\right|_{h=0}
\left(F_{Y^{(t+h)}|\overline{Z}_t^{(n)}}^{-1}\circ
F_{Y^{(t)}|\overline{Z}_t^{(n)}}\right)\left(y\right)
\label{Dn}
\end{equation}
and we define $X^{(n)}$ as the continuous solution to the
differential equation
\begin{equation}
\frac{d}{dt}X^{\left(n\right)}(t)
=D^{(n)}\Bigl(X^{(n)}(t),t;\overline{Z}^{(n)}_t\Bigr)
\label{Xn}
\end{equation}
with final condition $X^{(n)}\left(\tau\right)=Y$. Section~\ref{Dnsec}
shows existence of $D^{(n)}$ and provides two expressions for
$D^{(n)}$. Section~\ref{Sdisc} shows that the conditions of the
discrete-time result are satisfied for the discretized situation, so
that Proposition~\ref{detth} guarantees that there exists a continuous
solution $X^{(n)}(t)$ to the differential equation (\ref{Xn}), with
final condition $X^{(n)}\left(\tau\right)=Y$ and with the same
distribution as $Y^{(t)}$ given $\overline{Z}_t^{(n)}$.

%To prove that $X^{(n)}(t)$ converges almost surely to $X(t)$ as
%$n\rightarrow \infty$, we bound the difference between $X^{(n)}(t)$
%and $X(t)$ in terms of $D$ and $D^{(n)}$, using a result from the
%theory of differential equations, Theorem~\ref{difb} in
%Appendix~\ref{\difeqap}. After bounding the difference between
%$X^{(n)}(t)$ and $X(t)$ in terms of $D$ and $D^{(n)}$, we show that
%$D^{(n)}\bigl(y,t;\overline{Z}_t^{(n)}\bigr)$ converges almost surely
%to $D\left(y,t;\overline{Z}_t\right)$ as $n\rightarrow \infty$. We
%conclude this section by proving that this implies that $X^{(n)}(t)$
%converges almost surely to $X(t)$.

Sections~\ref{XXnaf}--\ref{XnX} then prove that
$X^{(n)}\left(t\right)$ converges almost surely to $X(t)$ as $n$ tends
to infinity, using a result from differential equation theory which
bounds the difference between solutions to differential equations. The
proof is concluded in Section~\ref{Concl}, which shows that $X(t)$
has the same distribution as $Y^{(t)}$ given $\overline{Z}_t$ because
$X^{(n)}(t)$ has the same distribution as $Y^{(t)}$ given
$\overline{Z}_t^{(n)}$ and $X^{(n)}(t)$ converges almost surely to
$X(t)$.

%The proof can easily be adapted to the situation where the smoothness
%conditions only hold for $h$ in a neighbourhood $\left[0,
%\delta\right]$ of $0$, for some fixed $\delta>0$, by starting the
%discretization with $n$ satisfying e.g.\ $\tau/(2^n)<\delta/2$ instead
%of with $n=1$.
Section~\ref{Disc2} indicates how the proof can be adapted to
include situations where the probability that $Z$ jumps at time $t$ is
zero except for at finitely many times $t$.

\section{Proof of main result}
\label{XY}

\subsection{Introduction}

The purpose of the current article is to prove that $X\left(t\right)$
mimics $Y^{(t)}$. This result is proved in this section for
non-survival outcomes. Section~\ref{result} states the assumptions and
the precise statement of the result, and
Sections~\ref{Dsec}--\ref{Disc2} provide the proof.

\subsection{Mimicking counterfactual non-survival outcomes: assumptions
and result}
\label{result}

This section provides precise conditions under which $X(t)$ mimics
$Y^{(t)}$. First, consider the definition of $D$,
equation~(\ref{defd}). Notice that $D$ involves an uncountable number
of distribution functions $F_{Y^{(t+h)}| \overline{Z}_t}$. In many
cases conditioning on $\overline{Z}_t$ means conditioning on a
null-event, so that these conditional distributions are not unique.
Every single conditional distribution is almost surely unique (see
Web-Appendix~\ref{Cond}), but because an uncountable number of
them is used ($t$ and $h$ are continuous) this is not sufficient for overall
almost sure uniqueness.  Therefore the regularity conditions below
should be read as: there exists a collection of conditional
distribution functions $F_{Y^{(t+h)}| \overline{Z}_t}$ such that all
these regularity conditions are satisfied. These versions of
$F_{Y^{(t+h)}| \overline{Z}_t}$ are chosen in the definition of $D$ as
well as everywhere else in this article. We only consider $h\geq 0$, so
the derivative with respect to $h$ at $h=0$ is always the right hand
derivative.
%As we will show in Section~\ref{XYin}, all assumptions can
%be relaxed to $h$ in a neighbourhood $\left[0,\delta\right]$ of $0$,
%provided that this neighbourhood does not depend on $\overline{Z}$.

With the support of a random variable $X$ this article means those $x$
such that for every open set $U_x$ containing $x$, $P\left(X\in
  U_x\right)>0$. Let $y_1$ and $y_2$ be the lower- and upper limit of
the support of the outcome of interest $Y$. In this article, these are
assumed to be finite, and moreover it is assumed that
\begin{assu}\label{sup}\emph{(support)}.
\begin{enumerate}[a)]
\item All $F_{Y^{(t+h)}|\overline{Z}_t}$ for all $t\geq 0$ and for $h\geq 0$
  have the same bounded support $\left[y_1,y_2\right]$.
\item All $F_{Y^{(t+h)}|\overline{Z}_t}\left(y\right)$ for all $t\geq 0$ and for $h\geq 0$
have a continuous non-zero density
$f_{Y^{(t+h)}|\overline{Z}_t}\left(y\right)$ on
$y\in\left[y_1,y_2\right]$.
\item There exists an $\varepsilon>0$ such that
$f_{Y^{(t)}|\overline{Z}_t}\left(y\right)\geq \varepsilon$ for all
$y\in\left[y_1,y_2\right]$, $\omega\in\Omega$ and $t\in\left[0,\tau\right]$.
\end{enumerate}
\end{assu}
The support condition may be restrictive for certain applications.
Nevertheless, most real-life situations can be approximated this way,
since $y_1$ and $y_2$ are can have arbitrary (finite) values and
$\varepsilon>0$ can be vary small. Although the support condition
may well be stronger than necessary, it simplifies the analysis
considerably and, for that reason, it is adopted here.

The remaining regularity conditions are smoothness conditions. They
allow for non-smooth-ness where the covariate- and treatment process
$Z()$ jumps. This is important since if the covariate- and treatment
process $Z()$ jumps this can lead to a different prognosis for the
person and thus to non-smoothness of the functions concerned.
\begin{assu} \label{cont1} \emph{(continuous derivatives)}.
For $\omega\in\Omega$ fixed,
\begin{enumerate}[a)]
\item $F_{Y^{(t+h)}|\overline{Z}_t}\left(y\right)$ is $C^1$ in
$\left(h,y\right)$ for $y\in\left[y_1,y_2\right]$ and $h\geq 0$.
%is differentiable
%with respect to $h$ and $y$ with derivatives which are continuous in
%$\left(h,y\right)$ on for $y\in\left[y_1,y_2\right]$.
\item If $Z$ does not jump in $\left(t_1,t_2\right)$ then both
$\left.\frac{\partial}{\partial h}\right|_{h=0} F_{Y^{(t+h)}|\overline{Z}_t}\left(y\right)$ and
$\frac{\partial}{\partial y}F_{Y^{(t)}|\overline{Z}_t}\left(y\right)$
are continuous in $\left(y,t\right)$ on
$\left[y_1,y_2\right]\times\left[t_1,t_2\right)$
and can be continuously extended to
$\left[y_1,y_2\right]\times\left[t_1,t_2\right]$.
\end{enumerate}
\end{assu}

Structural nested models in continuous time are meant to estimate the
effect of a continuous treatment, for which the effect of a small
duration is small. Then, Assumption~\ref{bdd+} is a regularity
condition:
\begin{assu} \label{bdd+} \emph{(bounded derivatives)}.
\begin{enumerate}[a)]
\item There exists a constant $C_1$ such that for all $\omega\in\Omega$,
 $t$, $h\geq 0$ and $y\in [y_1,y_2]$,
\begin{equation*} \frac{\partial}{\partial y}
F_{Y^{(t+h)}|\overline{Z}_t}\left(y\right)\leq C_1.
\end{equation*}
\item There exists a constant $C_2$ such that for all $\omega\in\Omega$,
$t$, $h\geq 0$ and $y\in [y_1,y_2]$,
\begin{equation*} \Bigl|\frac{\partial}{\partial h}
F_{Y^{(t+h)}|\overline{Z}_t}\left(y\right)\Bigr|\leq C_2.
\end{equation*}
\end{enumerate}
\end{assu}

\begin{assu} \label{Lipc}
\emph{(Lipschitz continuity)}.
\begin{enumerate}[a)]
\item There exists a constant $L_1$ such that for all $\omega\in\Omega$ and $t$
and $y,z\in\left[y_1,y_2\right]$,
\begin{equation*}
\Bigl|\frac{\partial}{\partial y}F_{Y^{(t)}|\overline{Z}_t}\left(y\right)-
\frac{\partial}{\partial z}F_{Y^{(t)}|\overline{Z}_t}\left(z\right)\Bigr| \leq
L_1\left|y-z\right|. \end{equation*}
\item There exists a constant $L_2$ such that for all $\omega\in\Omega$ and
$t$ and $y,z\in\left[y_1,y_2\right]$,
\begin{equation*}
\left|\left.\frac{\partial}{\partial h}\right|_{h=0} F_{Y^{(t+h)}|\overline{Z}_t}\left(y\right)-
\left.\frac{\partial}{\partial h}\right|_{h=0}
F_{Y^{(t+h)}|\overline{Z}_t}\left(z\right)\right| \leq
L_2\left|y-z\right|.\end{equation*}
\end{enumerate}
\end{assu}
%Furthermore, as we will see in Section~\ref{XT}, we could have
%restricted all these assumptions to $h\geq 0$. The derivative with
%respect to $h$ at $h=0$ should then be read as the right hand
%derivative.

The main theorem of this article is
\begin{thm} \label{thm1}\emph{(mimicking counterfactual outcomes)}.
  Suppose that Regularity Conditions \ref{sup}--\ref{Lipc} are
  satisfied. Then $D\left(y,t;\overline{Z}_t\right)$ exists.
  Furthermore for every $\omega\in\Omega$ there exists exactly one
  continuous solution $X(t)$ to
  $dX(t)/dt=D\left(X(t),t;\overline{Z}_t\right)$ with final condition
  $X\left(\tau\right)=Y$. If also Consistency Assumption~\ref{inst} is
  satisfied, then this $X(t)$ has the same distribution as $Y^{(t)}$
  given $\overline{Z}_t$ for all $t\in\left[0,\tau\right]$.
\end{thm}
%and $P\left(t \;{\rm is} \;{\rm a}\; {\rm
%jump}\;{\rm time}\;{\rm of}\;\overline{Z}\right)= 0$ for every $t$,
%Theorem~\ref{thm1} extends to the case of finitely many specific times
%$t$ with $P\left(t \;{\rm is} \;{\rm a}\; {\rm jump}\;{\rm time}\;{\rm
%    of}\;\overline{Z}\right)> 0$; see Section~\ref{Disc2}.

\subsubsection{Simpler regularity conditions}\label{simpler}

I state some more restrictive but simpler conditions implying all the
conditions in Section~\ref{result}:
\begin{assu} \emph{(regularity condition)}.
\begin{itemize}
\item \emph{(support)}.
\begin{enumerate}[a)]
\item There exist finite numbers $y_1$ and $y_2$
such that all $F_{Y^{(t+h)}|\overline{Z}_t}$ have the same bounded
support $\left[y_1,y_2\right]$.
\item All $F_{Y^{(t+h)}|\overline{Z}_t}\left(y\right)$ have a
continuous non-zero density
$f_{Y^{(t+h)}|\overline{Z}_t}\left(y\right)$ on
$y\in\left[y_1,y_2\right]$.
\item There exists an $\varepsilon>0$ such that
$f_{Y^{(t)}|\overline{Z}_t}\left(y\right)\geq \varepsilon$ for all
$y\in\left[y_1,y_2\right]$, $\omega\in\Omega$ and $t\in\left[0,\tau\right]$.
\end{enumerate}
\item \emph{(smoothness)}. For every $\omega\in\Omega$
\begin{enumerate}[a)]
\item $\left(y,t,h\right)\rightarrow F_{Y^{(t+h)}|\overline{Z}_t}\left(y\right)$ is
differentiable with respect to $t$, $y$ and $h$ with continuous
derivatives on $\left[y_1,y_2\right]\times\left[t_1,t_2\right)\times
\left[0,\infty\right)$ if $Z$ does not jump in $\left(t_1,t_2\right)$, with a
continuous extension to $\left[y_1,y_2\right]\times\left[t_1,t_2\right]\times
\left[0,\infty\right)$.
\item The derivatives of $F_{Y^{(t+h)}|\overline{Z}_t}\left(y\right)$
with respect to $y$ and $h$ are bounded by constants $C_1$ and $C_2$,
respectively.
\item $\frac{\partial}{\partial y} F_{Y^{(t)}|\overline{Z}_t}(y)$ and
$\left.\frac{\partial}{\partial h}\right|_{h=0}
F_{Y^{(t+h)}|\overline{Z}_t}(y)$ have derivatives with respect to $y$
which are bounded by constants $L_1$ and $L_2$, respectively.
\end{enumerate}
\end{itemize}
\end{assu}

%\subsubsection{Outline}
%\label{disco}

%\subsection[Continuous time]{Mimicking counterfactual outcomes:
%continuous time}
%\label{Disc}

\subsection{Existence of and a different expression for $D$}
\label{Dsec}

The lemma below can be used to prove existence of $D$ and to find a
useful formula for $D$ (and later two useful formulas for $D^{(n)}$ in
Section~\ref{det}):

%\begin{lem} \label{quot}
%Suppose that $F_h$ is a family of non-decreasing
%functions and $F_h\left(y\right)$ is
%differentiable with respect to $y$ and $h$ in a neighbourhood
%$U_{h_0,y_0}$ of
%$\left(h_0,y_0\right)$, with derivatives which are continuous in
%$\left(h,y\right)$. If also $F_{h_0}'\left(y_0\right)$ is non-zero
%then $F_h$ is invertible in a neighbourhood of
%$\left(h_0,y_0\right)$. Moreover $\left(\frac{\partial}{\partial h}
%F_h^{-1}\right)\left(F_h\left(y\right)\right)$ exists and satisfies
%\begin{equation*}
% \frac{\partial}{\partial h} F_h\left(y\right)+
%F_h'\left(y\right)\cdot
%\Bigl(\frac{\partial}{\partial h}
%F^{-1}_h\Bigr)\left(F_h\left(y\right)\right)=0\end{equation*}
%in a neighbourhood of $\left(h_0,y_0\right)$.
%\end{lem}
%
%\noindent{\bf Proof.}
%
%
%$\;$
%
%\noindent We want to use Lemma~\ref{quot} to prove that $D$ exists and
%to find a nice expression for $D$, but $D$ was defined as a right hand
%derivative with respect to $h$. Hence we will use the corollary below,
%which requires differentiability etcetera for $h\geq 0$ only.

\begin{lem} \label{quotc}
Suppose that $F_h$ is a family of non-decreasing functions. Suppose
that there exists a neighbourhood $U_{0,y_0}$ of $\left(0,y_0\right)$
so that $F_h\left(y\right)$ is differentiable with respect to $y$ and
$h$ on $U_{0,y_0}\cap \left\{h\geq 0\right\}$. For $h=0$, the right hand derivative is meant. Suppose furthermore that these derivatives
are continuous in $\left(h,y\right)$. If also
$F_{0}'\left(y_0\right)$ is non-zero then there exists a neighbourhood
$V_{0,y_0}$ of $\left(0,y_0\right)$ such that on the restriction of
this neighbourhood to $h\geq 0$, $F_h$ is invertible. Moreover,
$\left(\frac{\partial}{\partial h}
F_h^{-1}\right)\left(F_h\left(y\right)\right)$ exists and satisfies
\begin{equation*}
 \frac{\partial}{\partial h} F_h\left(y\right)+
F_h'\left(y\right)\cdot
\left(\frac{\partial}{\partial h}
F^{-1}_h\right)\left(F_h\left(y\right)\right)=0.\end{equation*}
\end{lem}

\noindent{\bf Proof.}
Define an extension of $F$ to
\begin{equation*}
\tilde{U}_{0,y_0}=\left\{\left(y,h\right): h\geq 0 \; {\rm and}\;
\left(y,h\right)\in
U_{0,y_0}\right\} \cup \left\{\left(y,h\right): h< 0 \;{\rm and}\;
\left(y,-h\right)\in U_{0,y_0}\right\},\end{equation*}
an open neighbourhood of $\left(0,y_0\right)$, in the following way:
\begin{equation*}
\tilde{F}_h\left(y\right)=\left\{\begin{array}{ll}F_h\left(y\right)&{\rm
if}\;h\geq 0\\
2F_0\left(y\right)- F_{-h}\left(y\right)&{\rm if}\;
h<0.\end{array}\right. \end{equation*}
Define $\phi:U_{h_0,y_0}\rightarrow {\mathbb R}^2$ as
$\phi\left(h,y\right)=\bigl(h,\tilde{F}_h\left(y\right)\bigr)$. The result
follows from the Local Inverse Function Theorem and direct
calculation, after noticing that $D\bigl(\phi\circ\phi^{-1}\bigr)$ is
the identity mapping; see Web-Appendix~\ref{LIFT} for details. \hfill $\Box$

$ $

Because of Assumptions~\ref{cont1}a and~\ref{sup}c,
Lemma~\ref{quotc} can be applied to
$F_h\left(y\right)=F_{Y^{(t+h)}|\overline{Z}_t}\left(y\right)$ with
$y_0=y$. Thus $D$ as defined in equation~(\ref{defd}) exists and
\begin{equation} D\left(y,t;\overline{Z}_t\right)
%&=&\left.\frac{\partial}{\partial h}\right|_{h=0}\left(
%F^{-1}_{Y^{(t+h)}|\overline{Z}_t}\circ
%F_{Y^{(t)}|\overline{Z}_t}\left(y\right)\right)\nonumber\\
%&=&\left(\left.\frac{\partial}{\partial h}\right|_{h=0}
%F^{-1}_{Y^{(t+h)}|\overline{Z}_t}\right)
%\left(F_{Y^{(t)}|\overline{Z}_t}\left(y\right)\right)\nonumber\\
=-\frac{\left.\frac{\partial}{\partial h}\right|_{h=0}
F_{Y^{(t+h)}|\overline{Z}_t}\left(y\right)}{\frac{\partial}{\partial
y}F_{Y^{(t)}|\overline{Z}_t}\left(y\right)}.
\label{Deq}
\end{equation}

\subsection{Existence and uniqueness of $X(t)$} \label{Dexu}

This section shows that the differential equation
$dX(t)/dt=D\left(X(t),t;\overline{Z}_t\right)$ with final condition
$X\left(\tau\right)=Y$ has a unique continuous solution. Fix $\omega$
for the rest of Section~\ref{Dexu}. Since $D$ may be discontinuous at
the jump times of the covariate- and treatment process $Z$, we
consider the intervals between jumps of $Z$ separately.  It suffices
to prove existence and uniqueness of $X(t)$ with final condition on
any interval between jumps of $Z$, because with probability one $Z$
only jumps finitely many times.

Hence suppose that $Z$ does not jump in $\left(t_1,t_2\right)$ and
that $t_1$ is either a jump time of $Z$ or $0$ and that $t_2$ is
either a jump time of $Z$ or $\tau$. From equation~(\ref{Deq}) I
conclude that $D\left(y,t;\overline{Z}_t\right)$ is continuous on
$\left[y_1,y_2\right]\times\left[t_1,t_2\right)$ because of
  Assumptions~\ref{cont1}b and~\ref{sup}c. The differential equation
  has a final condition at the upper end of the interval
  $\left[t_1,t_2\right)$. Therefore we define $\tilde{D}$ on
    $\left[y_1,y_2\right]\times\left[t_1,t_2\right]$ as
\begin{equation*} \tilde{D}\left(y,t\right)=\left\{
\begin{array}{ll}
D\left(y,t;\overline{Z}_t\right)&{\rm if} \;t\in \left[t_1,t_2\right)\\
\lim_{t\uparrow t_2}D\left(y,t;\overline{Z}_t\right)&{\rm if} \; t=t_2.
\end{array}\right.
\end{equation*}
This limit exists because of Assumption~\ref{sup}c and the
extension-assumption in Assumption~\ref{cont1}b. It makes $\tilde{D}$
continuous on $\left[y_1,y_2\right]\times\left[t_1,t_2\right]$. When
calculating the continuous solution to
$dX(t)/dt=D\left(X(t),t;\overline{Z}_t\right)$ on $\left[t_1,t_2\right]$,
one means to use $\tilde{D}$ on $\left[t_1,t_2\right]$ if $D$ jumps
at $t_2$.

To prove existence and uniqueness of $X$ on $\left[t_1,t_2\right]$, we
apply Theorem~\ref{difbc2} to the differential equation with
$\tilde{D}$. We check the conditions of Theorem~\ref{difbc2} for
$\tilde{D}$. Continuity of $\tilde{D}$ was shown in the previous
paragraph. $F^{-1}_{Y^{(t+h)}|\overline{Z}_t}\circ
F_{Y^{(t)}|\overline{Z}_t}\left(y_1\right)=y_1$ for all $h$ because of
Assumption~\ref{sup}a and~b, so that
$D\left(y_1,t;\overline{Z}_t\right)=0$. Similarly,
$D\left(y_2,t;\overline{Z}_t\right)=0$. To show that
equation~(\ref{Lipvw}) holds, notice that global Lipschitz continuity
of $\tilde{D}$ in $y$ on
$\left[y_1,y_2\right]\times\left[t_1,t_2\right)$ with Lipschitz
constant $C=L_2/\varepsilon+L_1 C_2/\varepsilon^2$ follows from
equation~(\ref{Deq}), since the numerator is bounded by $C_2$ and is
Lipschitz with Lipschiz constant $L_2$ and also the denominator is
Lipschitz with Lipschitz constant $L_1$ and bounded away from $0$ by
$\epsilon$ (Assumptions~\ref{Lipc}, \ref{bdd+}b and~\ref{sup}c; see
Web-Appendix~\ref{Lips}). This same constant works on
$\left[y_1,y_2\right]\times\left[t_1,t_2\right]$ by continuity.  By
Theorem~\ref{difbc2}, the differential equation~(\ref{Xdef}) with
$\tilde{D}$ has a unique solution, and this solution stays in
$\left[y_1,y_2\right]$.

\subsection[Discrete time]{Mimicking counterfactual outcomes: discrete time}
\label{det}

This section considers the situation where $\overline{Z}$, the available
information on the treatment- and covariate process, can be fully
described by its values at finitely many fixed time points
$0=\tau_0<\tau_1<\tau_2<\ldots<\tau_K<\tau_{K+1}=\tau$. At these time
points, $Z(t)$ may jump with probability greater than zero. We prove that
in this situation, $X(t)$ mimics $Y^{(t)}$.

We assume that there exist conditional distribution functions (\cite{Bau}, \cite{Pol2})
$F_{Y^{(t)}| \overline{Z}_{\tau_k}}$ satisfying the following
regularity condition:
\begin{assu} \label{cond} \emph{(smoothness)}.
Suppose that for $k=0,\ldots,K$ and $t\in\left[\tau_k,\tau_{k+1}\right]$
there exist conditional distribution functions
$F_{Y^{(t)}|\overline{Z}_{\tau_{k}}}$ such that
\begin{enumerate}[a)]
\item{For all $t\in\left[\tau_k,\tau_{k+1}\right]$,
$F_{Y^{(t)}|\overline{Z}_{\tau_{k}}}\left(y\right)$ is continuous in
$y$.}
\item{For all $t\in\left[\tau_k,\tau_{k+1}\right]$, the support of
$F_{Y^{(t)}|\overline{Z}_{\tau_{k}}}\left(y\right)$ is an interval.}
\item{For $x\in\left[0,1\right]$ fixed,
$F^{-1}_{Y^{(t)}|\overline{Z}_{\tau_k}}\left(x\right)$ is differentiable
with respect to $t$ on
$\left[\tau_k,\tau_{k+1}\right]$.}
\end{enumerate}
\end{assu}
Throughout Section~\ref{det}, fixed versions of
$F_{Y^{(t)}|\overline{Z}_{\tau_k}}\left(y\right)$ are used satisfying
Assumption~\ref{cond}. Since $\overline{Z}_t$ contains the same
information as $\overline{Z}_{\tau_k}$ for
$t\in\left[\tau_k,\tau_{k+1}\right)$, we can and will choose the same
versions when conditioning on $\overline{Z}_t$.

\begin{prop} \label{detth} \emph{(mimicking counterfactual outcomes in discrete time).} Suppose that the treatment- and
covariate process $Z$ can be fully described by its values at finitely
many fixed points
$0=\tau_0<\tau_1<\tau_2<\ldots<\tau_K<\tau_{K+1}=\tau$, and suppose
also that Smoothness Assumption~\ref{cond} is satisfied. Then
$D\left(y,t;\overline{Z}_t\right)$ as defined in equation~(\ref{defd})
exists for all $t$. Furthermore if also Assumption~\ref{inst}
(consistency) is satisfied, then there exists a continuous solution
$X(t)$ to $dX(t)/dt=D\left(X(t),t;\overline{Z}_t\right)$ with final
condition $X\left(\tau\right)=Y$ for which $X\left(t\right)$ has the
same distribution as $Y^{(t)}$ given $\overline{Z}_t$.
\end{prop}

\noindent{\bf Proof.} For $t\in\left[\tau_k,\tau_{k+1}\right)$,
$D\left(y,t;\overline{Z}_t\right)
%&=&\left.\frac{\partial}{\partial
%h}\right|_{h=0} \left(F_{Y^{(t+h)}| \overline{Z}_t}^{-1}\circ
%F_{Y^{(t)}|\overline{Z}_t}\right)\left(y\right)\\
=\left.\frac{\partial}{\partial h}\right|_{h=0} \bigl(F_{Y^{(t+h)}|
\overline{Z}_{\tau_k}}^{-1}\circ
F_{Y^{(t)}|\overline{Z}_{\tau_k}}\bigr)\left(y\right)$,
so existence of $D\left(y,t;\overline{Z}_t\right)$ on each interval
$\left[\tau_k,\tau_{k+1}\right)$ follows from Assumption~\ref{cond}c.

Next, define $\tilde{X}$ as follows. $\tilde{X}\left(\tau\right)=Y$,
and for $t\in\left[\tau_k,\tau_{k+1}\right)$ ($k=0,\ldots,K-1$),
\begin{eqnarray*} \tilde{X}(t)&=&
F^{-1}_{Y^{(t)}|\overline{Z}_{\tau_k}}\circ
F_{Y^{\left(\tau_{k+1}\right)}|\overline{Z}_{\tau_{k}}} \circ \ldots
\circ
F^{-1}_{Y^{\left(\tau_{K-1}\right)}|\overline{Z}_{\tau_{K-1}}}\circ
F_{Y^{\left(\tau_K\right)}|\overline{Z}_{\tau_{K-1}}}\circ\\ &&\hspace{5.2cm}
\circ F^{-1}_{Y^{\left(\tau_K\right)}|\overline{Z}_{\tau_{K}}}\circ
F_{Y^{\left(\tau\right)}|\overline{Z}_{\tau_{K}}}
\left(Y\right).%\label{x}
\end{eqnarray*}
$\tilde{X}(t)$ is well-defined because of Assumption~\ref{cond}a and
b.  First we show that $\tilde{X}=X$: it is a continuous solution to
$\tilde{X}'(t)=D\left(\tilde{X}(t),t;\overline{Z}_t\right)$ with
$\tilde{X}\left(\tau\right)=Y$. Next we show that
$\tilde{X}\left(t\right)$ has the same distribution as $Y^{(t)}$ given
$\overline{Z}_t$.

Continuity of $\tilde{X}$ on $\left[\tau_k,\tau_{k+1}\right)$ is clear
from Assumption~\ref{cond}c. Moreover,
\begin{eqnarray*}\lim_{t\uparrow \tau_{k+1}} \tilde{X}\left(t\right)
&=&\lim_{t\uparrow \tau_{k+1}}F_{Y^{(t)}|\overline{Z}_{\tau_k}}^{-1}\circ
F_{Y^{\left(\tau_{k+1}\right)}|\overline{Z}_{\tau_k}}
\left(\tilde{X}\left(\tau_{k+1}\right)\right)\\
&=&F_{Y^{\left(\tau_{k+1}\right)}|\overline{Z}_{\tau_k}}^{-1}\circ
F_{Y^{\left(\tau_{k+1}\right)}|\overline{Z}_{\tau_k}}
\left(\tilde{X}\left(\tau_{k+1}\right)\right)
\end{eqnarray*}
because of Assumption~\ref{cond}c, which is equal to
$\tilde{X}\left(\tau_{k+1}\right)$ because of Assumption~\ref{cond}b. Thus,
$\tilde{X}(t)$ is also continuous from the left at $t=\tau_{k+1}$.
For $t\in\left[\tau_k,\tau_{k+1}\right)$, $\tilde{X}$ satisfies the
differential equation:
\begin{eqnarray*}
\tilde{X}'\left(t\right)
%&=& \left.\frac{\partial}{\partial h}\right|_{h=0} X(t+h)\\
&=&\left.\frac{\partial}{\partial h}\right|_{h=0}
F^{-1}_{Y^{(t+h)}|\overline{Z}_{\tau_k}}\circ
F_{Y^{\left(\tau_{k+1}\right)}|\overline{Z}_{\tau_{k}}}
\left(\tilde{X}\left(\tau_{k+1}\right)\right)\\
&=& \left(\left.\frac{\partial}{\partial h}\right|_{h=0}
F^{-1}_{Y^{(t+h)}|\overline{Z}_t} \circ F_{Y^{(t)}|\overline{Z}_t}\right)
\circ F^{-1}_{Y^{(t)}|\overline{Z}_{\tau_k}}\circ
F_{Y^{\left(\tau_{k+1}\right)}|\overline{Z}_{\tau_{k}}}
\left(\tilde{X}\left(\tau_{k+1}\right)\right)\\
&=&D\left(\tilde{X}(t),t;\overline{Z}_t\right),
\end{eqnarray*}
where in the second line it is used that conditioning on $\overline{Z}_t$
is the same as conditioning on $\overline{Z}_{\tau_k}$, so that
$F_{Y^{(t)}|\overline{Z}_t} \circ
F^{-1}_{Y^{(t)}|\overline{Z}_{\tau_k}}$ is the identity because of
Assumption~\ref{cond}a and b. Thus indeed $\tilde{X}$ is a continuous
solution to
$\tilde{X}'(t)=D\left(\tilde{X}(t),t;\overline{Z}_t\right)$ with
$\tilde{X}\left(\tau\right)=Y$.

Next, we prove that $\tilde{X}(t)$ has the same distribution as
$Y^{(t)}$ given $\overline{Z}_t$ by induction, starting at $t=\tau$,
then $t\in\left[\tau_K,\tau\right)$, etcetera.  For $t=\tau$,
  $\tilde{X}\left(\tau\right)=Y$, so that $\tilde{X}(\tau)$ has the
  same distribution as $Y^{\left(\tau\right)}$ given
  $\overline{Z}_\tau$ because of Assumption~\ref{inst}. For the
  induction step, suppose that for $t\in\left[\tau_k,\tau\right]$ (for
  $k=K+1$ read $t=\tau$), $\tilde{X}(t)$ has the same distribution as
  $Y^{(t)}$ given $\overline{Z}_t$. Thus,
  $\tilde{X}\left(\tau_k\right)$ has the same distribution as
  $Y^{\left(\tau_k\right)}$ given $\overline{Z}_{\tau_{k}}$, and hence
  $\tilde{X}\left(\tau_k\right)$ also has the same distribution as
  $Y^{\left(\tau_k\right)}$ given
  $\overline{Z}_{\tau_{k-1}}$. Therefore Assumption~\ref{cond}a
  implies that $F_{Y^{\left(\tau_k\right)}|\overline{Z}_{\tau_{k-1}}}
  \left(\tilde{X}\left(\tau_k\right)\right)$ is uniformly distributed
  on $\left[0,1\right]$ given $\overline{Z}_{\tau_{k-1}}$
  (Lemma~\ref{un} has a formal proof). Then
  $\tilde{X}(t)=F^{-1}_{Y^{(t)}|\overline{Z}_{\tau_{k-1}}}\circ
  F_{Y^{\left(\tau_k\right)}|\overline{Z}_{\tau_{k-1}}}
  \left(\tilde{X}\left(\tau_k\right)\right)$ has distribution function
  $F_{Y^{(t)}|\overline{Z}_{\tau_{k-1}}}$ given
  $\overline{Z}_{\tau_{k-1}}$ (Lemma~\ref{gZ} has a formal proof), so
  also given $\overline{Z}_t$.  That finishes the induction step, so
  that indeed $\tilde{X}(t)$ mimics $Y^{(t)}$ for all
  $t\in\left[0,\tau\right]$.  \hfill $\Box$

\subsection{Discretization and choices of conditional distributions}
\label{discsch}

We return to the continuous-time setting and define a discretization of
the covariate- and treatment process $Z$. Later, we will apply the
result of the previous section to this discretized continuous-time
setting. This section also chooses versions of the conditional distribution
functions given this discretized process.

For $n$ fixed define $\tau_0^{(n)}=0$,
$\tau_1^{(n)}=\frac{1}{2^n}\tau$, $\tau_2^{(n)}=\frac{2}{2^n}\tau$,
..., $\tau_{2^n}^{(n)}=\frac{2^n}{2^n}\tau=\tau$. Consider the grid at
stage $n$ consisting of these points. This way the interval
$\left[0,\tau\right]$ is split up into $2^n$ intervals of equal
length, and when $n$ increases points are added in the middle of these
intervals. For ease of notation, the superscript $^{(n)}$ in
$\tau_k^{(n)}$ is dropped if it is clear which $n$ is meant. Define
$\overline{Z}_t^{(n)}=\left(Z\bigl(\tau_k^{(n)}\bigr):
0\leq\tau_k^{(n)}\leq t\right)$ if $Z$ takes values in a discrete
space,\\ $\overline{Z}_t^{(n)}=\left(1_{\left[\frac{i}{2^n},\frac{i+1}{2^n}\right)}\Bigl(Z\bigl(\tau_k^{(n)}\bigr)\Bigr):0\leq\tau_k^{(n)}\leq
  t,i\in{\mathbb Z}\right)$ if $Z$ takes values in ${\mathbb R}$ and
  \\ $\overline{Z}_t^{(n)}=
  \left(1_{\left[\frac{i}{2^n},\frac{i+1}{2^n}\right)}\Bigl(Z\bigl(\tau_k^{(n)}\bigr)_j\Bigr):0\leq\tau_k^{(n)}\leq
    t,i\in{\mathbb Z},j=1,\ldots,m\right)$ if $Z$ takes values in
    ${\mathbb R}^m$.\\ With this discretization, the information about
    $\overline{Z}_t$ contained in $\overline{Z}_t^{(n)}$ increases
    with $n$: once a grid point is added it stays on the grid for $n$
    larger, and the information about $Z$ in a fixed grid point also
    increases with $n$. Note also that $\overline{Z}_t^{(n)}$ depends
    deterministically on $\overline{Z}_t$, so that no extra randomness
    is necessary to construct $\overline{Z}_t^{(n)}$. Thus
    $\overline{Z}_t^{(n)}$ has the properties promised in the outline
    of the proof, Section~\ref{XYin}.

Next, versions of conditional distributions are chosen. Recall
$\overline{Z}_{\tau_k}$ takes values in the space of cadlag functions
on $[0,\tau_k]$ with the projection $\sigma$-algebra, which is the
same as the Skorohod-$\sigma$-algebra (\cite{Bill68} Theorem~14.5). This space is Polish
(\cite{Bill68} Chapter~3). Therefore, there
exists a conditional distribution
$P_{\overline{Z}_{\tau_k}|\overline{Z}_{\tau_k}^{(n)}}$
(\cite{Bau} Section 10.3 or \cite{Pol2}).  Moreover, $ P\bigl(Y^{(t+h)}\leq
y\big|\overline{Z}_{\tau_k}^{(n)}\bigr)$ $=\int
F_{Y^{(t+h)}|\overline{Z}_{\tau_k}=z}\left(y\right)
dP_{\overline{Z}_{\tau_k}|\overline{Z}_{\tau_k}^{(n)}}\left(z\right)
\;\mbox{\rm a.s.}$. This is a conditional distribution function: it is
non-decreasing in $y$ since all
$F_{Y^{(t+h)}|\overline{Z}_{\tau_k}=z}\left(y\right)$ are
non-decreasing because they are conditional distribution functions,
and because of Lebesgue's Dominated Convergence Theorem
%~\ref{Leb}
the limit for $y\rightarrow -\infty$ equals $0$ and
the limit
for $y\rightarrow \infty$ equals $1$. Therefore, the following choices can
be made:

\begin{nota} \label{choice} We choose fixed conditional distributions
$P_{\overline{Z}_{\tau_{k}}|\overline{Z}_{\tau_k}^{(n)}}$.
I also choose
\begin{equation*} F_{Y^{(t)}|\overline{Z}^{(n)}_{\tau_k}}\left(y\right)=\int
F_{Y^{(t)}|\overline{Z}_{\tau_k}=z}\left(y\right)
dP_{\overline{Z}_{\tau_{k}}|\overline{Z}_{\tau_k}^{(n)}}\left(z\right),
\end{equation*}
with $F_{Y^{(t)}|\overline{Z}_{\tau_k}=z}$ as in Section~\ref{result},
to be the version of the conditional distribution function of
$Y^{(t)}$ given $\overline{Z}_{\tau_k}^{(n)}$ which is used in the
rest of the proof. If $s\in\left(\tau_k,\tau_{k+1}\right)$, the same
version for $F_{Y^{(t)}|\overline{Z}_s^{(n)}}$ is chosen; this is
possible since for $s\in\left(\tau_k,\tau_{k+1}\right)$,
$\overline{Z}_s^{(n)}=\overline{Z}_{\tau_k}^{(n)}$.
\end{nota}
Notice that $Z^{(n)}$ has been constructed with values in a discrete
space. This will assure that the two different expressions for
$D^{(n)}$ in Section~\ref{Dnsec} below are equal except for at a null set
which does not depend on $y$ and $t$.

\subsection{Existence of and two expressions for $D^{(n)}$}
\label{Dnsec}

This section proves existence of $D^{(n)}$ as defined in
equation~(\ref{Dn}), Section~\ref{XYin}. Moreover, two useful formulas
for $D^{(n)}$ are proven. One is used to prove smoothness of
$D^{(n)}$, the other formula is used to prove that $D^{(n)}$ converges
to $D$.

First, existence of $D^{(n)}$ is shown. Fix $n$ and $t$, and choose
$\tau_k^{(n)}$ such that $t\in\bigl[\tau_k^{(n)},\tau_{k+1}^{(n)}\bigr)$. Define
\begin{equation*}
F_h(y)=F_{Y^{(t+h)}|\overline{Z}_{\tau_k}^{(n)}}\left(y\right)
=\int
F_{Y^{(t+h)}|\overline{Z}_{\tau_k}=z}\left(y\right)
dP_{\overline{Z}_{\tau_k}|\overline{Z}_{\tau_k}^{(n)}}\left(z\right).
\end{equation*}
To apply Lemma~\ref{quotc} on $F_h\left(y\right)$,
in $\left(h_0,y_0\right)=\left(0,y\right)$, we check the conditions.
Clearly, $F_h(y)$ is non-decreasing.  We show
that $F_h\left(y\right)$ is differentiable with respect to $y$ with
derivative\linebreak $\int \frac{\partial}{\partial y}
F_{Y^{(t+h)}|\overline{Z}_{\tau_k}=z}\left(y\right)
dP_{\overline{Z}_{\tau_k}|\overline{Z}_{\tau_k}^{(n)}}\left(z\right)$.
For $\omega$ fixed,
$P_{\overline{Z}_{\tau_k}|\overline{Z}_{\tau_k}^{(n)}}$ is a
probability measure on $\overline{{\cal Z}}_{\tau_k}$.  Moreover,
$\frac{\partial}{\partial y}
F_{Y^{(t+h)}|\overline{Z}_{\tau_k}=z}\left(y\right)$ is bounded by
$C_1$, %(Assumptions~\ref{cont1}a and~\ref{bdd+}a)
which is integrable
with respect to
$P_{\overline{Z}_{\tau_k}|\overline{Z}_{\tau_k}^{(n)}}$, and also
$F_{Y^{(t+h)}|\overline{Z}_{\tau_k}=z}\left(y\right)$ is integrable
with respect to
$P_{\overline{Z}_{\tau_k}|\overline{Z}_{\tau_k}^{(n)}}$, since bounded by
$1$. Therefore, $F_h\left(y\right)$ is differentiable with
respect to $y$ with derivative $\int \frac{\partial}{\partial y}
F_{Y^{(t+h)}|\overline{Z}_{\tau_k}=z}\left(y\right)
dP_{\overline{Z}_{\tau_k}|\overline{Z}_{\tau_k}^{(n)}}\left(z\right)$.
With the same reasoning (but with Assumption~\ref{bdd+}b instead of
\ref{bdd+}a), $F_h\left(y\right)$ is differentiable with
respect to $h$ with derivative $\int \frac{\partial}{\partial h}
F_{Y^{(t+h)}|\overline{Z}_{\tau_k}=z}\left(y\right)
dP_{\overline{Z}_{\tau_k}|\overline{Z}_{\tau_k}^{(n)}}\left(z\right)$.
That these derivatives of $F_h(y)$ with respect to $y$ and $h$ are
continuous in $\left(y,h\right)$ follows from Lebesgue's Dominated
Convergence Theorem applied on the expressions we just derived (the
conditions are satisfied because of Assumptions~\ref{cont1}a
and~\ref{bdd+}). Furthermore, $F_0'\left(y\right)=\int
\frac{\partial}{\partial y}
F_{Y^{(t)}|\overline{Z}_{\tau_k}=z}\left(y\right)
dP_{\overline{Z}_{\tau_k}|\overline{Z}_{\tau_k}^{(n)}}\left(z\right)$
is non-zero (Assumption~\ref{sup}b).  Thus the conditions of
Lemma~\ref{quotc} are satisfied for $F_h(y)$, and therefore
$\frac{\partial}{\partial h}
F^{-1}_{Y^{(t+h)}|\overline{Z}^{(n)}_{\tau_k}}\left(y\right)$ exists,
and $D^{(n)}\bigl(y,t;\overline{Z}_t^{(n)}\bigr)$ exists and
satisfies
\begin{eqnarray}
D^{(n)}\left(y,t;\overline{Z}_t^{(n)}\right)
&=&\left.\frac{\partial}{\partial h}\right|_{h=0}
\left(F^{-1}_{Y^{(t+h)}|\overline{Z}_{\tau_k}^{(n)}}\circ
F_{Y^{(t)}|\overline{Z}_{\tau_k}^{(n)}}\right)\left(y\right)\nonumber\\
&=&-\frac{\left.\frac{\partial}{\partial h}\right|_{h=0} \int
F_{Y^{(t+h)}|\overline{Z}_{\tau_k}=z}\left(y\right)
dP_{\overline{Z}_{\tau_k}|\overline{Z}_{\tau_k}^{(n)}}\left(z\right)}
{\frac{\partial}{\partial y}\int
F_{Y^{(t)}|\overline{Z}_{\tau_k}=z}\left(y\right)
dP_{\overline{Z}_{\tau_k}|\overline{Z}_{\tau_k}^{(n)}}\left(z\right)}\nonumber\\
&=&-\frac{\int\left.\frac{\partial}{\partial h}\right|_{h=0}
F_{Y^{(t+h)}|\overline{Z}_{\tau_k}=z}\left(y\right)
dP_{\overline{Z}_{\tau_k}|\overline{Z}_{\tau_k}^{(n)}}
\left(z\right)}{\int\frac{\partial}{\partial y}
F_{Y^{(t)}|\overline{Z}_{\tau_k}=z}\left(y\right)
dP_{\overline{Z}_{\tau_k}|\overline{Z}_{\tau_k}^{(n)}}\left(z\right)}.
\label{Dneq2}
\end{eqnarray}
%where in the second line the previous reasoning is used once
%again. This expression for $D^{(n)}$ will be used in
%Sections~\ref{Sdisc} and~\ref{XXnaf} to prove smoothness of $D^{(n)}$
%in $\left(y,t\right)$, which we need in order to use results about
%differential equations.

Next, the second expression for $D^{(n)}$ is derived. We show that there
exists an $\Omega'\subset \Omega$ with probability one such that
\begin{equation}
D^{(n)}\Bigl(y,t;\overline{Z}_t^{(n)}\Bigr)=
-\frac{E\left[\left.\frac{\partial}{\partial h}\right|_{h=0}
F_{Y^{(t+h)}|\overline{Z}_t}\left(y\right)\big|\overline{Z}_t^{(n)}\right]}
{E\left[\frac{\partial}{\partial y}
F_{Y^{(t)}|\overline{Z}_t}\left(y\right)\big|\overline{Z}_t^{(n)}\right]}
\rule{5ex}{0ex}
\forall\omega\in\Omega' \;\forall y\;\forall t\;\forall n
.\label{Dneq}
\end{equation}

First we choose this $\Omega'$, in such a way that on $\Omega'$
conditional probabilities given $\overline{Z}^{(n)}_{\tau_k}$ are
unique, for all $n$ and $\tau_k$.  Fix $n$ and $\tau_k$ for a moment.
It is known from general theory about conditioning that conditional
probabilities given $\overline{Z}_{\tau_k}^{(n)}=z$ can be written as
a measurable function of $z$. It is also known that conditional
probabilities given $\overline{Z}_{\tau_k}^{(n)}=z$ are almost surely
unique.  Combining these two facts, it follows that conditional
probabilities given $\overline{Z}_{\tau_k}^{(n)}=z$ are unique except
for at $\omega$'s for which
$\overline{Z}_{\tau_k}^{(n)}\left(\omega\right)$ has probability zero,
that is, except for $\omega$'s in
\begin{equation*}
\bigcup_{z:P(\overline{Z}_{\tau_k}^{(n)}=z)=0}
\left\{\omega\in\Omega:
\overline{Z}_{\tau_k}^{(n)}\left(\omega\right)=z\right\}.
\end{equation*}
Since, by construction, $\overline{Z}_{\tau_k}^{(n)}$ takes only
countably many values, this is a countable union of null sets and thus
a null set. Define
\begin{equation}
\Omega'=\Omega\setminus \bigcup_{n\in\mathbb{N}}\bigcup_{k\in\left\{0,\ldots,2^n\right\}}
\bigcup_{z:P(\overline{Z}_{\tau_k}^{(n)}=z)=0}\left\{\omega\in\Omega:
\overline{Z}_{\tau_k}^{(n)}\left(\omega\right)=z\right\}.
\label{omeac}
\end{equation}
This set has probability one since its complement is a countable
union of null sets: ${\mathbb N}$ is countable and for each $n$ there
are only finitely many $k$. On this $\Omega'$ conditional
probabilities given $\overline{Z}^{(n)}_{\tau_k}$ are unique, for all
$n$ and $\tau_k$.

Next, it is shown that equation~(\ref{Dneq}) holds for $\Omega'$ as defined in
equation~(\ref{omeac}). As shown in Section~\ref{discsch}, there
exists a conditional distribution
$P_{\overline{Z}_{t}|\overline{Z}_{\tau_k}^{(n)}}$. For $t\geq
\tau_k$ and $h\geq 0$,
$F_h(y):=P\bigl(Y^{(t+h)}\leq y|\overline{Z}_{\tau_k}^{(n)}\bigr)
=\int F_{Y^{(t+h)}|\overline{Z}_{t}=z}\left(y\right)
dP_{\overline{Z}_{t}|\overline{Z}_{\tau_k}^{(n)}}\left(z\right)
\;\mbox{\rm a.s.}$.  On $\Omega'$ this version is the same as the one
used in the definition of $D^{(n)}$ of equation~(\ref{Dn}), since
conditional probabilities given $\overline{Z}_{\tau_k}^{(n)}$ are
unique on $\Omega'$. Verifying the conditions of Lemma~\ref{quotc}
can be done in exactly the same way as for the first expression for
$D^{(n)}$. Therefore, Lemma~\ref{quotc} implies that for
$\omega\in\Omega'$ and $t\in\left[\tau_k,\tau_{k+1}\right)$,
\begin{eqnarray*}
D^{(n)}\Bigl(y,t;\overline{Z}_t^{(n)}\Bigr)
&=&\left.\frac{\partial}{\partial h}\right|_{h=0}
\left(F^{-1}_{Y^{(t+h)}|\overline{Z}_{\tau_k}^{(n)}}\circ
F_{Y^{(t)}|\overline{Z}_{\tau_k}^{(n)}}\right)\left(y\right)\\
&=&-\frac{\left.\frac{\partial}{\partial h}\right|_{h=0} \int
F_{Y^{(t+h)}|\overline{Z}_{t}=z}\left(y\right)
dP_{\overline{Z}_{t}|\overline{Z}_{\tau_k}^{(n)}}\left(z\right)}
{\frac{\partial}{\partial y}\int
F_{Y^{(t)}|\overline{Z}_{t}=z}\left(y\right)
dP_{\overline{Z}_{t}|\overline{Z}_{\tau_k}^{(n)}}\left(z\right)}\\
&=&-\frac{E\left[\left.\frac{\partial}{\partial h}\right|_{h=0}
F_{Y^{(t+h)}|\overline{Z}_{t}}\left(y\right)\big|
\overline{Z}_{\tau_k}^{(n)}\right]}
{E\left[\frac{\partial}{\partial y}
F_{Y^{(t)}|\overline{Z}_{t}}\left(y\right)\big|
\overline{Z}_{\tau_k}^{(n)}\right]}.
\end{eqnarray*}
Equation~(\ref{Dneq}) follows.

%This second expression for $D^{(n)}$ will be used in
%Section~\ref{DDnaf} to prove convergence of $D^{(n)}$ to $D$.

\subsection{Applying the discrete-time result}
\label{Sdisc}
%$X^{(n)}(t)$ mimics $Y^{(t)}$ given $\overline{Z}_t^{(n)}$.

\begin{lem} \label{detlem} Suppose that Regularity
  Conditions~\ref{sup}--\ref{Lipc} and Consistency Assumption~\ref{inst}
  are satisfied. Then for every $n$ there exists a
  continuous solution $X^{\left(n\right)}(t)$ to the differential
  equation with $D^{(n)}$ with final condition
  $X^{(n)}\left(\tau\right)=Y$.  $X^{(n)}\left(t\right)$ is unique on
  $\Omega'$ of equation~(\ref{omeac}). Furthermore,
  $X^{(n)}\left(t\right)$ has the same conditional distribution as
  $Y^{(t)}$ given $\overline{Z}_t^{(n)}$.
\end{lem}

\noindent{\bf Proof.} Fix $n$. First, we show
that there exists a continuous solution $X^{(n)}$ for which
$X^{(n)}\left(t\right)$ has the same conditional distribution as
$Y^{(t)}$ given $\overline{Z}_t^{(n)}$, using Proposition~\ref{detth}.
Thus we check that the conditional distributions
$F_{Y^{(t)}|\overline{Z}^{(n)}_{\tau_k}}$ of $Y^{(t)}$ given
$\overline{Z}^{(n)}_{\tau_k}$ chosen in Notation~\ref{choice} satisfy
Assumption~\ref{cond}. In the second paragraph of Section~\ref{Dnsec},
we showed that $F_{Y^{(t)}|\overline{Z}^{(n)}_{\tau_k}}\left(y\right)$
is strictly increasing and differentiable with respect to $y$ on
$\left[y_1,y_2\right]$, which accounts for Assumption~\ref{cond}a and
b. %In the last paragraph of Section~\ref{Dnsec}
Just before equation~(\ref{Dneq2}), it was concluded that for
$x\in\left[0,1\right]$ fixed,
$F^{-1}_{Y^{(t)}|\overline{Z}^{(n)}_{\tau_k}}\left(x\right)$ is
differentiable with respect to $t$ on
$\left[\tau_k,\tau_{k+1}\right]$, which accounts for
Assumption~\ref{cond}c. Hence Proposition~\ref{detth} guarantees
existence of a continuous solution $X^{(n)}$ to
$X^{(n)}(t)'=D^{(n)}\left(X^{(n)}(t),t\right)$ with final condition
$X^{(n)}(\tau)=Y$ and with $X^{(n)}(t)\sim Y^{(t)}$ given
$\overline{Z}_t^{(n)}$.

Proposition~\ref{detth} does not imply that $X^{(n)}$ is unique.
Almost sure uniqueness of $X^{(n)}$ follows with Theorem~\ref{difbc2}
in the Appendix along the same lines as uniqueness of $X$ (see
Section~\ref{Dexu}), but using equations~(\ref{Dneq2})
and~(\ref{Dneq}) for $D^{(n)}$ instead of equation~(\ref{Deq}) for
$D$, as follows. Fix $n$ and suppose that
$t\in\left[\tau_k,\tau_{k+1}\right)$. First, it is proven that
$D^{(n)}$ is continuous on
$\left[y_1,y_2\right]\times\left[\tau_k,\tau_{k+1}\right)$ with a
continuous extension to
$\left[y_1,y_2\right]\times\left[\tau_k,\tau_{k+1}\right]$, using
equation~(\ref{Dneq2}). Expression~(\ref{Dneq2}) for $D^{(n)}$ has an
obvious extension $\tilde{D}^{(n)}$ to
$\left[\tau_k,\tau_{k+1}\right]$.  We prove that this
$\tilde{D}^{(n)}$ is continuous on
$\left[y_1,y_2\right]\times\left[\tau_k,\tau_{k+1}\right]$.  To show
that $\int\left.\frac{\partial}{\partial h}\right|_{h=0}
F_{Y^{(t+h)}|\overline{Z}_{\tau_k}=z}\left(y\right)
dP_{\overline{Z}_{\tau_k}|\overline{Z}_{\tau_k}^{(n)}} \left(z\right)$
and
%\linebreak
$\int\frac{\partial}{\partial y}
F_{Y^{(t)}|\overline{Z}_{\tau_k}=z}\left(y\right)
dP_{\overline{Z}_{\tau_k}|\overline{Z}_{\tau_k}^{(n)}}\left(z\right)$
are continuous in $\left(y,t\right)$ Lebesgue's Dominated
Convergence Theorem can be used, as follows.
\begin{equation*}\left.\frac{\partial}{\partial h}\right|_{h=0}
F_{Y^{(t+h)}|\overline{Z}_{\tau_k}=z}\left(y\right)=
\left.\frac{\partial}{\partial
h}\right|_{h=t-\tau_k}F_{Y^{(\tau_k+h)}|\overline{Z}_{\tau_k}=z}\left(y\right)
\end{equation*}
and
\begin{equation*}
\frac{\partial}{\partial y}
F_{Y^{(t)}|\overline{Z}_{\tau_k}=z}\left(y\right)
=\frac{\partial}{\partial y}
F_{Y^{(\tau_k+\left(t-\tau_k\right))}|\overline{Z}_{\tau_k}=z}\left(y\right)
\end{equation*}
are continuous in $\left(y,t\right)$ because of
Assumption~\ref{cont1}a. Both these derivatives are bounded because of
Assumption~\ref{bdd+}. Therefore Lebesgue's Dominated Convergence
Theorem implies that the integrals of these derivatives with respect
to the measure
$\mu=P_{\overline{Z}_{\tau_k}|\overline{Z}_{\tau_k}^{(n)}}$ are
continuous in $\left(y,t\right)$. Because of Assumption~\ref{sup}b the
denominator of $\tilde{D}^{(n)}$ is non-zero for
$y\in\left[y_1,y_2\right]$, so that indeed $\tilde{D}^{(n)}$ is
continuous in $\left(y,t\right)$ on
$\left[y_1,y_2\right]\times\left[\tau_k,\tau_{k+1}\right]$.

Next, it is shown that $D^{(n)}$ is Lipschitz continuous in $y$ on
$\left[y_1,y_2\right]\times\left[\tau_k,\tau_{k+1}\right]$ with
Lipschitz constant $L_2/\varepsilon + C_2L_1/\varepsilon^2$ for all
$\omega\in\Omega'$, with $\Omega'$ as in
equation~(\ref{omeac}). Expression~(\ref{Dneq2}) for $D^{(n)}$ on
$\Omega'$ has an obvious extension $\tilde{D}^{(n)}$ to
$\left[\tau_k,\tau_{k+1}\right]$. That this $\tilde{D}^{(n)}$ is
Lipschitz continuous in $y$ on
$\left[y_1,y_2\right]\times\left[\tau_k,\tau_{k+1}\right]$ with
Lipschitz constant $L_2/\varepsilon + C_2L_1/\varepsilon^2$ on
$\Omega'$ follows the same way as for $D$ in Section~\ref{Dexu}.
Because of Assumption~\ref{sup}c, the denominator is bounded away from
$0$ for $y\in\left[y_1,y_2\right]$, and because of
Assumption~\ref{sup}a, the numerator is equal to zero for $y=y_1$ and
for $y=y_2$. Hence, on $\Omega'$,
$\tilde{D}^{(n)}(y_1,t)=\tilde{D}^{(n)}(y_2,t)=0$. Therefore
Theorem~\ref{difbc2} implies that, on $\Omega'$, there exists a unique
solution to the differential equation with $\tilde{D}^{(n)}$ on
$\left[\tau_k,\tau_{k+1}\right]$, and this solution stays in
$\left[y_1,y_2\right]$. Since for $n$ fixed there are only finitely
many $\tau_k$, the same is true on $\left[0,\tau\right]$.  \hfill
$\Box$

%$ $
%
%\noindent Notice that $X^{(n)}(t)\sim Y^{(t)}$ given
%$\overline{Z}_{\tau_k}^{(n)}$ implies that for every $t$,
%$X^{(n)}(t)\in\left[y_1,y_2\right]$ almost surely. In fact, the proof
%of Proposition~\ref{detth} combined with the proof of the uniqueness
%statement of Lemma~\ref{detlem} leads to a stronger statement:
%$X^{(n)}(t)\in \left[y_1,y_2\right]$ for every $t$ for every
%$\omega\in\Omega'$, with $\Omega'$ the set of probability one which
%was defined in equation~(\ref{omeac}).

\subsection{Bounding the difference between $X$ and $X^{(n)}$ in
terms of $D$ and $D^{(n)}$}
\label{XXnaf}

To bound the difference between $X$ and $X^{(n)}$ in terms of $D$ and
$D^{(n)}$, Theorem~\ref{difbc2} is applied on $y=X\left(t\right)$ and
$z=X^{(n)}\left(t\right)$. Since we need that both $D$ and $D^{(n)}$
are continuous, we apply Theorem~\ref{difbc2} on the intervals between
the jumps of $Z$ and the grid points $\tau_k^{(n)}$. Fix $n$ and
restrict $\omega$ to $\omega\in\Omega'$, with $\Omega'$ the set of
probability one as defined in equation~(\ref{omeac}), so that the
expression for $D^{(n)}$ of equation~(\ref{Dneq}) can be used. The
bound will thus hold almost surely. To focus attention on the
differential equations, the $\overline{Z}_t$'s and
$\overline{Z}_t^{(n)}$'s in $D$ and $D^{(n)}$ are skipped below.

Suppose that $\left(t_1,t_2\right)$ is such an interval including no
jumps of $Z$ and no grid points at stage $n$. We check the conditions
of Theorem~\ref{difbc2} for $y=X\left(t\right)$ and
$z=X^{(n)}\left(t\right)$. Section~\ref{Dexu} already showed that
$D:\left[y_1,y_2\right]\times \left[t_1,t_2\right)\rightarrow {\mathbb
  R}$ has a continuous extension $\tilde{D}:\left[y_1,y_2\right]\times
\left[t_1,t_2\right]\rightarrow {\mathbb R}$ which satisfies the
conditions of Theorem~\ref{difbc2}, with $C$ the constant function
$L_2/\varepsilon+C_2 L_1/\varepsilon^2$, and in the proof of
Lemma~\ref{detlem} in Section~\ref{Sdisc}, it was shown that on $\Omega'$ the
same is true for $D^{(n)}$. Therefore Theorem~\ref{difbc2} implies
that for $t\in\left[t_1,t_2\right]$, with $C=L_2/\varepsilon+C_2
L_1/\varepsilon^2$ as above,
\begin{eqnarray} \bigl|X^{(n)}(t) -X(t)\bigr|
&\leq& e^{\int_{t}^{t_2} C\; ds}\;
\bigl|X^{(n)}\left(t_{2}\right)-X\left(t_{2}\right)\bigr|\nonumber\\
&&\;\;\;+\int_{t}^{t_2} e^{\int_{t}^{s} C\; d\eta}\;
\bigl|D\bigl(X^{(n)}\left(s\right),s\bigr)
-D^{(n)}\bigl(X^{(n)}\left(s\right),s\bigr)\bigr|\;ds\nonumber\\
&=& e^{C\cdot\left(t_{2}-t\right)}\;
\bigl|X^{(n)}\left(t_{2}\right)-X\left(t_{2}\right)\bigr|\nonumber\\
&&\;\;\;+\int_{t}^{t_{2}} e^{C\cdot\left(s-t\right)}\;
\bigl|D\bigl(X^{(n)}\left(s\right),s\bigr)
-D^{(n)}\bigl(X^{(n)}\left(s\right),s\bigr)\bigr|\;ds. \label{afsch}
\end{eqnarray}

If $Z$ does not jump in $\left[\left(1-1/2^n\right)\tau,\tau\right]$, (\ref{afsch}) can be applied on
$\left[\left(1-1/2^n\right)\tau,\tau\right]$, and since
$X^{(n)}\left(\tau\right)=X\left(\tau\right)=Y$ it follows that on
$\left[\left(1-1/2^n\right)\tau,\tau\right]$,
\begin{equation} \bigl|X^{(n)}(t) -X(t)\bigr|
\leq \int_{t}^{\tau} e^{C\cdot \left(s-t\right)}\;
\bigl|D\bigl(X^{(n)}\left(s\right),s\bigr)-
D^{(n)}\bigl(X^{(n)}\left(s\right),s\bigr)\bigr|\; ds.
\label{afsch1}
\end{equation}
If $Z$ does not jump after $\left(1-2/2^n\right)\tau$ one can also
apply~(\ref{afsch}) on $\left[\left(1-2/2^n\right)\tau,
\left(1-1/2^n\right)\tau\right]$, and using equation~(\ref{afsch1})
for $t=\left(1-1/2^n\right)\tau$, it follows that equation~(\ref{afsch1})
also holds on $\left[\left(1-2/2^n\right)\tau,
\left(1-1/2^n\right)\tau\right]$:
\begin{eqnarray*} \lefteqn{\bigl|X^{(n)}(t) -X(t)\bigr|}\\
&\leq& e^{C\cdot\left(\left(1-\frac{1}{2^n}\right)\tau-t\right)}
\int_{\left(1-\frac{1}{2^n}\right)\tau}^{\tau}
e^{C\cdot\left(s-\left(1-\frac{1}{2^n}\right)\tau\right)}
\bigl|D\bigl(X^{(n)}\left(s\right),s\bigr)-
D^{(n)}\bigl(X^{(n)}\left(s\right),s\bigr)\bigr| ds\\
&&+\int_{t}^{\left(1-\frac{1}{2^n}\right)\tau}
e^{C\cdot\left(s-t\right)} \;
\bigl|D\bigl(X^{(n)}\left(s\right),s\bigr)-
D^{(n)}\bigl(X^{(n)}\left(s\right),s\bigr)\bigr|\; ds\\
&=&\int_{t}^{\tau}
e^{C\cdot\left(s-t\right)} \;
\bigl|D\bigl(X^{(n)}\left(s\right),s\bigr)-
D^{(n)}\bigl(X^{(n)}\left(s\right),s\bigr)\bigr|\; ds.
\end{eqnarray*}
If $Z$ does not jump in $\left(\left(1-m/2^n\right)\tau,\tau\right]$
and $t\in\left(\left(1-m/2^n\right)\tau,\tau\right]$ then, with the
same reasoning, equation~(\ref{afsch1}) holds on
$t\in\left(\left(1-m/2^n\right)\tau,\tau\right]$.
Suppose now that $Z$ jumps in $\left(\left(1-(m+1)/2^n\right)\tau,
\left(1-m/2^n\right)\tau\right]$. Then this interval can be split up
into the part before and the part after the jump, so that, again with
the same reasoning as before and since both $X^{(n)}$ and $X$ are
continuous in $t$, equation~(\ref{afsch1}) still holds.

With probability one there are at most finitely many jump times of
$Z$, so that equation~(\ref{afsch1}) holds almost surely for all $t$, and even
\begin{eqnarray}
\sup_{t\in\left[0,\tau\right]}\bigl|X^{(n)}(t)
-X(t)\bigr|
%}\nonumber\\
%&&\hspace{1.2cm}
&\leq& \sup_{t\in\left[0,\tau\right]}
\int_t^\tau e^{C\cdot\left(s-t\right)}
\bigl|D\bigl(X^{(n)}\left(s\right),s\bigr)-
D^{(n)}\bigl(X^{(n)}\left(s\right),s\bigr)\bigr| ds\nonumber\\
%&&\hspace{1.2cm}
&=& \int_0^\tau e^{C\cdot s}
\bigl|D\bigl(X^{(n)}\left(s\right),s\bigr)-
D^{(n)}\bigl(X^{(n)}\left(s\right),s\bigr)\bigr| ds\;\;\;\mbox{\rm a.s.}.
\label{Xschat}
\end{eqnarray}

\subsection{Convergence of $D^{(n)}$ to $D$}
\label{DDnaf}

This section proves that
$D^{(n)}\bigl(y,t;\overline{Z}_t^{(n)}\bigr)$ converges almost surely
to $D\left(y,t;\overline{Z}_t\right)$, for fixed
$\left(y,t\right)\in\left[y_1,y_2\right]\times\left[0,\tau\right]$. From
equations~(\ref{Deq}) and (\ref{Dneq}) it follows that
\begin{equation*} D\left(y,t;\overline{Z}_t\right)=
-\frac{\left.\frac{\partial}{\partial h}\right|_{h=0}
F_{Y^{(t+h)}|\overline{Z}_t}\left(y\right)}{\frac{\partial}{\partial
y}F_{Y^{(t)}|\overline{Z}_t}\left(y\right)}
\end{equation*}
and
\begin{equation*} D^{(n)}\Bigl(y,t;\overline{Z}_t^{(n)}\Bigr)
=-\frac{E\left[\left.\frac{\partial}{\partial h}\right|_{h=0}
    F_{Y^{(t+h)}|\overline{Z}_t}\left(y\right)\big|\overline{Z}_t^{(n)}\right]}
{E\left[\frac{\partial}{\partial y}
    F_{Y^{(t)}|\overline{Z}_t}\left(y\right)\big|\overline{Z}_t^{(n)}\right]}
\;\;\;\;\mbox{\rm a.s.}.\end{equation*} L\'evy's Upward Theorem (see
e.g.~\cite{Wil} page 134) can be applied to the
denominator and the numerator of $D^{(n)}$, since both
$\left.\frac{\partial}{\partial h}\right|_{h=0}
F_{Y^{(t+h)}|\overline{Z}_t}\left(y\right)$ and
$\frac{\partial}{\partial y} F_{Y^{(t)}|\overline{Z}_t}\left(y\right)$
are bounded (Assumption~\ref{bdd+}). L\'evy's Upward Theorem leads to
\begin{equation*} E\left[\left.\frac{\partial}{\partial h}\right|_{h=0}
F_{Y^{(t+h)}|\overline{Z}_t}\left(y\right)\Big|
\overline{Z}_t^{(n)}\right]\rightarrow
E\left[\left.\frac{\partial}{\partial h}\right|_{h=0}
F_{Y^{(t+h)}|\overline{Z}_t}\left(y\right)\bigg|
\sigma\left(\cup_{n=1}^\infty \overline{Z}_t^{(n)}\right) \right]
\;\mbox{\rm a.s.}
%\label{h1}
\end{equation*}
and
\begin{equation*} E\left[\frac{\partial}{\partial y}
F_{Y^{(t)}|\overline{Z}_t}\left(y\right)\Big|\overline{Z}_t^{(n)}\right]
\rightarrow E\left[\frac{\partial}{\partial
y}F_{Y^{(t)}|\overline{Z}_t}\left(y\right)\bigg|
\sigma\left(\cup_{n=1}^\infty\overline{Z}_t^{(n)}\right) \right]
\;\mbox{\rm a.s.}
%\label{h2}
\end{equation*}
as $n\rightarrow\infty$. The conditioning on
$\sigma\bigl(\cup_{n=1}^\infty\overline{Z}_t^{(n)}\bigr)$ can be
replaced by conditioning on $\overline{Z}_t$ in both expressions,
because of Lemma~\ref{Zt-} in the Appendix. Since moreover the
denominators are bounded away from $0$ (Assumption~\ref{sup}c), the
Continuous Mapping Theorem implies that, for fixed
$\left(y,t\right)\in\left[y_1,y_2\right]\times\left[0,\tau\right]$,
\begin{equation}
D^{(n)}\Bigl(y,t;\overline{Z}_t^{(n)}\Bigr)
\rightarrow D\left(y,t;\overline{Z}_t\right) \label{Dconv}\;\mbox{\rm a.s.}.
\end{equation}

\subsection{$X^{(n)}(t)$ converges to $X(t)$ and $X(t)$ is measurable}
\label{XnX}

To show that $X^{(n)}(t)$ converges almost surely to $X(t)$ and
that $X(t)$ is measurable, the bound
of equation~(\ref{Xschat}) and almost sure convergence of
$D^{(n)}(y,t)$ to $D(y,t)$ for $(y,t)$ fixed of
equation~(\ref{Dconv}) are the starting point.

First it is proven that for $s$ fixed,
$D^{(n)}\left(X^{(n)}\left(s\right),s\right)-
D\left(X^{(n)}\left(s\right),s\right)$ converges almost surely to $0$.
Recall from Section~\ref{Dexu} that
$D:\left[y_1,y_2\right]\times \left[t_1,t_2\right)\rightarrow {\mathbb
  R}$ has a continuous extension $\tilde{D}:\left[y_1,y_2\right]\times
\left[t_1,t_2\right]\rightarrow {\mathbb R}$ which is Lipschitz
continuous in $y$ with Lipschitz constant $L_2/\varepsilon+C_2
L_1/\varepsilon^2$. Recall also that in the proof of
Lemma~\ref{detlem} in Section~\ref{Sdisc} it was shown that on $\Omega'$,
the set of probability one of equation~(\ref{omeac}), the same is true
for $D^{(n)}$. Therefore, the pointwise almost sure convergence of
$D^{(n)}(y,t)$ to $D(y,t)$ of equation~(\ref{Dconv}) implies that for fixed $s$
indeed
\begin{equation}\label{Dnfix}\bigl|D^{(n)}\bigl(X^{(n)}\left(s\right),s\bigr)-
D\bigl(X^{(n)}\left(s\right),s\bigr)\bigr|\rightarrow 0
\;\;\;\;\;\mbox{\rm a.s.}\end{equation}
(for details see Web-Appendix~\ref{Conv}).

To show that equation~(\ref{Dnfix}) implies that the bound of
(\ref{Xschat}) converges almost surely to $0$, define
\begin{equation*}
A=\left\{\left(s,\omega\right)\in\left[0,\tau\right]\times \Omega:\bigl|
D^{(n)}\bigl(X^{(n)}\left(s\right),s\bigr)-
D\bigl(X^{(n)}\left(s\right),s\bigr)\bigr|\rightarrow 0\right\},
\end{equation*}
with $A_s$ its section at $s$ and $A_\omega$ its section at
$\omega$. Then
\begin{equation*}
A_s=\left\{\omega\in\Omega:\bigl|D^{(n)}\bigl(X^{(n)}\left(s\right),s\bigr)-
D\bigl(X^{(n)}\left(s\right),s\bigr)\bigr|\rightarrow 0\right\}
\end{equation*}
has probability one because equation~(\ref{Dnfix}).
Therefore, using Fubini's Theorem, with $\lambda$ the Lebesgue-measure
on
$\left[0,\tau\right]$,
\begin{eqnarray*} \left(\lambda \times P\right)\left(A\right)
&=&\int_{\left(0,\tau\right)}
P\left(A_s\right)d\lambda\left(s\right)\\
&=&\int_{\left(0,\tau\right)}1d\lambda\left(s\right)=\tau.
\end{eqnarray*}
Also by Fubini's Theorem,
\begin{equation*} \left(\lambda \times P\right)\left(A\right)=\int
\lambda\left(A_\omega\right)dP\left(\omega\right) ,\end{equation*} so
that since $\lambda\left(A_\omega\right)\leq \tau$,
$\lambda\left(A_\omega\right)=\tau$ $P$-almost everywhere. This shows
that for $P$-almost all $\omega$, $A_\omega$ has measure $\tau$. So
for $P$-almost all $\omega$,
$\left|D\left(X^{(n)}\left(s\right),s\right)-
D^{(n)}\left(X^{(n)}\left(s\right),s\right)\right|$ converges to $0$
for $\lambda$-almost all $s$. Moreover, because of
expression~(\ref{Deq}) for $D$ and expression~(\ref{Dneq}) for
$D^{(n)}$ and Assumptions~\ref{bdd+}b and~\ref{sup}c, $e^{C\cdot
s}\left|D\left(\cdot,s\right)-D^{(n)}\left(\cdot,s\right)\right|$ is
bounded by $2 e^{C\cdot \tau} C_2/\varepsilon$ on $\Omega'$. Therefore
for almost all $\omega$ Lebesgue's Dominated Convergence Theorem can
be applied on the integral of $e^{C\cdot
s}\left|D\left(X^{(n)}\left(s\right),s\right)-
D^{(n)}\left(X^{(n)}\left(s\right),s\right)\right|$ with respect to
$\lambda$, $\int_{\left[0,\tau\right]} e^{C\cdot
s}\left|D\left(X^{(n)}\left(s\right),s\right)-
D^{(n)}\left(X^{(n)}\left(s\right),s\right)\right| ds$, implying that
for almost all $\omega$ this integral converges to $0$ as
$n\rightarrow \infty$. With equation~(\ref{Xschat}), this
implies that
\begin{equation} \sup_{t\in\left[0,\tau\right]}\bigl|X^{(n)}(t) -X(t)\bigr|
\rightarrow 0 \;\;\;\;\mbox{\rm a.s.}.
\end{equation}

Since the almost sure limit of a sequence of random variables is
measurable if the $\sigma$-algebra is complete, measurability of
$X(t)$ follows immediately from measurability of the $X^{(n)}$.

\subsection{Conclusion}
\label{Concl}

This section shows that since $X^{(n)}(t)\sim Y^{(t)}$ given
$\overline{Z}_t^{(n)}$ (see Section~\ref{Sdisc}) and
$X^{(n)}(t)\rightarrow X(t)$ $\;\mbox{\rm a.s.}$ (see
Section~\ref{XnX}), $X(t)\sim Y^{(t)}$ given $\overline{Z}_t$. This
completes the proof.

It is well-known (see e.g.~\cite{Vaart}; Lemma~\ref{anpm} provides a
formal proof for this conditional version) that $X(t)\sim Y^{(t)}$
given $\overline{Z}_t$ if
\begin{equation*}
E\left[f\left(X(t)\right)|\overline{Z}_t\right]-
E\bigl[f\bigl(Y^{(t)}\bigr)|\overline{Z}_t\bigr]=0\;\mbox{\rm a.s.}
\end{equation*}
for every bounded Lipschitz continuous function
$f:\mathbb{R}\rightarrow\mathbb{R}$. Suppose without loss of
generality that $f$ is bounded by $1$ and has Lipschitz constant
$L$. Then, using the triangle inequality,
\begin{eqnarray*}
\bigl|E\left[f\left(X(t)\right)|\overline{Z}_t\right]-
E\bigl[f\bigl(Y^{(t)}\bigr)|\overline{Z}_t\bigr]\bigr|
&\leq&
\left|E\left[f\left(X(t)\right)|\overline{Z}_t\right]-
E\left[f\left(X(t)\right)\big|\overline{Z}_t^{(n)}\right]\right|\\
&&\;\;\; + \left|E\left[f\left(X(t)\right)\big|\overline{Z}_t^{(n)}\right]-
E\left[f\bigl(X^{(n)}(t)\bigr)\big|\overline{Z}_t^{(n)}\right]\right|\\
&&\;\;\; + \left|E\left[f\bigl(X^{(n)}(t)\bigr)\big|\overline{Z}_t^{(n)}\right]-
E\bigl[f\bigl(Y^{(t)}\bigr)|\overline{Z}_t\bigr]\right|.
\end{eqnarray*}
Because of Jensen's inequality, the second term is bounded by
$E\left[\bigl|f\left(X(t)\right)- f\bigl(X^{(n)}(t)\bigr)\bigr|
\,\big|\overline{Z}_t^{(n)}\right]$, which is bounded by
$E\left[L\bigl|X(t)- X^{(n)}(t)\bigr|\wedge
2\big|\overline{Z}_t^{(n)}\right]$ since $f$ is Lipschitz continuous with
Lipschitz constant $L$ and bounded by $1$. Because $X^{(n)}(t)\sim
Y^{(t)}$ given $\overline{Z}_t^{(n)}$, the third term is equal to
$\left|E\left[f\bigl(Y^{(t)}\bigr)\big|\overline{Z}_t^{(n)}\right]-
E\bigl[f\bigl(Y^{(t)}\bigr)|\overline{Z}_t\bigr]\right|$. Therefore,
\begin{eqnarray}
\bigl|E\left[f\left(X(t)\right)|\overline{Z}_t\right]-
E\bigl[f\bigl(Y^{(t)}\bigr)|\overline{Z}_t\bigr]\bigr|
&\leq&
\left|E\left[f\left(X(t)\right)|\overline{Z}_t\right]-
E\left[f\left(X(t)\right)\big|\overline{Z}_t^{(n)}\right]\right|\nonumber\\
&&\;\;\; + E\left[L\bigl|X(t)-
X^{(n)}(t)\bigr|\wedge 2\big|\overline{Z}_t^{(n)}\right]\nonumber\\
&&\;\;\; + \left|E\left[f\bigl(Y^{(t)}\bigr)\big|\overline{Z}_t^{(n)}\right]-
E\bigl[f\bigl(Y^{(t)}\bigr)|\overline{Z}_t\bigr]\right| \;\mbox{\rm a.s.}.
\label{ffen}
\end{eqnarray}
We show that the right hand side converges in probability to zero. On
the first and the last term, L\'evy's Upward Theorem (see
e.g.~\cite{Wil} page 134) can be applied, since the integrands are
bounded by 1. L\'evy's Upward Theorem leads to
\begin{equation*}
E\left[f\left(X(t)\right)\big|\overline{Z}_t^{(n)}\right]
\rightarrow E\left[f\left(X(t)\right)\big|
\sigma\left(\cup_{n=1}^\infty\overline{Z}_t^{(n)}\right)\right]
\end{equation*}
and
\begin{equation*}
E\left[f\bigl(Y^{(t)}\bigr)\big|\overline{Z}_t^{(n)}\right]
\rightarrow E\left[f\bigl(Y^{(t)}\bigr)\big|
\sigma\left(\cup_{n=1}^\infty\overline{Z}_t^{(n)}\right)\right]
\end{equation*}
as $n\rightarrow\infty$. Thus, with Lemma~\ref{Zt-} in the Appendix,
both the first and the last term of equation~(\ref{ffen}) converge to
$0$ almost surely. The second term converges to $0$ in probability
since it is almost surely non-negative and its expectation converges
to $0$:
\begin{equation*} E\left(E\left[L\bigl|X(t)-
X^{(n)}(t)\bigr|\wedge 2\big|\,\overline{Z}_t^{(n)}\right]\right) =
E\left(L\bigl|X(t)- X^{(n)}(t)\bigr|\wedge 2\right) \rightarrow 0
\end{equation*}
because of Lebesgue's Dominated Convergence Theorem and the fact
that $X^{(n)}(t)$ converges almost surely to $X(t)$.

Thus $\left|E\left[f\left(X(t)\right)|\overline{Z}_t\right]-
E\left[f\left(Y^{(t)}\right)|\overline{Z}_t\right]\right|$ is bounded
by a random variable which converges in probability to $0$. Hence,
this first random variable is almost surely equal to $0$. Therefore,
indeed $X(t)$ mimics $Y^{(t)}$ in the sense that $X(t)$ has the same
distribution as $Y^{(t)}$ given $\overline{Z}_t$.

\subsection[Discrete-continuous time]{Mimicking counterfactual outcomes:
discrete-continuous time}
\label{Disc2}

In certain situations there are specific times $t$ with $P\left(t
\;{\rm is} \;{\rm a}\; {\rm jump}\;{\rm time}\;{\rm of}\;Z\right)>
0$. For finitely many such times $t$, the proof in Section~\ref{XY}
can be adapted by adding these finitely many times to the grid, for
each $n$.
%Enumerate the times $t$ for which $P\left(t \;{\rm
%is} \;{\rm a}\; {\rm jump}\;{\rm point}\;{\rm of}\;Z\right)> 0$
%and add them e.g.\ one in each step (indicated by $n$ in the preceding
%proof) to the grid (and never leave them out afterwards). That way
%these $t$ will all be on the grid in the limit.
%Just the proof of Lemma~\ref{Zt-} has to be adapted. But that is easy
%since if all $t$ for which $P\left(t \;{\rm is} \;{\rm a}\; {\rm
%jump}\;{\rm time}\;{\rm of}\;Z\right)> 0$ are
%eventually
%on the grid, then also for these $t$,
%$\sigma\bigl(\overline{Z}_t^{(n)}\bigr)\uparrow
%\sigma\left(\overline{Z}_t\right)$.

%thought countably many at first, but now: Because of
%differentiability conditions I think finitely many, not
%countably many. With countably many there is a limit point and if it's
%also a limit from the right natural differentiability conditions on
%$F_{Y^{(t+h)}|\overline{Z}_t}$ get really weird there I'm afraid. I'd
%rather stick to finitely many for security's sake; then I'm pretty
%convinced it is ok.

\section{Mimicking counterfactual survival outcomes}
\label{XT}

\subsection{Introduction}
\label{tintro}

This section indicates how to prove that $X(t)$ mimics $Y^{(t)}$
in the sense that $X(t)$ has the same distribution as $Y^{(t)}$ given
the covariate- and treatment history $\overline{Z}_t$, under
conditions aimed at survival. The conditions are similar to the ones
in Section~\ref{XY}, but adapted to survival as the outcome of
interest. The proof also follows roughly the same lines as the one for
other outcomes, but some changes are necessary. A full proof can be
found in Web-Appendix~\ref{XTapp}.

If covariates and treatment were measured at time $t$, it cannot be
avoided to include in $\overline{Z}_t$ whether or not a person was
alive at time $t$: what are a person's covariates if he or she is dead?
Therefore we include in $Z(t)$ an indicator for whether or not a
person is alive at time $t$.
%To be more precise: in the discrete-time case, where
%$\overline{Z}=\left(Z(t): t\in\left[0,\tau\right]\right)$ can be fully
%described by its values at finitely many fixed times
%$0=\tau_0<\tau_1<\tau_2<\ldots<\tau_K<\tau_{K+1}=\tau$, we add an
%indicator of whether or not the person is alive at time $\tau_k$ to
%$Z\left(\tau_k\right)$. In the continuous-time case we add an
%indicator of whether or not a person is alive at time $t$ to $Z(t)$
%for every $t$.
Thus if a person died at or before time $t$, the survival time can be
read from $\overline{Z}_t$.

The conditions in Section~\ref{XY} usually exclude survival as the
outcome of interest, since if the outcome is survival the Support
Condition~\ref{sup}, saying that all $F_{Y^{(t+h)}|\overline{Z}_t}$
have the same bounded support $\left[y_1,y_2\right]$, will not hold:
$\overline{Z}_t$ includes the covariate-measurements and treatment
until time $t$, and given that a person is dead at time $t$ and given
his or her survival time, the distribution of this survival time
cannot have the fixed support $\left[y_1,y_2\right]$, independent of
$t$. Also given that a person is alive at time $t$, the survival time
often does not have the fixed support $\left[y_1,y_2\right]$: one
often expects that $t$ is the left limit of the support, and obviously
the left limit of the support should be greater than or equal to $t$.

I make two extra assumptions. The first is a straightforward
consistency assumption, stating that stopping treatment after death
does not change the survival time. The second extra assumption states
that there is no instantaneous effect of treatment at the time the
person died (notice that the difference between $Y^{\left(Y\right)}$,
the outcome with treatment stopped at the survival time $Y$, and $Y$
is in treatment \emph{at} time $Y$).

\begin{assu} \label{ttcons}
\emph{(consistency)}.  $Y^{(t)}=Y$ on $\left\{\omega:Y\leq
t\right\}\cup\left\{\omega:Y^{(t)}\leq t\right\}$.
\end{assu}

\begin{assu} \label{ttinst}
\emph{(no instantaneous effect of treatment at the
time the person died)}. $Y^{(t)}=Y\;\;{\rm on}\;\left\{\omega:Y=
t\right\}\cup\left\{\omega:Y^{(t)}=t\right\}.$
\end{assu}
\noindent As can be expected, these assumptions imply that treatment
in the future does not cause or prevent death at present, see
Web-Appendix~\ref{XTapp}.

%\begin{lem} \label{conslem} Under Assumptions~\ref{ttcons} (consistency)
%and~\ref{ttinst} (no instantaneous effect of treatment at the time the
%person died),
%\begin{enumerate}[a)]
%\item For all $h\geq 0$: $Y^{(t+h)}=Y$ on $\left\{\omega:Y\leq t\right\}\cup\cup_{h\geq
%0}\left\{\omega:Y^{(t+h)}\leq t\right\}$.
%\item For all $\left(y,t,h\right)$ with $y\leq t+h$ and $h\geq 0$:
%$\left\{\omega: Y^{(t+h)}\leq y\right\}=\left\{\omega: Y\leq
%y\right\}$.
%\end{enumerate}
%\end{lem}

%\noindent We will only use versions of conditional
%distributions which are consistent with Lemma~\ref{conslem} in the
%sense that $F_{Y^{(t+h)}|\overline{Z}_t}\left(y\right)=
%F_{Y|\overline{Z}_t}\left(y\right)$ for all $y\leq t+h$, $h\geq 0$ and
%$\omega\in\Omega$.
%was: are ``natural'' considering Lemma.. in the sense the e.g.... I
%checked that this can be
%``in the sense that'' instead of ``for example''

For survival outcomes this article uses the following minor adaptation of the
definition of $D$,
\begin{equation} \label{tdefd}
D\left(y,t;\overline{Z}_t\right)=\left\{\begin{array}{l} 0
\;\;\;\;\;\;\;\;\;\;\;\;\;\;\;\;\;\; {\rm if}\; \overline{Z}_t\;
{\rm indicates}\; {\rm the}\; {\rm person}\; {\rm is}\;{\rm dead}\;{\rm
at}\; t\;\; {\rm or}\;y<t\\
\left.\frac{\partial}{\partial h}\right|_{h=0}
\Bigl(F_{Y^{(t+h)}| \overline{Z}_t}^{-1}\circ
F_{Y^{(t)}|\overline{Z}_t}\Bigr)\left(y\right)\;\;\;\;\;\;\;\;\;\;\;\;\,{\rm otherwise,}\;
{\rm for}\;y>t\\
\lim_{y\downarrow t}D\left(y,t;\overline{Z}_t\right)
\;\;\;\;\;\;\;\;\;\;\;\;\;\;\;\;\;\;\;\;\;\;\;\;\;\;\;\;\;\;\;\;\;\;\;
{\rm otherwise,}\;{\rm for}\;y=t,
\end{array}\right.\end{equation}
as we explain now.
%Because of the reasons just below
First remark that considering the interpretation of
$D\left(y,t;\overline{Z}_t\right)$ as the infinitesimal effect of a
short duration of treatment directly after $t$ on survival,
$D\left(y,t;\overline{Z}_t\right)$ should be zero if $\overline{Z}_t$
indicates the person is dead at time $t$. Although in that case
indeed $F_{Y^{(t+h)}|\overline{Z}_t}$ and $F_{Y^{(t)}|\overline{Z}_t}$
are almost surely the same for every $h\geq 0$, since withholding
treatment after death does not change the survival time,
$F_{Y^{(t+h)}|\overline{Z}_t}^{-1}$ will often not exist. Therefore if
$\overline{Z}_t$ indicates the person is dead at time $t$, this article just
formally defines $D\left(y,t;\overline{Z}_t\right)$ to be zero.  Next
consider $y< t$. Notice that considering the interpretation of
$D\left(y,t;\overline{Z}_t\right)$ as the infinitesimal effect of
treatment directly after time $t$ on the survival-quantile $y$,
$D\left(y,t;\overline{Z}_t\right)$ should be zero for $y< t$ since
treatment at or after time $t$ should not cause or prevent death at or
before time $t$, so it should not affect quantiles of the survival curve
before time $t$. Indeed if $\overline{Z}_t$ indicates that the person is
alive at time $t$, $F_{Y^{(t+h)}|\overline{Z}_t}\left(y\right)=
F_{Y|\overline{Z}_t}\left(y\right)=0$ for $y\leq t$ for all $h\geq 0$,
but also for these $y$, $F_{Y^{(t+h)}|\overline{Z}_t}^{-1}\left(y\right)$ often does not
exist. Therefore, this article defines $D\left(y,t;\overline{Z}_t\right)=0$ for
$y<t$. In order to make $D$ continuous on $y\geq t$ in between the
jump times of $Z$, we define
$D\left(t,t;\overline{Z}_t\right)=\lim_{y\downarrow
  t}D\left(y,t;\overline{Z}_t\right)$. This limit exists under the
conditions in Section~\ref{tresult}. It is not necessarily equal to
zero.

Notice that the area where $D$ is possibly non-zero is
$\left(y,t\right)\in\left[0,\infty\right)\times
\left[0,\min\left\{Y,\tau\right\}\right]:y\geq t$. Therefore if
$Y<\tau$, the solution to the differential equation $X(t)$ is equal to
$Y$ for $t\in\left[Y,\tau\right]$. An example of such $X(t)$ is shown
in Figure~\ref{Dfig} (right).

In the case of a survival outcome, right censoring is common. For right censoring, \cite{Enc} proposed the artificial censoring estimator. A slight adaptation of this estimator is presented in Section~\ref{cens}.

\subsection{Mimicking counterfactual survival outcomes: assumptions and result}
\label{tresult}
%\subsubsection{Regularity conditions}\label{tcondc}

This section presents precise conditions under which $X(t)$ mimics
$Y^{(t)}$, for survival outcomes, following Section~\ref{simpler}
(Web-Appendix~\ref{XTapp} provides conditions similar to
Section~\ref{result}). We choose versions of
$F_{Y^{(t+h)}|\overline{Z}_t}$ that (a) are consistent with the fact
that treatment after death is irrelevant, and (b) satisfy all
regularity conditions below.  These versions are used in the
definition of $D$ for survival outcomes, and everywhere in the proof.

\begin{assu}\emph{(Regularity conditions)}.
\label{tassu}
\begin{itemize}
\item \emph{(support)}. There exists a finite number $y_2\geq \tau$ such that
\begin{enumerate}[a)]
\item If $Y>t$, all $F_{Y^{(t+h)}|\overline{Z}_t}$, for $h\geq 0$ and
  $t\in\left[0,\tau\right]$, have support $\left[t,y_2\right]$.
\item If $Y>t$, all $F_{Y^{(t+h)}|\overline{Z}_t}$, for $h\geq 0$ and
  $t\in\left[0,\tau\right]$, have a continuous non-zero density
  $f_{Y^{(t+h)}|\overline{Z}_t}\left(y\right)$ on
  $y\in\left[t+h,y_2\right]$.
\item There
exists an $\varepsilon>0$ such that for all $\omega\in\Omega$ and
$t$ with $Y>t$,
$f_{Y^{(t)}|\overline{Z}_t}\left(y\right)>\varepsilon$ for $y\in\left[t,y_2\right]$.
\end{enumerate}
\item \emph{(smoothness)}.
For every $\omega\in\Omega$
\begin{enumerate}[a)]
\item If $Z$ does not jump in $\left(t_1,t_2\right)$ and $Y>t_1$, the
restriction of $\left(y,t,h\right)\rightarrow
F_{Y^{(t+h)}|\overline{Z}_t}\left(y\right)$ to
$\left\{\left(y,t,h\right)\in\left[t_1,y_2\right]\times\left[t_1,t_2\right)\times
{\mathbb R}_{\geq 0}:y\geq t+h\right\}$ is $C^1$ in $\left(y,t,h\right)$.
\item The derivatives of $F_{Y^{(t+h)}|\overline{Z}_t}\left(y\right)$
($y> t+h$) with respect to $y$ and $h$ are bounded by constants $C_1$ and $C_2$,
respectively.
\item $\frac{\partial}{\partial y} F_{Y^{(t)}|\overline{Z}_t}(y)$ and
$\left.\frac{\partial}{\partial h}\right|_{h=0}
F_{Y^{(t+h)}|\overline{Z}_t}(y)$ ($y>t$) have derivatives with respect
to $y$ which are bounded by constants $L_1$ and $L_2$, respectively.
\item For all $\omega\in\Omega$ and $t$ with $Y>t$,
$F_{Y|\overline{Z}_t}(y)$ is continuous and strictly increasing on its
support $[t,y_2]$.
\end{enumerate}
\end{itemize}
\end{assu}

\begin{thm} \label{tthm1}
Suppose that Regularity Condition~\ref{tassu} is
satisfied. Then $D\left(y,t;\overline{Z}_t\right)$ as defined in
equation~(\ref{tdefd}) exists. Furthermore, for every $\omega\in\Omega$
there exists exactly one continuous solution $X(t)$ to
$dX(t)/dt=D\left(X(t),t;\overline{Z}_t\right)$ with final condition
$X\left(\tau\right)=Y$. If also Assumptions~\ref{inst}, \ref{ttcons}
and~\ref{ttinst} (consistency and no instantaneous treatment effect at
time of death) are satisfied then this $X(t)$ has the same
distribution as $Y^{(t)}$ given $\overline{Z}_t$ for all
$t\in\left[0,\tau\right]$.
\end{thm}

\subsection{Outline of the proof}

The proof of Theorem~\ref{tthm1} follows the same lines as the proof
of Theorem~\ref{thm1}.  The one essential difference between survival
outcomes and non-survival outcomes is: if $\overline{Z}_t$ indicates
the person is alive at time $t$, $X(t)$ should be greater than $t$,
since we want $X(t)$ to have the same distribution as $Y^{(t)}$ given
$\overline{Z}_t$ ($Y^{(t)}>t$ in that case because of Consistency
Assumption~\ref{ttcons}). This leads to an additional problem in the
proof for the continuous-time case, namely: how to prove that the
solution stays above the line $y=t$ for $t\in\left[0,Y\right]$? I
solve this additional problem in Web-Appendix~\ref{XTapp} by showing
that, under the assumptions of Section~\ref{tresult},
$D\left(t,t;\overline{Z}_t\right)\leq 1$. In addition, extra technical
problems arise because the smoothness conditions have to be adapted to
the survival setting; see Web-Appendix~\ref{XTapp} for details.

\subsection{Survival outcomes and right censoring}\label{cens}

In the case of a survival outcome, right censoring is common. \cite{Enc} proposed the artificial censoring estimator for
  administrative censoring. That is censoring due to end-of-follow-up because
  the study ends.
  The idea behind artificial censoring is that,
  instead of adding $X(t)$ or $X(0)$ to the model for predicting treatment
  changes (see Theorem~5.2), %\ref{thec},
  one could add a function $\tilde{X}(0)$ of $X(0)$ and the
  censoring time $C$, which is observed for all patients. The artificial censoring estimator treats the censoring time $C$ as a baseline covariate. This is
  justified in the case of censoring due to study closure, because in this case $C$ only
  depends on the date a patient enrolled in the study. Conditional on the value of $\overline{Z}_{t-}$, functions of $X(t)$ and
  $\overline{Z}_{t-}$ are not predictive of treatment changes (Theorem~5.2). %~\ref{thec}
  Therefore, conditional on $\overline{Z}_{t-}$, $\tilde{X}(0)$ is not predictive of treatment changes either. This produces an
  estimation procedure for $\psi$ analogous to that in Theorem~5.2, but that allows for right censoring.% \ref{thec}.

We slightly adapt this procedure, and propose to add a function of
$X(t)$ and $C$ to the model for the prediction of treatment changes.
%, which comes down to less artificial censoring and therefore most likely increased precision.
In particular, for $D$ as in equation~(3) %(\ref{Dexa1}),
and for $min(Y,C)\geq t$, we propose to add to the prediction model of treatment changes the function $\tilde{X}(t,\psi)=\min\left(X_\psi(t),C(t,\psi)\right)$, with
\begin{equation*}
C(t,\psi)=\left\{\begin{array}{ll} C & {\rm if}\;\psi\geq 0\\
t +e^\psi(C-t) & {\rm if}\;\psi<0.
\end{array}\right.
\end{equation*}
As required, $\tilde{X}(t,\psi)$ is a function of $X_\psi(t)$ and $C$. In addition, we will show that both for the case that $\psi\geq 0$ and for the case that $\psi< 0$, $\tilde{X}(t,\psi)$ is
observed for all patients. This follows from the fact that
\begin{equation}\label{Xstar}
\tilde{X}(t,\psi)=min(X^*(t,\psi),C(t,\psi)),\; {\rm with}\; X^*(t,\psi)=t+\int_t^{min(Y,C)}e^{\psi 1_{{\rm no}\;{\rm prophylaxis}\;{\rm at}\; s}}ds,
\end{equation}
which is observed for all patients. For $\psi\geq 0$, equation~(\ref{Xstar}) follows from
\begin{eqnarray*}
\tilde{X}(t,\psi)&=& min(t+\int_t^Ye^{\psi 1_{{\rm no}\;{\rm prophylaxis}\;{\rm at}\; s}}ds,C)\\
&=&min(t+\int_t^Ye^{\psi 1_{{\rm no}\;{\rm prophylaxis}\;{\rm at}\; s}}ds,t+\int_t^Ce^{\psi 1_{{\rm no}\;{\rm prophylaxis}\;{\rm at}\; s}}ds,C)\\
&=&min(X^*(t,\psi),C(t,\psi)),
\end{eqnarray*}
where for the second equality we used that for $\psi\geq 0$, $t+\int_t^Ce^{\psi 1_{{\rm no}\;{\rm prophylaxis}\;{\rm at}\; s}}ds\geq C$. For $\psi<0$, equation~(\ref{Xstar}) follows from
\begin{eqnarray*}
\tilde{X}(t,\psi)&=& min(t+\int_t^Ye^{\psi 1_{{\rm no}\;{\rm prophylaxis}\;{\rm at}\; s}}ds,t+e^\psi(C-t))\\
&=&min(t+\int_t^Ye^{\psi 1_{{\rm no}\;{\rm prophylaxis}\;{\rm at}\; s}}ds,t+\int_t^Ce^{\psi 1_{{\rm no}\;{\rm prophylaxis}\;{\rm at}\; s}}ds,t+e^\psi(C-t))\\
&=&min(X^*(t,\psi),C(t,\psi)),
\end{eqnarray*}
where for the second equality we used that for $\psi<0$, $t+\int_t^Ce^{\psi 1_{{\rm no}\;{\rm prophylaxis}\;{\rm at}\; s}}ds\geq t+e^\psi(C-t)$.

For the case
that $\psi<0$, some patients are ``artificially'' censored, since if
$C>t$, $C(t,\psi)=t+e^{\psi}(C-t)<C$. Artificial censoring produces a
subclass of the estimators considered in Theorem~\ref{thec}
allowing $h_t$ to depend on $\psi$:
$h_{t,\psi}(X_\psi(t),\overline{Z}_{t-})=1_{min(Y,C)\geq
  t}\tilde{h}_t(min\left(X_\psi(t),C(t,\psi)\right),\overline{Z}_{t-})$
(notice that $1_{min(Y,C)\geq t}$ is a function of
$\overline{Z}_{t-}$). In general, one could add to the prediction model for treatment changes any function of
$X_\psi(t)$ and $C$ that is observed for all patients. \cite{Enc} suggests to also consider adding
$\Delta(t,\psi)=1_{\tilde{X}(t,\psi)\leq C(t,\psi)}$ to the model for the prediction of treatment changes. Since both $\tilde{X}(t,\psi)$ and $C(t,\psi)$ are observed for all patients, so is $\Delta(t,\psi)$. Thus, the above reasoning shows that this procedure leads to consistent estimation of the treatment effect as well.

The procedure above can easily be adapted to for example model~(5), %\ref{Dexa2},
by replacing $C(t,\psi)$ accordingly. To be more specific, for that case one could use
\begin{equation*}
C(t,\psi)=t+e^{mim(\psi_1,0)+min(\psi_2,0)+min(\psi_3,0)}(C-t).
\end{equation*}

\section{Simulation study}\label{sim}

In the simulation study, we calibrated the distributions of the variables and the parameter values to HIV/AIDS data, perhaps the most salient example of application of structural nested models in the empirical literature. We focus on the first two years since HIV diagnosis. Time zero is the time of HIV diagnosis. The outcome variable is the CD4 count, a commonly used marker of the state of the immune system of HIV-positive patients. The usual treatment for HIV-positive patients is ART, antiretroviral treatment. ART is not always initiated immediately after diagnosis. ART initiation time often depends on the last measured CD4 count. When the CD4 count is at or below $350$ copies/ml, HIV-positive patients are much more likely to initiate ART than when the CD4 count is above $350$ copies/ml. Web-Appendix~\ref{simapp} describes how we generated the data for the simulation study in detail, including distributions and parameter values. This section provides an overview.

In this simulation study no one is treated at time zero, and once treatment is initiated, it is never stopped.
$Y^{(t)}$ is the counterfactual outcome had treatment been as given in reality until time $t$, and continued or initiated after that.
For example, if treatment was initiated by time $t$ for a particular patient, $Y^{(t)}$ is the observed outcome for that patient,
since he or she was already treated at time $t$ and treatment is never stopped. On the other hand, if treatment was not initiated by time $t$,
$Y^{(t)}$ is the outcome had treatment been initiated at time $t$. Thus, in the definition of $Y^{(t)}$ in Section~\ref{sets},
the switch at time $t$ to ``some kind of baseline treatment regime $\overline{0}$'' is, in this case,
``treat continuously'' from time $t$ onwards. In the simulations, we study a setting with $t\in[0,2]$.
The subscript $_{t}$ indicates the treatment initiation time, so for example $L_{1,t}$ indicates (counterfactual) covariates at time $1$
under ``treatment started at time $t$''. Similarly, the subscript $_{\infty}$ indicates (counterfactual) variables under no treatment.
For example, $L_{2,\infty}$ indicates (counterfactual) covariates at time $2$ under no treatment. In the simulation design, the counterfactual covariates $L$ are as follows:
\begin{eqnarray*}
L_0&=&\tilde{L}_0+e_0,\\
L_{1,\infty}&=&\tilde{L}_0-\beta_0+e_{1,\infty},\\
L_{2,\infty}&=&\tilde{L}_0-2\beta_0+e_{2,\infty}\\
L_{1,t}&=&\tilde{L}_0-\beta_0+\theta (1-t)+e_{1,t} \;{\rm for}\;t\in\left[0,1\right], \;{\rm and}\; L_{1,\infty}\;{\rm otherwise}\\
L_{2,t}&=&\tilde{L}_0-2\beta_0+\psi (2-t)+e_{2,t},
\end{eqnarray*}
 where $\tilde{L}_0$ and the $e_{j,t}$ are random variables with values in $\mathbb{R}$. Notice that $(1-t)$ and $(2-t)$ are simply the durations of treatment until the respective
 covariate measurements. We assume that the $e_{j,t}$ ($j=0,1,2$) are independent of $\tilde{L}_0$, and that the $e_{2,t}$ have a distribution function which does not depend on $t$.
 We also assume that the $e_{2,t}$ are independent of all previous variables (and of the treatment initiation time, $T$, described below).
 In the simulations, $\psi\geq 0$ (a similar study could have been done for $\psi<0$). We define $Y_{t}=L_{2,t}$, the counterfactual outcome with treatment initiated at time $t$, which could potentially be observed at time $2$.

We show in Web-Appendix~\ref{simapp} that the outcome processes adopted in our simulation study are not rank preserving. This is easily seen because
with probability one, two patients with the same observed data do not have the same value of $\tilde{L}_0$.

Suppose that the hazard of the treatment initiation time, $T$, given the covariate history at time $t$ and given that treatment was not initiated before time $t$, is piecewise constant as follows:
\begin{equation*}
\lambda_T(t)=\left\{\begin{array}{ll}
\lambda^{(0)}_0 & {\rm if}\; L_0> c_0 \;{\rm and}\;t\in[0,1]\\
\lambda^{(0)}_1 & {\rm if}\; L_0\leq c_0 \;{\rm and}\;t\in[0,1]\\
\lambda^{(1)}_0 & {\rm if}\; L_{1,\infty} > c_1 \;{\rm and}\;t\in(1,2]\\
\lambda^{(1)}_1 & {\rm if}\; L_{1,\infty} \leq c_1 \;{\rm and}\;t\in(1,2],
\end{array}\right.
\end{equation*}
for constants $c_0$ and $c_1$ in $\mathbb{R}$. Notice that $T$ depends on $\tilde{L}_0$, $e_{0,\infty}$, and $e_{1,\infty}$, if $\lambda^{(0)}_0\neq \lambda^{(0)}_1$ or $\lambda^{(1)}_0\neq \lambda^{(1)}_1$.

In the simulation study, treatment can be initiated in continuous time, but the covariates are only measured at times $0$, $1$, and $2$, so that the treatment and covariate
history up to time $t$, $\overline{Z}_t$, consists of the treatment information up to time $t$ and $L_0$, $(L_0,L_1)$, or $(L_0,L_1,L_2)$,
depending on whether $t\in[0,1)$, $t\in[1,2)$, or $t=2$. In the simulations, treatment affects later outcomes, and time-dependent covariates ($L_1$) which
depend on previous treatment also predict future treatment and the outcome of interest. This is the type of setting structural nested models were developed for.

Web-Appendix~\ref{simapp} shows that for this data generating mechanism,
\begin{equation*}
D(y,t;\overline{Z}_t)=-\psi 1_{{\rm untreated}\;{\rm at}\;t}.
\end{equation*}
Then, it follows from the definition of $X_\psi$ that
\begin{equation*}
X_{\psi}(t)=Y+\psi (min(T,2)-t)1_{T>t},
\end{equation*}
where $(min(T,2)-t)1_{T>t}$ is the duration of the patient not being on treatment between time $t$ and time $2$.

As shown in Web-Appendix~\ref{simapp}, a consistent estimator of $\psi$ can be defined as follows. In the first step, the nuisance parameters
$\left(\lambda^{(0)}_0,\lambda^{(0)}_1,\lambda^{(1)}_0,\lambda^{(1)}_1\right)$ are estimated using maximum likelihood theory.
In the second step, $\psi$ is estimated as $\hat{\psi}=- \sum_{i=1}^n A_{1i}/ \sum_{i=1}^n A_{2i}$, where
\begin{eqnarray*}
A_{1i}&=&-Y_i\left(Z_i(0)\hat{\lambda}^{(0)}_{1}+(1-Z_i(0))\hat{\lambda}^{(0)}_{0}\right)min(T_i,1)\nonumber\\
&&-Y_i\left(Z_i(1)\hat{\lambda}^{(1)}_{1}+(1-Z_i(1))\hat{\lambda}^{(1)}_{0}\right)\left(1-\delta_i^{(0)}\right)(min(T_i,2)-1)\nonumber\\
&&+Y_i\delta_i^{(0)}+ Y_i\delta_i^{(1)},
\end{eqnarray*}
\begin{eqnarray*}
A_{2i}&=&-\left(Z_i(0)\hat{\lambda}^{(0)}_{1}+(1-Z_i(0))\hat{\lambda}^{(0)}_{0}\right)min(T_i,1)min(T_i,2)\nonumber\\
&&-\left(Z_i(1)\hat{\lambda}^{(1)}_{1}+(1-Z_i(1))\hat{\lambda}^{(1)}_{0}\right)\left(1-\delta_i^{(0)}\right)(min(T_i,2)-1)^2\nonumber\\
&&+\delta_i^{(0)}min(T_i,2) + \delta_i^{(1)}(min(T_i,2)-1),
\end{eqnarray*}
$\delta_i^{(0)}=1_{T_i\leq 1}$, $\delta_i^{(1)}=1_{1<T_i\leq 2}$, $Z_i(0)=1_{L_0\leq c_0}$, and $Z_i(1)=1_{L_1\leq c_1}$.

We ran a simulation study with $n=500$, $1000$, $2000$, $5000$, and $10000$, with 5000 repetitions each. The results are presented in Table~1. As detailed in Web-Appendix~\ref{simapp}, setting~1 has the least noise around the signals, and setting~3 the most.

\begin{table}[htb!]\caption{Simulations. Mean Squared Errors (MSE) and bias. $5000$ repetitions each.}
\begin{tabular}{|l|lll|lll|lll|}
&setting 1 &&& setting 2 &&& setting 3\\
n&MSE & $\mathlarger{\frac{{\rm MSE}\times n}{1000}}$ &bias & MSE & $\mathlarger{\frac{{\rm MSE}\times n}{1000}}$ & bias & MSE & $\mathlarger{\frac{{\rm MSE}\times n}{1000}}$ & bias\\
\hline
&&&&&&&&&\\
100&747&75&-0.28&1907&191&-0.040&2875&287&-0.39\\
500&146&73&-0.22&356&178&-0.29&542&271&-0.61\\
1000&72&72&-0.10&176&176&-0.068&268&268&-0.23\\
2000&35&70&-0.11&89&179&0.051&138&275&-0.08\\
5000&14&69&-0.067&35&175&-0.0040&54&268&-0.05\\
10000&6.6&66&-0.066&18&178&-0.022&27&270&-0.06\\
\end{tabular}
\end{table}

In this simulation study, both for small and large samples, the bias of the estimators is small. In all three settings and for all sample sizes considered
(including the small sample size $n=100$), the MSE of the estimators arises mostly from the variance, not from the bias. Also, if the true parameter $\psi$ equals
$300$ as in this simulation study, for $n=500$, $\sqrt{MSE}/\psi=0.04$ in setting~1, and~$0.08$ in setting~3. Thus, the estimates are already precise in
relatively small samples. Because, as shown in the Web-Appendix, the MSE in this simulation study does not depend on the true parameter, $\psi$, a larger sample
size would be required to obtain precise estimators of small true parameter values $\psi$. We conclude that in this simulation study, continuous-time structural
nested models perform extremely well.

\section{Discussion}
\label{ext}

Structural nested models have become a major part of statistical tools
for estimation of the effect of time varying treatments, in the
presence of time-dependent confounding by indication; see
e.g.\ \cite{hyp} for a discrete-time application, and
e.g.\ \cite{Enc}, \cite{Aids}, \cite{smoke}, \cite{RGPCP}, \cite{R92},
\cite{Tilling}, and \cite{Keiding, Kei} for continuous-time
applications. Structural nested models in continuous time are useful
to estimate the effect of a treatment that can be initiated at any
point in time, and for which a short duration of treatment has a small
effect on the outcome of interest. In contrast with discrete-time
structural nested models, in the case of survival outcomes, the
resulting parameter estimates can often be interpreted as
rates. So-far, continuous-time analyses relied on (local) rank
preservation. The main result of the current article is to prove that
for continuous-time structural nested models, assumptions about the
joint distributions of counterfactuals or deterministic treatment
effects/ (local) rank preservation are not necessary to ``mimic
counterfactual outcomes'', and, based on that, to consistently
estimate treatment effects. This article provides a proof for outcomes
that are measured at the end of the study as well as a proof for
survival outcomes. Important public health decisions are based on
analyses with continuous-time structural nested models, so it is
important to relax unverifiable and disputable assumptions underlying
these analyses.
%should be done right People base

An interesting topic for future research is to investigate whether the
Support Conditions~\ref{sup} or~\ref{tassu} can be weakened, for
example to an assumption about the support varying in a differentiable
way between the jump times of the covariate- and treatment process
$Z$. We expect that in that case one has to assume that where $Z$ jumps,
the support of $Y^{(t)}$ given $\overline{Z}_t$ gets smaller or stays
the same as $t$ increases (see Figure~\ref{supfig}). Otherwise, $X(t)$
may move out of the support of $Y^{(t)}$ given $\overline{Z}_t$
(recall that $X$ is the solution to a differential equation with
\emph{final} condition). It is reasonable to assume that the support
of $Y^{(t)}$ given $\overline{Z}_t$ gets smaller or stays the same as
$t$ increases, since more information about $\overline{Z}$ should not
enlarge the range of $Y^{(t)}$.

\begin{figure}[htb!]
\begin{picture}(380,50)
\put(43,45){\makebox(0,0){support}}
\put(43,35){\makebox(0,0){of $Y^{(t)}$}}
\put(43,25){\makebox(0,0){given $\overline{Z}_t$}}
\put(65,5){\vector(1,0){315}}
\put(65,5){\line(0,1){45}}
\put(373,-1){\makebox(0,0){$t$}}
\put(185,-2){\makebox(0,0){jump time of $Z$}}
\put(185,3){\line(0,1){4}}
\qbezier(65,45)(135,45)(185,40)
\qbezier(185,35)(290,26)(380,27)
\qbezier(65,7)(135,10)(185,10)
\qbezier(185,17)(290,18)(380,23)
\end{picture}
\caption{Example of support of $Y^{(t)}$ given $\overline{Z}_t$.}
\label{supfig}
\end{figure}

A problem which may occur without a support condition is that the
denominator in equation~(\ref{Deq}) (the quotient expression for $D$)
or in equation (\ref{Dneq2}) or (\ref{Dneq}) (the quotient expression
for $D^{(n)}$) may tend to $0$, which may ``blow up'' $D$ or
$D^{(n)}$. In that case it might help to assume that there exists a
constant $C$ such that (a) for all $\omega\in\Omega$, $t$ and $y$, $
F_{Y^{(t+h)}|\overline{Z}_t}^{-1}\circ
F_{Y^{(t)}|\overline{Z}_t}\left(y\right)-y\leq C\cdot h $, and (b) for all $t$, $y$ and $B\subset \overline{Z}_t$ with
$P\left(\overline{Z}_t\in B\right)>0$,
$
F_{Y^{(t+h)}|\overline{Z}_t\in B}^{-1}\circ
F_{Y^{(t)}|\overline{Z}_t\in B}\left(y\right)-y\leq C\cdot
h$.
%That is, quantiles of the outcome of interest do not move more than
%$C\cdot h$ when treatment is stopped at time $t+h$ instead of $t$,
%both given that the treatment- and covariate process is $\overline{Z}_t$ and
%given that it satisfies $\overline{Z}_t\in B$ for some $B$ with
%$P\left(\overline{Z}_t\in B\right)>0$.
This assumption does not look unreasonable if there is no
``instantaneous treatment effect''. It is to be expected that under
this assumption both $D$ and $D^{(n)}$ are bounded by $C$.
%The problem of a denominator
%which might tend to $0$, blowing up $D$ or $D^{(n)}$, would be solved
%this way.
%Nog netter maken dit, maar lijkt een leuk idee!

Based on the results of the current article, \cite{Lok} shows that also if a semiparametric Cox model is used to predict treatment changes in Theorem~\ref{thec}, the resulting estimating equations for the treatment effect are unbiased. However, the estimating equations are no longer of the form of an average of terms that are independent for the different persons. Thus, consistency and asymptotic normality for this situation constitute interesting topics for future research.

\section{Acknowledgements}

I am indebted to Richard Gill and
  Aad van der Vaart for their support, insight and encouragement on
  this project. I also thank James Robins for fruitful discussions,
  and Nell Sedransk for constructive comments on the writing. I thank Susan Little for allowing me to use the AIEDRP data to calibrate the distributions in the simulation study. This
  work was sponsored by the Netherlands Organization for Scientific
  Research (NWO) with a Talent scholarship, and by the National Institutes of Health, NIAID R01AI100762.
  The content is solely the responsibility of the author and does not necessarily represent the official views of the National Institutes of Health.

\pagebreak

\appendix

\noindent Appendix~A contains results that are frequently used in the main
article. Appendix~B describes Mimicking counterfactual survival outcomes. Appendix~C describes details of the simulation study. Appendix~D describes some facts about conditioning.
Appendix~E describes a corollary of the Local Inverse Function Theorem. Appendix~F describes some facts about Lipschitz continuity and differentiability. Appendix~G describes some theory about differential equations. Appendix~H describes convergence theorems.

\section{Results that are frequently used}
%
%\section{Some facts about conditioning and differential equations}\label{Cond}
%
The first theorem is a corollary of a theorem in \cite{anDe} Chapter 2, see Web-Appendix~\ref{difeqapp}.
%It is a consequence of Gronwall's lemma.
\begin{thm}
\label{difbc2} Suppose that $I$ is a closed interval in ${\mathbb R}$,
$f:I\times \left[y_1,y_2\right]\rightarrow {\mathbb R}$ is continuous
with for all $t\in I$, $f\left(t,y_1\right)=f\left(t,y_2\right)=0$ and
$C:I\rightarrow\left[0,\infty\right)$ is continuous, and suppose that
\begin{equation}
\left| f\left(t,y\right)-f\left(t,z\right)\right|\leq
C\left(t\right) \left| y-z\right|
\label{Lipvw}
\end{equation}
for all $t\in I$ and $y,z\in \left[y_1,y_2\right]$.  Then, for every
$t_0\in I$ and $y_0\in\left[y_1,y_2\right]$, there exists a unique solution
$y\left(t\right)$ of $y'\left(t\right)=f\left(t,y(t)\right)$ with
$y\left(t_0\right)=y_0$, and this solution is defined for all $t\in
I$. Furthermore, $y\left(t\right)\in\left[y_1,y_2\right]$ for all $t\in
I$. Suppose that $g:I\times \left[y_1,y_2\right]\rightarrow {\mathbb
R}$ is continuous and
$z:I\rightarrow \left[y_1,y_2\right]$ is a
solution of
$z'\left(t\right)=g\left(t,z(t)\right)$. Then
\begin{eqnarray*} \lefteqn{\left| y\left(t\right)
-z\left(t\right)\right|}\\
&&\hspace{0.7cm}\leq
e^{\int_{t}^{t_0} C\left(s\right)ds}
\left| y\left(t_0\right)-z\left(t_0\right)\right|
+\int_{t}^{t_0} e^{\int_{t}^s C\left(\eta\right)d\eta}
\left| f\left(t,z\left(t\right)\right)-
g\left(s,z\left(s\right)\right)\right| ds
\end{eqnarray*}
for all $t,t_0\in I$ with $t\leq t_0$.
\end{thm}

The proof of the following lemma can be found in Appendix~\ref{Cond}.
\begin{lem} \label{Zt-} Let $X$ be a random variable with
$E\left|X\right|<\infty$, and let $\overline{Z}_t$ be a random variable with
values in $\overline{{\cal Z}}_t$, the space of cadlag functions on
$\left[0,t\right]$ provided with the projection $\sigma$-algebra, with
$P\left(Z \;{\rm jumps}\;{\rm at}\;t\right)=0$. Then any version of
$E\bigl[X|\sigma\bigl(\cup_{n=1}^{\infty}\overline{Z}^{(n)}_t\bigr)\bigr]$,
with $\overline{Z}^{(n)}_t$ as defined in Section~\ref{discsch}, is
also a version of $E\left[X|\overline{Z}_t\right]$.
\end{lem}

%\appendix\chapter{}

%\begin{center}
%{\huge Web-Appendix with Mimicking counterfactual outcomes to estimate causal effects}\vspace*{0.2cm}\\
%
%{{\Large Judith J.~Lok}\vspace*{0.2cm}\\
%{\it Department of Biostatistics, Harvard School of Public Health}}
%\end{center}

\section{Web-Appendix: Mimicking counterfactual survival outcomes}
\label{XTapp}

\subsection{Introduction}
\label{tintroapp}

For the definition of the infinitesimal shift function $D$ for
survival outcomes see Section~\ref{tintro}. Also the additional
Consistency Assumption~\ref{ttcons} and Assumption of No instantaneous
treatment effect \ref{ttinst} can be found in Section~\ref{tintro}.
Assumption~\ref{ttcons} implies the obvious fact that
treatment in the future does not cause or prevent death at present:

\begin{lem} \label{conslem} Under Assumptions~\ref{ttcons} (consistency),
%and~\ref{ttinst} (no instantaneous effect of treatment at the time the
%person died),
\begin{enumerate}[a)]
\item For all $h\geq 0$: $Y^{(t+h)}=Y$ on $\left\{\omega:Y\leq t\right\}\cup\cup_{h\geq
0}\left\{\omega:Y^{(t+h)}\leq t\right\}$.
\item For all $\left(y,t,h\right)$ with $y\leq t+h$ and $h\geq 0$:
$\left\{\omega: Y^{(t+h)}\leq y\right\}=\left\{\omega: Y\leq
y\right\}$.
\end{enumerate}
\end{lem}

\noindent{\bf Proof.} a): From Assumption~\ref{ttcons}, $Y^{(t)}=Y$ on
$\left\{\omega:Y\leq t\right\}\cup\left\{\omega:Y^{(t)}\leq
t\right\}$. Thus if $Y\leq t$ then $Y^{(t)}=Y$, and moreover for all
$h>0$, $Y<t+h$, so that, again from Assumption~\ref{ttcons},
$Y^{(t+h)}=Y$. If $Y^{(t)}\leq t$ the same reasoning can be used to prove that
$Y^{(t+h)}=Y$. If for some $h>0$, $Y^{(t+h)}\leq t$ then
also $Y^{(t+h)}\leq t+h$, so that, again from Assumption~\ref{ttcons},
$Y^{(t+h)}=Y$.\\ b):
For $y\leq t+h$ and $h\geq 0$, $\left\{\omega: Y^{(t+h)}\leq
y\right\}=\left\{\omega: Y^{(t+h)}\leq
y\right\}\cap\left\{Y^{(t+h)}\leq t+h\right\}=\left\{\omega:Y\leq
y\right\}$ because of Assumption~\ref{ttcons}.  \hfill $\Box$

$ $
%c) is a special case of b).

\noindent Henceforth, this article will only use versions of conditional
distributions which are consistent with Lemma~\ref{conslem} in the
sense that $F_{Y^{(t+h)}|\overline{Z}_t}\left(y\right)=F_{Y|\overline{Z}_t}\left(y\right)$ for all $y\leq t+h$, $h\geq 0$, and
$\omega\in\Omega$.
%was: are ``natural'' considering Lemma.. in the sense the e.g.... I
%checked that this can be
%``in the sense that'' instead of ``for example''

Notice that with definition~(\ref{tdefd}) of $D$, the area where $D$ is
possibly non-zero is $\left(y,t\right)\in\left[0,\infty\right)\times
  \left[0,\min\left\{Y,\tau\right\}\right]:y\geq t$. Therefore if
  $Y<\tau$, the solution to the differential equation $X(t)$ is equal
  to $Y$ for $t\in\left[Y,\tau\right]$. An example of such $X(t)$ is
  shown in Figure~\ref{Dfig}, right panel.

\subsection{Mimicking counterfactual survival outcomes: assumptions and result}
\label{tresultapp}
%\subsubsection{Regularity conditions}\label{tcondc}

This section provides precise conditions under which $X(t)$ mimics $Y^{(t)}$ in the
sense that $X(t)$ has the same distribution as $Y^{(t)}$ given
$\overline{Z}_t$, for survival outcomes. We choose versions of
$F_{Y^{(t+h)}|\overline{Z}_t}$ which are consistent with
Lemma~\ref{conslem} and which satisfy all regularity conditions
below. These versions are used in the definition of $D$ of
equation~(\ref{tdefd}), and everywhere else in this section.

%As already mentioned in Section~\ref{XTo} we assume that $Z(t)$
%includes an indicator for whether or not a person is alive at time
%$t$. $\overline{Z}_t$ then automatically includes the time of death of
%a person in case the person died at or before time $t$. Among other
%things this implies that $\overline{Z}_t$ jumps at the time of death
%$Y$, and therefore $D\left(y,t;\overline{Z}_t\right)$ does not
%necessarily exist at $t=Y$. As always when $Z$ jumps at time $t$ and we
%consider the differential equation (\ref{tXdef}) on the interval below
%$t$ we use the limit from the left for $D$ there.

Section~\ref{tintro} indicated why it is not reasonable
to assume that the conditional distribution of the survival time has
the fixed support $\left[y_1,y_2\right]$ given any covariate- and
treatment history $\overline{Z}_t$. If a person is alive at time $t$,
one often expects that $t$ is the left limit of the
support. Therefore, this article assumes that

\begin{assu}\label{tsup}\emph{(support)}.
There exists a finite number $y_2\geq\tau$ such that
\begin{enumerate}[a)]
\item For all $\omega\in\Omega$ and $t$ with $Y>t$, all
$F_{Y^{(t+h)}|\overline{Z}_t}$ for $h\geq 0$ have support $\left[t,y_2\right]$.
\item For all $\omega\in\Omega$ and $t$ with $Y>t$, all
$F_{Y^{(t+h)}|\overline{Z}_t}$ for $h\geq 0$ have a continuous non-zero
density $f_{Y^{(t+h)}|\overline{Z}_t}\left(y\right)$ on
$y\in\left[t+h,y_2\right]$.
\item There exists an $\varepsilon>0$ such that for all
$\omega\in\Omega$ and $t$ with $Y>t$,
$f_{Y^{(t)}|\overline{Z}_t}\left(y\right)>\varepsilon$ for
$y\in\left[t,y_2\right]$.
\end{enumerate}
\end{assu}

Next consider the smoothness conditions in Section~\ref{XY}. It does
not seem reasonable to assume that
$F_{Y^{(t+h)}|\overline{Z}_t}\left(y\right)$ is continuously
differentiable with respect to $h$ and $y$ on
$\left(h,y\right)\in\left[0,\infty\right)\times \left[t,y_2\right]$
since for $y\leq t+h$,
$F_{Y^{(t+h)}|\overline{Z}_t}\left(y\right)=F_{Y|\overline{Z}_t}\left(y\right)$
(Lemma~\ref{conslem}b). Thus, the derivative of
$F_{Y^{(t+h)}|\overline{Z}_t}\left(y\right)$ with respect to $h$ is
likely not to exist at $y=t+h$ (and is equal to zero for
$y<t+h$). Also the derivative of
$F_{Y^{(t+h)}|\overline{Z}_t}\left(y\right)$ with respect to $y$ may
not exist at $y=t+h$, because of the different treatment before and
after $t+h$. For survival outcomes, the smoothness
conditions~\ref{cont1}--\ref{Lipc} are therefore replaced by:

\begin{assu} \label{tcont1} \emph{(continuous derivatives)}.
For $\omega\in\Omega$ fixed,
\begin{enumerate}[a)]
\item If $Y>t$ then $F_{Y^{(t+h)}|\overline{Z}_t}\left(y\right)$
restricted to $\left\{\left(h,y\right)\in\left[0,\infty\right)\times
\left[t,y_2\right]: y\geq t+h\right\}$ is $C^1$ in $\left(h,y\right)$.
% and has a $C^1$ extension to
%$\left\{\left(h,y\right)\in\left[0,\infty\right)\times
%\left[t,y_2\right]: y\geq t+h\right\}$.
%\item If $Y>t$ then $F_{Y^{(t+h)}|\overline{Z}_t}\left(y\right)$ restricted to
%$\left\{\left(h,y\right)\in\left[0,\infty\right)\times
%\left[t,y_2\right]: y\geq t+h\right\}$ is continuous in
%$\left(h,y\right)$, and restricted to
%$\left\{\left(h,y\right)\in\left[0,\infty\right)\times
%\left[t,y_2\right]: y>t+h\right\}$ is differentiable with respect to
%$h$ and $y$ with derivatives which are continuous in
%$\left(h,y\right)$ and have a continuous extension to
%$\left\{\left(h,y\right)\in\left[0,\infty\right)\times
%\left[t,y_2\right]: y\geq t+h\right\}$.
\item If $Z$ does not jump in $\left(t_1,t_2\right)$ and $Y>t_1$ then both
$\left.\frac{\partial}{\partial h}\right|_{h=0} F_{Y^{(t+h)}|\overline{Z}_t}\left(y\right)$ and
$\frac{\partial}{\partial y}F_{Y^{(t)}|\overline{Z}_t}\left(y\right)$
are continuous in $\left(y,t\right)$ on
$\left\{\left(y,t\right)\in
\left[t_1,y_2\right]\times\left[t_1,t_2\right):y> t\right\}$
and can be continuously extended to
$\left\{\left(y,t\right)\in
\left[t_1,y_2\right]\times\left[t_1,t_2\right]:y\geq t\right\}$.
\end{enumerate}
\end{assu}

\begin{assu} \label{tbdd+} \emph{(bounded derivatives)}.
\begin{enumerate}[a)]
\item There exists a constant $C_1$ such that for all
 $t$, $h\geq 0$ and $y> t+h$, for $\omega\in\Omega$ with $Y>t$,
\begin{equation*} \frac{\partial}{\partial y}
F_{Y^{(t+h)}|\overline{Z}_t}\left(y\right)\leq C_1.
\end{equation*}
\item There exists a constant $C_2$ such that for all
 $t$, $h\geq 0$ and $y> t+h$, for $\omega\in\Omega$ with $Y>t$,
\begin{equation*} \Bigl|\frac{\partial}{\partial h}
F_{Y^{(t+h)}|\overline{Z}_t}\left(y\right)\Bigr|\leq C_2.
\end{equation*}
\end{enumerate}
\end{assu}

\begin{assu} \label{tLipc}
\emph{(Lipschitz continuity)}.
\begin{enumerate}[a)]
\item There exists a constant $L_1$ such that for all
$t$ and $y,z\in\left(t,y_2\right]$, for $\omega\in\Omega$ with $Y>t$,
\begin{equation*}\Bigl|\frac{\partial}{\partial y}F_{Y^{(t)}|\overline{Z}_t}\left(y\right)-
\frac{\partial}{\partial z}F_{Y^{(t)}|\overline{Z}_t}\left(z\right)\Bigr| \leq
L_1\left|y-z\right|. \end{equation*}
\item There exists a constant $L_2$ such that for all
$t$ and $y,z\in\left(t,y_2\right]$, for $\omega\in\Omega$ with $Y>t$,
\begin{equation*}\left|\left.\frac{\partial}{\partial h}\right|_{h=0} F_{Y^{(t+h)}|\overline{Z}_t}\left(y\right)-
\left.\frac{\partial}{\partial h}\right|_{h=0} F_{Y^{(t+h)}|\overline{Z}_t}\left(z\right)\right| \leq
L_2\left|y-z\right|.\end{equation*}
\end{enumerate}
\end{assu}
%\noindent The above smoothness conditions on
%$F_{Y^{(t+h)}|\overline{Z}_t}$ concentrate on $y\geq t+h$. For
%$y\in\left[t,t+h\right)$ we can choose
%$F_{Y^{(t+h)}|\overline{Z}_t}(y)=F_{Y|\overline{Z}_t}(y)$ because of
%Lemma~\ref{conslem}b. Because of Assumption~\ref{inst} (consistency),
%$F_{Y|\overline{Z}_t}$ has the same support as
%$F_{Y^{(\tau)}|\overline{Z}_t}$, so $F_{Y|\overline{Z}_t}$ has support
%$\left[t,y_2\right]$ if $Y>t$ under Assumption~\ref{tsup}a. We assume
%
\begin{assu}\label{Ytsup} \emph{(smoothness)}.
For all $\omega\in\Omega$ and $t$ with $Y>t$,
$F_{Y|\overline{Z}_t}(y)$ is continuous in $y$.
\end{assu}

\begin{thm} \label{tthm1app}
Suppose that Regularity Conditions~\ref{tsup}--\ref{tLipc} are
satisfied. Then $D\left(y,t;\overline{Z}_t\right)$ as defined in
equation~(\ref{tdefd}) exists. Furthermore for every $\omega\in\Omega$
there exists exactly one continuous solution $X(t)$ to
$dX(t)/dt=D\left(X(t),t;\overline{Z}_t\right)$ with final condition
$X\left(\tau\right)=Y$. If also Assumptions~\ref{inst}, \ref{ttcons}
and~\ref{ttinst} (consistency and no instantaneous treatment effect at
time of death) are satisfied then this $X(t)$ has the same
distribution as $Y^{(t)}$ given $\overline{Z}_t$ for all
$t\in\left[0,\tau\right]$.
\end{thm}

The simpler regularity conditions, comparable with Section~\ref{simpler} for non-survival outcomes, can be found in Section~\ref{tresult}.

\subsection{Existence of and a different expression for $D$}
\label{tDsec}

If $Y\leq t$, $D\left(y,t;\overline{Z}_t\right)=0$ by
definition~(\ref{tdefd}). Thus we can concentrate on $\omega\in
\Omega$ with $Y>t$. If $y>t$, Corollary~\ref{quotc} can be applied on
$F_h\left(y\right)=F_{Y^{(t+h)}|\overline{Z}_t}\left(y\right)$ with
$y_0=y$ and $U_{0,y_0}\cap\left\{h\geq
0\right\}=\left[0,y-t\right)\times \left(t,y_2\right]$, because of
Assumptions~\ref{tcont1}a
%, \ref{ttsmooth}a
and~\ref{tsup}c.
Thus for $y>t$, $D$ as defined in equation~(\ref{tdefd})
exists, and it is equal to
\begin{eqnarray} D\left(y,t;\overline{Z}_t\right)&=&\left(\left.\frac{\partial}{\partial h}\right|_{h=0}
F^{-1}_{Y^{(t+h)}|\overline{Z}_t}\right)\left(F_{Y^{(t)}|\overline{Z}_t}\left(y\right)\right)\nonumber\\
&=&-\frac{\left.\frac{\partial}{\partial h}\right|_{h=0}
F_{Y^{(t+h)}|\overline{Z}_t}\left(y\right)}{\frac{\partial}{\partial
y}F_{Y^{(t)}|\overline{Z}_t}\left(y\right)}.\label{tDeq}\end{eqnarray}
$D\left(t,t;\overline{Z}_t\right)$ is by definition (\ref{tdefd})
equal to the limit as $y\downarrow t$ of this
$D\left(y,t;\overline{Z}_t\right)$, which exists because of
Assumptions~\ref{tsup}c and~\ref{tcont1}a.
%ja

\subsection{Existence and uniqueness of $X(t)$} \label{tDexu}

If $\overline{Z}_t$ indicates the person is alive at time $t$ and $X(t)$
has the same distribution as $Y^{(t)}$ given $\overline{Z}_t$, we
should have that $X(t)$ stays above $t$ ($Y^{(t)}>t$ in that case
because of Consistency Assumption~\ref{ttcons}). In order to prove
that $X(t)$ stays indeed above $t$ if the person is alive at time $t$, we
prove that $D\left(t,t;\overline{Z}_t\right)\leq 1$.

\begin{lem} \label{lDl1} Under Assumptions~\ref{tsup} and~\ref{tcont1}a,
%and~\ref{ttsmooth}a,
if $Y>t$ then $D\left(t,t;\overline{Z}_t\right)\leq 1$.
\end{lem}

\noindent{\bf Proof.} We start with some ideas, which are made precise
below. Intuition says that $D$ not only measures the increase of
quantiles when treatment is prolonged but also the decrease of
quantiles when treatment is withheld. Thus quantiles $y$ seem to
approximately move to $y-h\,D\left(y,t;\overline{Z}_t\right)$ when
treatment is withheld between $t$ and $t+h$. If quantiles near $t$
move down to $t$ with speed greater than $1$ when treatment is
withheld starting from $h$ it seems like these quantiles will end up
below $t$ when treatment is withheld at time $t$. However, if treatment is
withheld starting from $t$ this does not cause death at or before time $t$,
so the quantiles above $t$ should stay above $t$. This leads to a
contradiction. The following makes this precise.
%although the intuitive idea may now be blurred by the technicalities.
%For quantiles just above $t$
%this suggests that this leads to $D$ being less than or
%equal to $1$ in the limit, since otherwise these quantiles seem to
%disappear below $t$.

Fix $\omega\in\Omega$ and fix $t$ for which $Y>t$. Recall from
Section~\ref{tDsec} that $a:=\lim_{y\downarrow
  t}D\left(y,t;\overline{Z}_t\right)$ exists. We need to prove that
$a\leq 1$. Suppose that $a>1$. It is shown that this leads to a
contradiction. Notice that because of the chain rule, for $y>t+h$,
\begin{eqnarray} \frac{\partial}{\partial h}\left(F_{Y^{(t)}|\overline{Z}_t}^{-1}\circ
F_{Y^{(t+h)}|\overline{Z}_t}\right) \left(y\right)
&=&\left(F_{Y^{(t)}|\overline{Z}_t}^{-1}\right)'\left(F_{Y^{(t+h)}|\overline{Z}_t}
\left(y\right)\right) \frac{\partial}{\partial h}
F_{Y^{(t+h)}|\overline{Z}_t} \left(y\right)\nonumber\\
&=&\frac{\frac{\partial}{\partial h} F_{Y^{(t+h)}|\overline{Z}_t}
\left(y\right)}{F'_{Y^{(t)}|\overline{Z}_t}\left(F_{Y^{(t)}|\overline{Z}_t}^{-1}\circ
F_{Y^{(t+h)}|\overline{Z}_t} \left(y\right)\right)}
\label{Dl1}\end{eqnarray} exists and is continuous in
$\left(h,y\right)$ for $y>t+h$ with a continuous extension to
\linebreak $\left\{\left(h,y\right)\in\left[0,\infty\right)\times
\left[t,y_2\right]: y\geq t+h\right\}$ because of
Assumptions~\ref{tcont1}a and~\ref{tsup}c. Notice that for $h=0$ this
expression is equal to $-D\left(y,t;\overline{Z}_t\right)$ because of
expression~(\ref{tDeq}) for $D$. Thus the limit of~(\ref{Dl1}) for
$h=0$ and $y\downarrow t$ is equal to $-a$. This can be compared with
the intuitive idea that quantiles $y$ approximately move to
$y-h\,D\left(y,t;\overline{Z}_t\right)$ when treatment is withheld
between $t$ and $t+h$.

Now choose $\delta=\frac{a-1}{2}$, which is greater than $0$ since we
assumed that $a>1$. By continuity of~(\ref{Dl1}) in $\left(h,y\right)$
there exists an open neighbourhood $U_{\left(0,t\right)}$ of
$\left(0,t\right)$ such that on $\left\{\left(h,y\right)\in
U_{\left(0,t\right)}: y> t+h\;{\rm and}\;h\geq 0\right\}$, the
expression~(\ref{Dl1}) above is not further than $\delta$ away from
$-a$. Thus there also exist $h_0>0$ and $y_0>t$ such that for
$h\in\left[0,h_0\right]$, $y\leq y_0$ and $y>t+h$, (\ref{Dl1}) is not
further than $\delta$ away from $-a$. Choose
$h_1\in\left[0,h_0\right]$ with $t+\left(1+\delta\right)h_1\leq y_0$,
and define $y_1=t+\left(1+\delta\right)h_1$.

Notice that since $y_1>t+h_1$,
\begin{equation*}t<F_{Y^{(t)}|\overline{Z}_t}^{-1}\circ F_{Y^{(t+h_1)}|\overline{Z}_t}\left(y_1\right)
\end{equation*}
(informally this is about withholding treatment in the future not
causing death at present, which we wanted to use; formally this
follows e.g.\ from Assumption~\ref{tsup}a and b). Moreover, for
$y=y_1$, the derivative~(\ref{Dl1}) exists on
$h\in\left[0,h_1\right]$, since for $h\in\left[0,h_1\right]$,
$y_1=t+\left(1+\delta\right)h_1>t+h_1\geq t+h$. Thus by Taylor
expansion there exist an $\tilde{h}_1\in\left[0,h_1\right]$ with
\begin{equation*}F_{Y^{(t)}|\overline{Z}_t}^{-1}\circ F_{Y^{(t+h_1)}|\overline{Z}_t}\left(y_1\right)=y_1+h_1\left.\frac{\partial}{\partial h}\right|_{h=\tilde{h}_1}
F_{Y^{(t)}|\overline{Z}_t}^{-1}\circ F_{Y^{(t+h)}|\overline{Z}_t}\left(y_1\right).\end{equation*}
Combining this it follows that
\begin{equation*}t<y_1+ h_1\left.\frac{\partial}{\partial
h}\right|_{h=\tilde{h}_1}
F_{Y^{(t)}|\overline{Z}_t}^{-1}\circ F_{Y^{(t+h)}|\overline{Z}_t}\left(y_1\right)\end{equation*}
for some $\tilde{h}_1\in\left[0,h_1\right]$. Rewriting this leads to
\begin{equation} -h_1\left.\frac{\partial}{\partial h}\right|_{h=\tilde{h}_1}
F_{Y^{(t)}|\overline{Z}_t}^{-1}\circ
F_{Y^{(t+h)}|\overline{Z}_t}\left(y_1\right)<y_1-t=\left(1+\delta\right)h_1.\label{Dl12}\end{equation} For
$\bigl(\tilde{h}_1,y_1\bigr)$, (\ref{Dl1}) is not further than $\delta$
away from $-a$, since $\tilde{h}_1\in\left[0,h_0\right]$, $y_1\leq y_0$
and $y_1>t+h_1\geq t+\tilde{h}_1$, so that
\begin{equation*}\left.\frac{\partial}{\partial h}\right|_{h=\tilde{h}_1}
F_{Y^{(t)}|\overline{Z}_t}^{-1}\circ
F_{Y^{(t+h)}|\overline{Z}_t}\left(y_1\right)\in\left(-a-\delta,-a+\delta\right).
\end{equation*}
Therefore the expression on the left hand side of equation~(\ref{Dl12}) lies in
the interval\linebreak
$\left(\left(a-\delta\right)h_1,\left(a+\delta\right)h_1\right)$.
Equation~(\ref{Dl12}) thus implies that
$\left(a-\delta\right)h_1<\left(1+\delta\right)h_1$, so $a<1+2\delta$,
so $\frac{a-1}{2}<\delta$. This is in contradiction with our choice of
$\delta$, which was $\delta=\frac{a-1}{2}$.  \hfill $\Box$

%[[Dit heeft te maken met dat $\lim_{y\downarrow
%t}D\left(y,t;\overline{Z}_t\right)\leq 1$ moet zijn dan en daar heb ik een los
%vel A t/m D over volgeschreven die wel wat ideeen geven maar nog geen
%oplossing. Wellicht is het niet nodig om te bewijzen dat dit zo
%is. Als het zonder bewijs hiervan kan dan volgt $\leq 1$ vanzelf uit
%het uiteindelijke bewijs, omdat $X(t)\sim Y^{(t)}$ given $\overline{Z}_t$.
%Laat ik
%daar nog maar even op hopen; in elk geval is de opl. $X$ uniek ok
%zonder dat Jew wet dat hi niet uit het interval wegloopt als Jew deze
%uitbreidingen steeds zo kiest. Helaas: het lijkt (zie verder) dat dit
%toch opgelost moet worden.]]

%The following remark is perhaps better postponed until Open problems
%in a remark about it being open which $D$'s are possible. No, we do
%not make a separate section open problems anymore.
%\begin{rema}
%Remark that
%A similar reasoning leads to $\lim_{y\downarrow
%t}D\left(y,t;\overline{Z}_t\right)\geq -1$, with intuition behind it that
%prolonging treatment after $t$ does not cause death at or before
%$t$.
%\end{rema}
%
%WAS FOUT BEWIJS WANT ONTBREEKT H IN LEFT HAND SIDE OF LAST INEQUALITY!!
%VB GEEFT OK AAN DAT DEZE RESTRICTIE NIET HOEFT TE GELDEN!!

$\;$

Fix $\omega$ for the rest of this section. Just as in
Section~\ref{Dexu}, it suffices to prove existence and uniqueness of
$X(t)$ with final condition on any interval between jumps of $Z$,
because with probability one $Z$ jumps only finitely many times.
Hence suppose that $Z$ does not jump in $\left(t_1,t_2\right)$ and
that $t_1$ is either a jump time of $Z$ or $0$ and that $t_2$ is
either a jump time of $Z$ or $\tau$.

If $\overline{Z}_{t_1}$ indicates that the person is dead at $t_1$,
$D$ is identically $0$ on $\left[t_1,\tau\right]$ and $X(t)$ exists,
$X(t)$ is unique, and $X(t)$ is identically $Y$ on
$\left[t_1,\tau\right]$.

If $\overline{Z}_{t_1}$ indicates that the person is alive at $t_1$,
I use Corollary~\ref{difbc3} to prove existence and uniqueness of
$X(t)$. Notice that $D\left(y,t;\overline{Z}_t\right)$ is continuous
on
$\left\{\left(y,t\right)\in\left[t_1,y_2\right]\times\left[t_1,t_2\right):y\geq
t\right\}$ because of equation~(\ref{tDeq}) and
Assumptions~\ref{tcont1}b and~\ref{tsup}c. However, the differential
equation with $X(t)$ has a final condition at the upper end of the
interval $\left[t_1,t_2\right)$. Just as in Section~\ref{Dexu}, we
define the continuous extension $\tilde{D}$ of $D$ on
$\left[t_1,y_2\right]\times\left[t_1,t_2\right]:y\geq t$, which exists
because of Assumption~\ref{tsup}c and the extension-assumption in
Assumption~\ref{tcont1}b. Just as in Section~\ref{Dexu}, $\tilde{D}$
is Lipschitz continuous in $y$ on $\left\{\left(y,t\right)\in
\left[t_1,y_2\right]\times\left[t_1,t_2\right):y>t\right\}$ with
Lipschitz constant $L_2/\varepsilon+C_2 L_1/\varepsilon^2$. The same
constant works on $\left\{\left(y,t\right)\in
\left[t_1,y_2\right]\times\left[t_1,t_2\right]:y\geq t\right\}$ by
continuity. Because of Lemma~\ref{lDl1} above,
$D\left(t,t;\overline{Z}_t\right)\leq 1$ for all $t$. Thus
Corollary~\ref{difbc3} gives existence and uniqueness of a
continuous solution $X(t)$ to the differential equation with
$\tilde{D}$, with $X\left(t\right)\geq t$ if $Y>t$.

\subsection{Mimicking counterfactual survival: discrete time}
\label{tdet}

This section considers the situation where $Z$, the available
information on the treatment- and covariate process, can be fully
described by its values at finitely many fixed points
$0=\tau_0<\tau_1<\tau_2<\ldots<\tau_K<\tau_{K+1}=\tau$. At any time at
which a person's covariates are measured, one has to include in the
covariates whether or not a person was alive at that time (otherwise,
the covariates would be ill-defined). Hence we assume that
$\overline{Z}_t$ includes whether or not a person was alive at
$\tau_1,\ldots ,\tau_{p\left(t\right)}$, with $\tau_{p\left(t\right)}$
the last $\tau_k$ before or at time $t$.
%We will not include the survival time
%itself. [[Ik kan ok nog kiezen dat als bekend is dat een person dood
%is ok de survivaltijd zelf bekend wordt, dat kijk ik nog. Decision:
%choose not to, for reasons of time]]
%Assumption~\ref{cond} stays almost the same; we just have to add the
%restriction to $\overline{Z}_{\tau_k}$ indicating that the person is alive at
%time $\tau_k$:

For simplicity we pose differentiability conditions and restrictions
on the support of $Y^{(t)}$ given $\overline{Z}_{\tau_k}$ for $t\geq
\tau_k$ that are similar to the continuous-time case. Notice that if $\tau_k$
is the lower support limit of $F_{Y|\overline{Z}_{\tau_k}}$, then
$\tau_k$ is also the lower support limit of
$F_{Y^{(t)}|\overline{Z}_{\tau_k}}$ for all $t>\tau_k$: if $\tau_k$ is
the lower support limit of $F_{Y|\overline{Z}_{\tau_k}}$ then, because
of Lemma~\ref{conslem}b, for all $h>0$ and all $\delta>0$,
$P\left(Y^{(\tau_k+h)}\leq
\tau_k-\delta\right)=P\left(Y\leq
\tau_k-\delta\right)=0$, and, again because of
Lemma~\ref{conslem}b, for all $h> 0$ and all $0<\delta<h$,
$P\left(Y^{(\tau_k+h)}\leq
\tau_k+\delta\right)=P\left(Y\leq
\tau_k+\delta\right)>0$. In most cases $\tau_k$ will then also
be the lower support limit of
$F_{Y^{(\tau_k)}|\overline{Z}_{\tau_k}}$, unless by stopping
treatment at time $\tau_k$ the person stays alive with probability
one for a certain period of time, while if treatment is not stopped
at time $\tau_k$ the hazard of dying is non-zero immediately after
time $\tau_k$.
%If time allows I can still look whether less conditions would
%suffice, but I think that's more of a mathematical challenge than
%really important since the discrete-time case is in fact just a help
%to prove the continuous-time one, so we will not give it priority for
%the moment. Decision: let it be.
For the same reasons as in the continuous-time case, differentiability
is only assumed for $y\geq t$. Assumption~\ref{cond} is replaced for
survival outcomes by:

\begin{assu} \label{tcond} \emph{(smoothness)}.
Suppose that there exists a $y_2>\tau$ such that for $k=0,\ldots,K$
and $t\in\left[\tau_k,\tau_{k+1}\right]$ there exist conditional
distribution functions $F_{Y^{(t)}|\overline{Z}_{\tau_{k}}}$ which are
consistent with Lemma~\ref{conslem} and such that if $\overline{Z}_{\tau_k}$
indicates that the person is alive at time $\tau_k$:
\begin{enumerate}[a)]
\item{For every $t\in \left[\tau_k,\tau_{k+1}\right]$,
$F_{Y^{(t)}|\overline{Z}_{\tau_k}}$ has support $\left[\tau_k,y_2\right]$.}
\item {$F_{Y^{(t)}|\overline{Z}_{\tau_{k}}}\left(y\right)$ is
continuous in $\left(y,t\right)$ on
$\left(y,t\right)\in\mathbb{R}\times
\left[\tau_k,\tau_{k+1}\right]$.}
\item{$F_{Y^{(t)}|\overline{Z}_{\tau_{k}}}\left(y\right)$ is $C^1$ in
$\left(y,t\right)$ on $\left(y,t\right)\in\left[\tau_k,y_2\right]
\times\left[\tau_k,\tau_{k+1}\right]: y> t$ with a $C^1$ extension to
$\left(y,t\right)\in\left[\tau_k,y_2\right]
\times\left[\tau_k,\tau_{k+1}\right]: y\geq t$ .}
%differentiable with respect to $y$ and $t$ on
%$\left(y,t\right)\in\left(\tau_{k},y_2\right]\times
%\left(\tau_k,\tau_{k+1}\right]$ \linebreak $: y>t$, with derivatives
%which are continuous in $\left(y,t\right)$ and have a continuous
%extension to $\left(y,t\right)\in\left(\tau_{k},y_2\right]\times
%\left(\tau_k,\tau_{k+1}\right]: y\geq t$.
\item{For $t\in\left[\tau_k,\tau_{k+1}\right)$,
$\frac{\partial}{\partial y}
F_{Y^{(t)}|\overline{Z}_{\tau_{k}}}\left(y\right)$ is strictly
positive on $y\in\left[t,y_2\right]$.}
\end{enumerate}
\end{assu}
Throughout Section~\ref{tdet}, fixed versions of
$F_{Y^{(t)}|\overline{Z}_{\tau_k}}\left(y\right)$ are used satisfying
Assumption~\ref{tcond}. Since $\overline{Z}_t$ contains the same
information as $\overline{Z}_{\tau_k}$ for
$t\in\left[\tau_k,\tau_{k+1}\right)$, we can and will choose the same
versions when conditioning on $\overline{Z}_t$.

In this discrete-time case $\overline{Z}_t$ contains no indicator for
death or alive at time $t$ except for if $t$ is one of the
$\tau_k$'s. However, also in this case $X(t)$ should be above $t$ for
$t<Y$: for such $Y$, $X(t)$ should not play the role of $Y^{(t)}$'s
less than $t$. The reason for this is, intuitively, that if
$Y^{(t)}<t$, $Y^{(t)}=Y<t$, and if also $Y$'s greater than $t$ would
play this role there would be too many of them. It will be shown
explicitly that there exists a solution $X(t)$ with $X(t)>t$ for
$Y>t$. Hence the following analogue of Theorem~\ref{detth} for
survival outcomes:

\begin{prop} \label{tdetth} Suppose that the treatment- and
covariate process $Z$ can be fully described by its values at finitely
many fixed points
$0=\tau_0<\tau_1<\tau_2<\ldots<\tau_K<\tau_{K+1}=\tau$, and suppose
also that Assumption~\ref{tcond} is satisfied. Then
$D\left(y,t;\overline{Z}_t\right)$ exists for all $t$. Furthermore if
also Assumptions~\ref{inst}, \ref{ttcons} and~\ref{ttinst}
(consistency and no instantaneous treatment effect at time of death)
are satisfied, then there exists a continuous solution $X(t)$ to
$dX(t)/dt=D\left(X(t),t;\overline{Z}_t\right)$ with final condition
$X\left(\tau\right)=Y$ and with $X(t)>t$ if $Y>t$, for which
$X\left(t\right)$ has the same distribution as $Y^{(t)}$ given
$\overline{Z}_t$.
\end{prop}

\noindent{\bf Proof.}  For $t\in\left[\tau_k,\tau_{k+1}\right)$ and
$y>t$, Lemma~\ref{quotc} can be applied on
$F_t(y)=F_{Y^{(t)}|\overline{Z}_{\tau_k}}\left(y\right)$ with
$y_0=y$, because of Assumptions~\ref{tcond} c and d. Thus
$D\left(y,t;\overline{Z}_t\right)$ as defined in equation~(\ref{tdefd}) exists
for $y>t$ and
\begin{equation} \label{tDdisc} D\left(y,t;\overline{Z}_t\right)=-\frac{\frac{\partial}{\partial t}
F_{Y^{(t)}|\overline{Z}_{\tau_k}}\left(y\right)}{\frac{\partial}{\partial
y}F_{Y^{(t)}|\overline{Z}_{\tau_k}}\left(y\right)}.\end{equation} By
definition, $D\left(t,t;\overline{Z}_t\right)$ is equal to limit of
(\ref{tDdisc}) for $y\downarrow t$, which exists because of
Assumptions~\ref{tcond} c and d.

Under Assumption~\ref{tcond} one can explicitly write down a solution to
the differential equation $dX(t)/dt=D\left(X(t),t;\overline{Z}_t\right)$
with final condition $X(\tau)=Y$, as follows. For the moment
consider $\omega\in\Omega$ fixed. Define $\tau_{p\left(\omega\right)}$
as the last $\tau_k$ before the survival time $Y$. In the following,
$F_{Y^{(\tau_{p\left(\omega\right)})}|
\overline{Z}_{\tau_{p\left(\omega\right)}}}$ will denote the
distribution function $F_{Y^{(\tau_k)}|\overline{Z}_{\tau_k}}$ for
which $k=p\left(\omega\right)$. It will \emph{not} denote a
distribution function conditional on $p\left(\omega\right)$. Define
$\tilde{X}(t)$ as follows. For $t\geq \tau_{p\left(\omega\right)+1}$, define
$\tilde{X}(t)=Y$. For $t<\tau_{p\left(\omega\right)+1}$,
$t\in\left[\tau_k,\tau_{k+1}\right)$, define
\begin{equation*} \tilde{X}(t)=F^{-1}_{Y^{(t)}|\overline{Z}_{\tau_k}}\circ F_{Y^{(\tau_{k+1})}|\overline{Z}_{\tau_{k}}}
\circ \ldots \circ
F^{-1}_{Y^{(\tau_{p\left(\omega\right)})}|\overline{Z}_{\tau_{p\left(\omega\right)}}}\circ F_{Y^{(\tau_{p\left(\omega\right)+1})}|\overline{Z}_{\tau_{p\left(\omega\right)}}}
\left(Y\right).\label{tx}
\end{equation*}
This $\tilde{X}(t)$ is well-defined under Assumption~\ref{tcond} a and b.

It is first shown that if $Y>t$ then also $\tilde{X}(t)$ is greater than
$t$. First, consider $t\in\left[\tau_{p\left(\omega\right)},
  Y\right)\cap\left[0,\tau\right]$. For such $t$, $Y>t$, so
\begin{equation*}
F_{Y^{(\tau_{k\left(\omega\right)+1})}|\overline{Z}_{p\left(\omega\right)}}\left(Y\right)
>F_{Y^{(\tau_{p\left(\omega\right)+1})}|
\overline{Z}_{p\left(\omega\right)}}\left(t\right),
\end{equation*}
since because of Assumption~\ref{tcond} a,
$F_{Y^{(\tau_{p\left(\omega\right)+1})}|\overline{Z}_{p\left(\omega\right)}}$
is strictly increasing on its support \linebreak
$\bigl[\tau_{p\left(\omega\right)},y_2\bigr]$, which includes both $t$
and $Y$. Because of Lemma~\ref{conslem}b, the right hand side of this
expression is equal to
$F_{Y^{(t)}|\overline{Z}_{p\left(\omega\right)}}\left(t\right)$. Hence
\begin{equation*}
F_{Y^{(\tau_{k\left(\omega\right)+1})}|\overline{Z}_{p\left(\omega\right)}}\left(Y\right)
>F_{Y^{(t)}|\overline{Z}_{p\left(\omega\right)}}\left(t\right),
\end{equation*}
so that since $F^{-1}_{Y^{(t)}|\overline{Z}_{\tau_{p\left(\omega\right)}}}$
is strictly increasing on $\left[0,1\right]$ (Assumption~\ref{tcond} a),
\begin{equation*}
\tilde{X}(t)=F^{-1}_{Y^{(t)}|\overline{Z}_{\tau_{p\left(\omega\right)}}}\circ
F_{Y^{(\tau_{p\left(\omega\right)+1})}|
\overline{Z}_{\tau_{p\left(\omega\right)}}} \left(Y\right)
>F^{-1}_{Y^{(t)}|\overline{Z}_{\tau_{p\left(\omega\right)}}}\circ
F_{Y^{(t)}|\overline{Z}_{p\left(\omega\right)}}\left(t\right).
\end{equation*}
The right hand side is equal to $t$, since
$t\in\left[\tau_{p(\omega)},y_2\right]$ and
$F_{Y^{(t)}|\overline{Z}_{\tau_{p\left(\omega\right)}}}$ is strictly
increasing on $\left[\tau_{p(\omega)},y_2\right]$
(Assumption~\ref{tcond} a). Thus indeed $\tilde{X}(t)>t$ for
$t\in\left[\tau_{p\left(\omega\right)}, Y\right)$.

To show that $\tilde{X}(t)$ is also greater than $t$ for other $t<Y$,
I use induction, starting with $k=p\left(\omega\right)-1$ and ending
with $k=0$. It thus needs to be proven that if $t\in\left[\tau_k,
  \tau_{k+1}\right)$ and $\tilde{X}\left(\tau_{k+1}\right)>\tau_{k+1}$
  then $\tilde{X}\left(t\right)>t$. So suppose that
  $t\in\left[\tau_k, \tau_{k+1}\right)$ and that
    $\tilde{X}\left(\tau_{k+1}\right)>\tau_{k+1}$. Notice that
    $\tilde{X}\left(\tau_{k+1}\right)>\tau_k\geq t$, so that
\begin{equation*}
F_{Y^{(\tau_{k+1})}|\overline{Z}_{\tau_k}}\left(\tilde{X}\left(\tau_{k+1}\right)\right)
>F_{Y^{(\tau_{k+1})}|\overline{Z}_{\tau_k}}\left(t\right),
\end{equation*}
since because of Assumption~\ref{tcond} a,
$F_{Y^{(\tau_{k+1})}|\overline{Z}_{\tau_k}}$ is strictly increasing on
its support $\left[\tau_{k},y_2\right]$, which includes both $t$ and
$\tilde{X}(\tau_{k+1})$. Because of Lemma~\ref{conslem}b, the right hand side of this expression is equal to $F_{Y^{(t)}|\overline{Z}_{\tau_k}}\left(t\right)$. Therefore
\begin{equation*}
F_{Y^{(\tau_{k+1})}|\overline{Z}_{\tau_k}}\left(\tilde{X}\left(\tau_{k+1}\right)\right)
>F_{Y^{(t)}|\overline{Z}_{\tau_k}}\left(t\right),
\end{equation*}
so that, because $F^{-1}_{Y^{(t)}|\overline{Z}_{\tau_{k}}}$ is strictly increasing on $\left[0,1\right]$ (Assumption~\ref{tcond} a),
\begin{equation*}
\tilde{X}(t)=F^{-1}_{Y^{(t)}|\overline{Z}_{\tau_{k}}}\circ F_{Y^{(\tau_{k+1})}|\overline{Z}_{\tau_{k}}}
\left(\tilde{X}\left(\tau_{k+1}\right)\right)
>F^{-1}_{Y^{(t)}|\overline{Z}_{\tau_k}}\circ
F_{Y^{(t)}|\overline{Z}_{\tau_k}}\left(t\right).
\end{equation*}
The right hand side of this expression is equal to $t$, since
$t\in\left[\tau_k,y_2\right]$ and
$F_{Y^{(t)}|\overline{Z}_{\tau_{k}}}$ is strictly increasing on
$\left[\tau_k,y_2\right]$ (Assumption~\ref{tcond} a). It follows that
indeed $\tilde{X}(t)>t$ if $Y>t$.

Next, it is shown that $\tilde{X}(t)$ is a continuous solution to
$\tilde{X}'(t)=D\left(\tilde{X}(t),t;\overline{Z}_t\right)$ with final
condition $\tilde{X}\left(\tau\right)=Y$. First consider $t\geq
\tau_{p(\omega)+1}$. For these $t$,
$D\left(y,t;\overline{Z}_t\right)=0$, so $\tilde{X}(t)$ should be
equal to $Y$; and indeed $\tilde{X}(t)$ is equal to $Y$. For
$t\in\left[Y,\tau_{p\left(\omega\right)+1}\right)$ it is also true
  that $D\left(y,t;\overline{Z}_t\right)=0$, so $\tilde{X}(t)$ should
  be equal to $Y$. Because of Lemma~\ref{conslem}b,
\begin{equation*}
F_{Y^{\left(t\right)}|\overline{Z}_{\tau_{p\left(\omega\right)}}}\left(Y\right)=F_{Y^{(\tau_{p\left(\omega\right)+1})}|\overline{Z}_{\tau_{p\left(\omega\right)}}}\left(Y\right),
\end{equation*}
so $\tilde{X}(t)$ is indeed equal to $Y$.

To show that $\tilde{X}(t)$ satisfies
$\tilde{X}'(t)=D\left(\tilde{X}(t),t;\overline{Z}_t\right)$ on
$\left[\tau_{p\left(\omega\right)},Y\right)$, notice that for $h\geq 0$
  small, since
  $F_{Y^{(t)}|\overline{Z}_{\tau_{p\left(\omega\right)}}}$ is
  continuous (Assumption~\ref{tcond} b),
\begin{eqnarray*} \tilde{X}\left(t+h\right)
&=&F^{-1}_{Y^{(t+h)}|\overline{Z}_{\tau_{p\left(\omega\right)}}}\circ
F_{Y^{(\tau_{p\left(\omega\right)+1})}|\overline{Z}_{\tau_{p\left(\omega\right)}}}
\left(Y\right)\\
&=&F^{-1}_{Y^{(t+h)}|\overline{Z}_{\tau_{p\left(\omega\right)}}}\circ
F_{Y^{(t)}|\overline{Z}_{\tau_{p\left(\omega\right)}}}\circ
F^{-1}_{Y^{(t)}|\overline{Z}_{\tau_{p\left(\omega\right)}}}\circ F_{Y^{(\tau_{p\left(\omega\right)+1})}|\overline{Z}_{\tau_{p\left(\omega\right)}}}
\left(Y\right)\\
&=&F^{-1}_{Y^{(t+h)}|\overline{Z}_{\tau_{p\left(\omega\right)}}}\circ
F_{Y^{(t)}|\overline{Z}_{\tau_{p\left(\omega\right)}}} \left(\tilde{X}(t)\right).
\end{eqnarray*}
This expression is differentiable at $h=0$ with derivative
$D\left(\tilde{X}(t),t;\overline{Z}_t\right)$, since, as shown
before, $\tilde{X}(t)>t$. Thus indeed $\tilde{X}(t)$ satisfies
$\tilde{X}'(t)=D\left(\tilde{X}(t),t;\overline{Z}_t\right)$ on
$\left[\tau_{p\left(\omega\right)},Y\right)$.

I still need to prove continuity of $\tilde{X}(t)$ at $t=Y$, but it is
easier to show continuity on
$\left[\tau_{p\left(\omega\right)},\tau_{p\left(\omega\right)+1}\right]$,
so we show continuity of $\tilde{X}(t)$ on
$\left[\tau_{p\left(\omega\right)},\tau_{p\left(\omega\right)+1}\right]$.
From \citet{Vaart} Lemma~21.2, $F_n$ converges weakly to $F$
if and only if $F_n^{-1}(t)\rightarrow F^{-1}(t)$ at every $t$ where
$F^{-1}$ is continuous. Notice that because of Assumption~\ref{tcond}
b, $F_{Y^{(t)}|\overline{Z}_{\tau_{p(\omega)}}}$ converges weakly to
$F_{Y^{(t_0)}|\overline{Z}_{\tau_{p(\omega)}}}$ as
$t\in\left[\tau_{p\left(\omega\right)},\tau_{p(\omega)+1}\right]\rightarrow
t_0$ for any
$t_0\in\left[\tau_{p\left(\omega\right)},\tau_{p(\omega)+1}\right]$. Moreover,
because of Assumption~\ref{tcond} a,
$F^{-1}_{Y^{(t_0)}|\overline{Z}_{\tau_k}}$ is continuous on
$\left(0,1\right]$. Therefore,
$F^{-1}_{Y^{(t)}|\overline{Z}_{\tau_{p(\omega)}}}(x)\rightarrow
F^{-1}_{Y^{(t_0)}|\overline{Z}_{\tau_{p(\omega)}}}(x)$ as
$t\in\left[\tau_{p\left(\omega\right)},\tau_{p(\omega)+1}\right]\rightarrow
t_0$, for every $x\in\left(0,1\right]$. Thus also
\begin{equation*}
\tilde{X}(t)= F^{-1}_{Y^{(t)}|\overline{Z}_{\tau_{p(\omega)}}}\circ
F_{Y^{(\tau_{p(\omega)+1})}|\overline{Z}_{\tau_{p(\omega)}}}(Y)
\rightarrow F^{-1}_{Y^{(t_0)}|\overline{Z}_{\tau_{p(\omega)}}}\circ
F_{Y^{(\tau_{p(\omega)+1})}|\overline{Z}_{\tau_{p(\omega)}}}(Y)
\end{equation*}
as
$t\in\left[\tau_{p\left(\omega\right)},\tau_{p(\omega)+1}\right]\rightarrow
t_0$ (notice that
$F_{Y^{(\tau_{p(\omega)+1})}|\overline{Z}_{\tau_{p(\omega)}}}(Y)>0$
since $Y^{(\tau_{p(\omega)+1})}$ given
$\overline{Z}_{\tau_{p(\omega)}}$ has support
  $\left[\tau_{p(\omega)},y_2\right]$ because of
  Assumption~\ref{tcond}a). For
  $t_0\in\left[\tau_{p\left(\omega\right)},\tau_{p(\omega)+1}\right)$,
    the right hand side of this expression is equal to
    $\tilde{X}(t_0)$, which implies continuity of $\tilde{X}(t)$ on
    $\left[\tau_{p\left(\omega\right)},\tau_{p(\omega)+1}\right)$. For
      $t_0=\tau_{p(\omega)+1}$, the right hand side of this expression
      is equal to
\begin{equation*}
F^{-1}_{Y^{(\tau_{p(\omega)+1})}|\overline{Z}_{\tau_{p(\omega)}}}\circ
F_{Y^{(\tau_{p(\omega)+1})}|\overline{Z}_{\tau_{p(\omega)}}}(Y),
\end{equation*}
which is equal to $Y$ since $Y$ is in the support of
$F_{Y^{(\tau_{p(\omega)+1})}|\overline{Z}_{\tau_{p(\omega)}}}$
(Assumption~\ref{tcond} a) and
$F_{Y^{(\tau_{p(\omega)+1})}|\overline{Z}_{\tau_{p(\omega)}}}$ is
strictly increasing on its support (Assumption~\ref{tcond} a). That
implies continuity of $\tilde{X}(t)$ at $\tau_{p(\omega)+1}$.

That also for $k<p(\omega)$, $\tilde{X}(t)$ satisfies
$\tilde{X}'(t)=D\left(\tilde{X}(t),t;\overline{Z}_t\right)$ on
$\left[\tau_{k},\tau_{k+1}\right)$ and that $\tilde{X}(t)$ is continuous on
  $\left[\tau_{k},\tau_{k+1}\right]$ follows the same way as in the
  previous paragraph, starting from the fact that for such $k$ and
  $t\in\left[\tau_{k},\tau_{k+1}\right)$,
    $\tilde{X}(t)=F_{Y^{(t)}|\overline{Z}_{\tau_k}}^{-1}\circ
    F_{Y^{(\tau_k)}|\overline{Z}_{\tau_k}}\left(\tilde{X}\left(\tau_{k+1}\right)\right)$.

Next it is proven that $\tilde{X}(t)$ has the same distribution as
$Y^{(t)}$ given $\overline{Z}_t$. For $t=\tau$,
$\tilde{X}\left(\tau\right)=Y$, so that $\tilde{X}(t)$ has the same
distribution as $Y^{\left(\tau\right)}$ given $\overline{Z}_\tau$
because of Assumption~\ref{inst} (consistency).  For the induction
step, suppose that for $t\in\left[\tau_k,\tau\right]$ (for $k=K+1$
read $t=\tau$), $\tilde{X}(t)$ has the same distribution as $Y^{(t)}$
given $\overline{Z}_t$. To show: for
$t\in\left[\tau_{k-1},\tau_k\right)$, $\tilde{X}(t)$ has distribution
  function
  $F_{Y^{(t)}|\overline{Z}_t}=F_{Y^{(t)}|\overline{Z}_{\tau_{k-1}}}$. If
  $\overline{Z}_{\tau_k}$ indicates that the person is dead at
  $\tau_k$ then $\tilde{X}(t)=Y=Y^{(t)}$ because of Lemma~\ref{conslem}a, so
  certainly $\tilde{X}(t)\sim Y^{(t)}$ given
  $\overline{Z}_{\tau_k}$. If $\overline{Z}_{\tau_k}$ indicates that
  the person is alive at $\tau_k$, then
\begin{equation*} \tilde{X}(t)=F^{-1}_{Y^{(t)}|\overline{Z}_{\tau_{k-1}}}\circ
F_{Y^{\left(\tau_k\right)}|\overline{Z}_{\tau_{k-1}}}
\left(\tilde{X}\left(\tau_k\right)\right)
\end{equation*}
and the rest of the proof can be copied from the proof of
Theorem~\ref{detth}.
%ja
\hfill $\Box$

%$ $
%
%\noindent Notice that $X(t)=Y=Y^{(t)}$ if $t\geq Y$
%(Lemma~\ref{conslem}a), so for $t\geq Y$ it follows that
%$X(t)=Y^{(t)}$.

\subsection{Discretization and choices of conditional distributions}
\label{tdiscsch}

The construction of the $\overline{Z}^{(n)}_{\tau_k}$ can be copied
from Section~\ref{discsch}. Notice that, by construction,
$\overline{Z}^{(n)}_{\tau_k}$ includes whether or not a person is
alive at $\tau_k$.
%ja.

\begin{nota} \label{tchoice} At this point we choose conditional distributions
$P_{\overline{Z}_{\tau_{k}}|\overline{Z}_{\tau_k}^{(n)}}$.
In addition, we choose
\begin{equation} F_{Y^{(t)}|\overline{Z}^{(n)}_{\tau_k}}\left(y\right)=\int
F_{Y^{(t)}|\overline{Z}_{\tau_k}=z}\left(y\right)
dP_{\overline{Z}_{\tau_{k}}|\overline{Z}_{\tau_k}^{(n)}}\left(z\right)
\end{equation}
to be the version of the conditional distribution function of
$Y^{(t)}$ given $\overline{Z}_{\tau_k}^{(n)}$ which is used in the
rest of the proof, and similarly for $Y$ instead of
$Y^{(t)}$. If $s\in\left(\tau_k,\tau_{k+1}\right)$ we
take the same version for $F_{Y^{(t)}|\overline{Z}_s^{(n)}}$; this is
possible because in that case
$\overline{Z}_s^{(n)}=\overline{Z}_{\tau_k}^{(n)}$.
\end{nota}

\noindent These distributions are consistent with Lemma~\ref{conslem}
in the sense that for $y\leq t$ and for all $\omega\in\Omega$,
$F_{Y^{(t)}|\overline{Z}_{\tau_k}}\left(y\right)=F_{Y|\overline{Z}_{\tau_k}}\left(y\right)$.
This follows immediately from the fact that all
$F_{Y^{(t)}|\overline{Z}_{\tau_k}}$ are consistent with
Lemma~\ref{conslem} in this sense.

\subsection{Existence of and two expressions for $D^{(n)}$}
\label{tDnsec}

We prove existence of $D^{(n)}$ as defined in equation~(\ref{tdefd})
for the discretized situation of Section~\ref{tdiscsch}. Moreover, as
in Section~\ref{Dnsec}, two useful formula's for $D^{(n)}$ are derived.

The same way as in Section~\ref{Dnsec} it follows that for $t\geq\tau_k$
and $y\in\left(t,y_2\right]$, if $\overline{Z}_{\tau_k}^{(n)}$
indicates the person is alive at $\tau_k$,
\begin{equation*}
%\lefteqn{
\left.\frac{\partial}{\partial h}\right|_{h=0}
\Bigl(F^{-1}_{Y^{(t+h)}|\overline{Z}_{\tau_k}^{(n)}}\circ
F_{Y^{(t)}|\overline{Z}_{\tau_k}^{(n)}}\Bigr)\left(y\right)
%}\nonumber\\ &&\;\hspace{1.3cm}
=-\frac{\int\left.\frac{\partial}{\partial h}\right|_{h=0}
F_{Y^{(t+h)}|\overline{Z}_{\tau_k}=z}\left(y\right)
dP_{\overline{Z}_{\tau_k}|\overline{Z}_{\tau_k}^{(n)}}
\left(z\right)}{\int\frac{\partial}{\partial y}
F_{Y^{(t)}|\overline{Z}_{\tau_k}=z}\left(y\right)
dP_{\overline{Z}_{\tau_k}|\overline{Z}_{\tau_k}^{(n)}}\left(z\right)}.
\end{equation*}
The limit for $y\downarrow t$ exists because of
Assumption~\ref{tcont1}a and Assumptions~\ref{tbdd+} and~\ref{tsup}b
(the proof is the same as the proof for continuity of this expression
in $\left(y,t\right)$ in Section~\ref{Sdisc}).  Hence with the
versions of $F_{Y^{(t)}|\overline{Z}^{(n)}_s}$ chosen in
Notation~\ref{tchoice},
\begin{equation} \label{tDn}
D^{(n)}\left(y,t;\overline{Z}_t^{(n)}\right)=\left\{\begin{array}{l} 0
\;\;\;\;\;\;\;\;\;\;\;\;\; {\rm if}\; \overline{Z}_t^{(n)}\;
{\rm indicates}\; {\rm the}\; {\rm person}\; {\rm is}\;{\rm dead}\;{\rm
at}\; t\;\;{\rm or}\; y<t\\
\left.\frac{\partial}{\partial h}\right|_{h=0}
\Bigl(F_{Y^{(t+h)}| \overline{Z}_t^{(n)}}^{-1}\circ
F_{Y^{(t)}|\overline{Z}_t^{(n)}}\Bigr)\left(y\right)\;\;\;\; {\rm
otherwise} \;{\rm for}\;y>t\\
\lim_{y\downarrow t}D\left(y,t;\overline{Z}_t\right)
\;\;\;\;\;\;\;\;\;\;\;\;\;\;\;\;\;\;\;\;\;\;\;\;\;\;\;\;\;\;\;\;\;\,
{\rm otherwise}\;{\rm for}\;y=t
\end{array}\right.
\end{equation}
exists for every $t$. Moreover, for
$t\in\left[\tau_k,\tau_{k+1}\right)$ and $y>t$,
\begin{equation} D^{(n)}\left(y,t;\overline{Z}_t^{(n)}\right)=-\frac{\int\left.\frac{\partial}{\partial h}\right|_{h=0}
F_{Y^{(t+h)}|\overline{Z}_{\tau_k}=z}\left(y\right)
dP_{\overline{Z}_{\tau_k}|\overline{Z}_{\tau_k}^{(n)}} \left(z\right)}{\int\frac{\partial}{\partial y}
F_{Y^{(t)}|\overline{Z}_{\tau_k}=z}\left(y\right)
dP_{\overline{Z}_{\tau_k}|\overline{Z}_{\tau_k}^{(n)}}\left(z\right)}.
\label{tDneq2}\end{equation}
This expression is similar to expression~(\ref{Dneq2}) for $D^{(n)}$
for non-survival outcomes.
%Remark that (\ref{tFform})
%is continuous in $\left(y,t\right)$ for $y>\tau_k$.
%ja

Expression~(\ref{Dneq}) for non-survival outcomes takes a different
form for survival outcomes. Just as for expression~(\ref{Dneq}), I
restrict to the $\Omega'$ defined in equation~(\ref{omeac}) in
Section~\ref{Dnsec}, a set of probability one on which conditional
probabilities given $\overline{Z}^{(n)}_t$ are uniquely defined.
%Contrarily to (\ref{tDneq2}) we naturally get an expression for all
%$t$, also for $t$ equal to one of the $\tau_k$'s.
For $n$ and $\omega\in\Omega'$ such
that $\overline{Z}_{t}^{(n)}$ indicates the person is alive at
the last $\tau_k$ at or before time $t$, it will be shown that for $y>t$,
\begin{equation}
D^{(n)}\left(y,t;\overline{Z}_t^{(n)}\right)=-\frac{E\left[1_{\left\{{\rm
alive}\;{\rm at}\;t\right\}}\,\left.\frac{\partial}{\partial
h}\right|_{h=0}
F_{Y^{(t+h)}|\overline{Z}_t}\left(y\right)\big|\overline{Z}_t^{(n)}\right]}
{E\left[1_{\left\{{\rm alive}\;{\rm at}\;t\right\}}\,
\frac{\partial}{\partial y}
F_{Y^{(t)}|\overline{Z}_t}\left(y\right)\big|\overline{Z}_t^{(n)}\right]}
.\label{tDneq}
\end{equation}
The indicator of being alive at time $t$ is new as compared to the
non-survival case of Section~\ref{XY}.

To prove equation~(\ref{tDneq}), first restrict to
$\omega\in\Omega'$. Suppose that $\overline{Z}_{\tau_k}^{(n)}$
indicates the person is alive at $\tau_k$ and suppose that $t\geq
\tau_k$. Then for $h\geq 0$ and $y\geq t+h$,
\begin{eqnarray*}  P\left(Y^{(t+h)}\leq
y\big|\overline{Z}_{\tau_k}^{(n)}\right) &=&P\left(Y^{(t+h)}\leq
y\big|\overline{Z}_{\tau_k}^{(n)}, Y\in\left(\tau_k,t\right]\right)
P\left(Y\in\left(\tau_k,t\right]\big|\overline{Z}_{\tau_k}^{(n)}\right)
\nonumber\\ &&+\int_{z_t:\;{\rm alive}\;{\rm at}\; t}
F_{Y^{(t+h)}|\overline{Z}_t=z_t}\left(y\right)dP_{\overline{Z}_t|\overline{Z}_{\tau_k}^{(n)}}\left(z_t\right).
\end{eqnarray*}
Given that $Y\leq t$, Lemma~\ref{conslem}a gives that $Y^{(t+h)}=Y\leq
t$, and, since $y\geq t$, also $Y^{(t+h)}\leq t\leq y$. It follows that
\begin{equation}
\label{ttFform}
P\left(Y^{(t+h)}\leq y\big|\overline{Z}_{\tau_k}^{(n)}\right)
=P\left(Y\in\left(\tau_k,t\right]\big|\overline{Z}_{\tau_k}^{(n)}\right)
+ \int_{z_t} 1_{\left\{{\rm
alive}\;{\rm at}\;t\right\}}\, F_{Y^{(t+h)}|\overline{Z}_t=z_t}\left(y\right)
dP_{\overline{Z}_t|\overline{Z}_{\tau_k}^{(n)}}\left(z_t\right).
\end{equation}

To derive equation~(\ref{tDneq}) from equation~(\ref{ttFform}),
Corollary~\ref{quotc} is applied on
\begin{equation*}F_h\left(y\right)=P\left(Y\in
\left(\tau_k,t\right]\big|\overline{Z}_{\tau_k}^{(n)}\right)
+\int_{z_t} 1_{\left\{{\rm alive}\;{\rm at}\;t\right\}}\,
F_{Y^{(t+h)}|\overline{Z}_{t}=z_t}\left(y\right)
dP_{\overline{Z}_{t}|\overline{Z}_{\tau_k}^{(n)}}\left(z_t\right),
\end{equation*}
for $y>t$, with $y_0=y$. We check the conditions of
Corollary~\ref{quotc}. Just as in Section~\ref{Dnsec}, on $y>t+h$,
$h\geq 0$, $F_h\left(y\right)$ is differentiable with respect to $y$
and $h$ with derivatives $\int\! 1_{\left\{{\rm alive}\;{\rm
    at}\;t\right\}}\,
F_{Y^{(t+h)}|\overline{Z}_{t}=z_t}'\left(y\right)
dP_{\overline{Z}_{t}|\overline{Z}_{\tau_k}^{(n)}}\left(z_t\right)$ and
$\int\!1_{\left\{{\rm alive}\;{\rm at}\;t\right\}}\,
\frac{\partial}{\partial
  h}\!F_{Y^{(t+h)}|\overline{Z}_{t}=z_t}\!\left(y\right)
dP_{\overline{Z}_{t}|\overline{Z}_{\tau_k}^{(n)}}\left(z_t\right)$,
respectively. Also the same way as in Section~\ref{Dnsec}, it follows
that these derivatives are continuous in $\left(y,h\right)$. That
$F'_0\left(y\right)$ is non-zero follows from Assumption~\ref{tsup}c,
if the probability that the person is alive at time $t$ given
$\overline{Z}_{\tau_k}^{(n)}$ is non-zero. Indeed the probability that
the person is alive at time $t$ given $\overline{Z}_{\tau_k}^{(n)}$ is
non-zero, since given any $\overline{Z}_{\tau_k}$ indicating that the
person is not dead at $\tau_k$, $Y$, which has the same distribution
as $Y^{\left(\tau\right)}$ given $\overline{Z}_{\tau_k}$ because of
Assumption~\ref{inst}, has support $\left[\tau_k,y_2\right]$
(Assumption~\ref{tsup}a).  Thus the conditions of
Corollary~\ref{quotc} are satisfied, and equation~(\ref{tDneq})
follows from equation~(\ref{ttFform}).

%[[This formula is only used to show convergence to $D$,
%so it seems that we have to bound $X^{(n)}-X$ in terms of $X$ instead of
%$X^{(n)}$ and prove that $X(t)>t$ if $Y>t$ anyway, which we did not
%manage to do in Section~\ref{tDexu}, otherwise this formula cannot be
%used. I'm afraid the limit there is easier than proving things
%here. THAT WAS INDEED THE CASE.]]

\subsection{Applying the discrete-time result}
\label{tSdisc}

\begin{lem} Suppose that Regularity Conditions~\ref{tsup}--\ref{tLipc}
and Assumptions~\ref{inst}, \ref{ttcons} and~\ref{ttinst} (consistency
and no instantaneous treatment effect at time of death) are
satisfied. Then for every $n$ there exists a continuous solution
$X^{(n)}(t)$ to the differential equation in the discretised setting,
\begin{equation*}
\frac{d}{dt}X^{(n)}(t)=D^{(n)}\Bigl(X^{(n)}(t),t;\overline{Z}^{(n)}_t\Bigr),
\end{equation*}
with final condition
$X^{(n)}\left(\tau\right)=Y$. $X^{(n)}\left(t\right)$ is almost surely
unique. Moreover, $X^{(n)}(t)>t$ if $Y>t$. Furthermore,
$X^{(n)}\left(t\right)$ has the same conditional distribution as
$Y^{(t)}$ given $\overline{Z}_t^{(n)}$.
\end{lem}

\noindent The proof of this lemma is different from the proof in
Section~\ref{Sdisc}, because of the different assumptions for the
discrete-time case if the outcome is survival.

$\;$

\noindent{\bf Proof.} Fix $n$. First, it is shown that there exists a
continuous solution $X^{(n)}$ with $X^{(n)}(t)>t$ if $Y>t$ for which
$X^{(n)}(t)$ has the same conditional distribution as $Y^{(t)}$ given
$\overline{Z}^{(n)}_t$, using Proposition~\ref{tdetth}. We thus need to
check that the versions of the conditional distributions
$F_{Y^{(t)}|\overline{Z}^{(n)}_{\tau_k}}$ of $Y^{(t)}$ given
$\overline{Z}^{(n)}_{\tau_k}$ chosen in Notation~\ref{tchoice} satisfy
Assumption~\ref{tcond}.

Section~\ref{tdiscsch} already showed that the conditional
distributions $F_{Y^{(t)}|\overline{Z}^{(n)}_{\tau_k}}\left(y\right)$
chosen in Notation~\ref{tchoice} are consistent with
Lemma~\ref{conslem}. We check Assumption~\ref{tcond} a--d. If
$\overline{Z}_{\tau_k}^{(n)}$ indicates that the person is alive at
time $\tau_k$ then:
\begin{enumerate}[a)]
\item{$F_{Y^{(t)}|\overline{Z}^{(n)}_{\tau_k}}$ has support
$\left[\tau_k,y_2\right]$ since all
$F_{Y^{(t)}|\overline{Z}_{\tau_k}=z}$ have support
$\left[\tau_k,y_2\right]$ (Assumption~\ref{tsup} a).}
\item{$F_{Y^{(t)}|\overline{Z}^{(n)}_{\tau_k}}(y)=F_{Y^{\left(\tau_k+\left(t-\tau_k\right)\right)}|
    \overline{Z}_{\tau_k}^{(n)}}(y)$ is continuous in $(y,t)$ on
  $(y,t)\in\mathbb{R}\times\left[\tau_k,\tau_{k+1}\right]:y\geq t$
  because of Assumption~\ref{tcont1} a and Lebesgue's dominated
  convergence theorem, since all $F_{Y^{(t)}|\overline{Z}_{\tau_k}}$
  are bounded by $1$. For $y\leq t$,
  $F_{Y^{(t)}|\overline{Z}_{\tau_k}}(y)=F_{Y|\overline{Z}_{\tau_k}}(y)$
  because of Lemma~\ref{conslem}, which is continuous in
  $\left(y,t\right)$ because of Assumption~\ref{Ytsup}. Therefore, the
  same is true for the version of
  $F_{Y^{(t)}|\overline{Z}^{(n)}_{\tau_k}}\left(y\right)$ chosen in
  Notation~\ref{tchoice}.}
\item{$F_{Y^{(t)}|\overline{Z}^{(n)}_{\tau_k}}(y)
=F_{Y^{\left(\tau_k+\left(t-\tau_k\right)\right)
|\overline{Z}^{(n)}_{\tau_k}}}\left(y\right)$ is $C^1$ in $(y,t)$
on $(y,t)\in\left[\tau_k,y_2\right]\times
\left[\tau_k,\tau_{k+1}\right]:y\geq t$ because all
$F_{Y^{(t)}|\overline{Z}_{\tau_k}}(y)$ are $C^1$ there
(Assumption~\ref{tcont1} a), and all
derivatives are bounded there (Assumption~\ref{tbdd+}), which follows with the same reasoning as for
$F_{Y^{(t+h)}|\overline{Z}^{(n)}_{\tau_k}}$ as in Section~\ref{Dnsec}:
integration and differentiation can be interchanged here.}
\item{Under c) it was shown that for $y\geq t$,
\begin{equation*}
\frac{\partial}{\partial y}F_{Y^{(t)}|\overline{Z}^{(n)}_{\tau_k}}\left(y\right)=\int
\frac{\partial}{\partial y}F_{Y^{(t)}|\overline{Z}_{\tau_k}=z}\left(y\right)
dP_{\overline{Z}_{\tau_{k}}|\overline{Z}_{\tau_k}^{(n)}}\left(z\right).
\end{equation*}
Because of Assumption~\ref{tsup} b this is greater than $0$. \hfill $\Box$}
\end{enumerate}

%Lemma~\ref{detlem} and its proof can be copied from Section~\ref{Sdisc}, with
%equation~(\ref{tFform}) replacing equation~(\ref{Fform}).

%$X^{(n)}$ is unique
%under the conditions in Section~\ref{tcondc} since $D^{(n)}$ is Lipschitz
%continuous in $y$ with Lipschitz constant $\frac{1}{\delta}L_2+\frac{C_2}{\varepsilon^2}L_1$ for every
%$\omega\in\Omega$. This follows by the same reasoning as
%$\tilde{D}$ in Section~\ref{Dexu}, but starting from (\ref{tDneq}) instead of
%(\ref{Deq}).

Continuity of $D^{(n)}$ on
$\left\{\left(y,t\right)\in\left(\tau_k,y_2\right]\times\left(\tau_k,\tau_{k+1}\right):
y>t\right\}$ with a continuous extension $\tilde{D}^{(n)}$ to
$\left\{\left(y,t\right)\in\left[\tau_k,y_2\right]
\times\left[\tau_k,\tau_{k+1}\right]: y\geq t\right\}$ follows from
%$\mathbb{R}\times\left[\tau_k,\tau_{k+1}\right]$ using
equation~(\ref{tDneq2}), similar to Section~\ref{Sdisc}.
%\begin{equation*} \tilde{D}^{(n)}\left(y,t;\overline{Z}_t^{(n)}\right)=-\frac{\int\left.\frac{\partial}{\partial h}\right|_{h=0}
%F_{Y^{(t+h)}|\overline{Z}_{\tau_k}=z}\left(y\right)
%dP_{\overline{Z}_{\tau_k}|\overline{Z}_{\tau_k}^{(n)}} \left(z\right)}{\int\frac{\partial}{\partial y}
%F_{Y^{(t)}|\overline{Z}_{\tau_k}=z}\left(y\right)
%dP_{\overline{Z}_{\tau_k}|\overline{Z}_{\tau_k}^{(n)}}\left(z\right)},
%\end{equation*}
Also similar to Section~\ref{Sdisc}, $\tilde{D}^{(n)}(y,t)$ is
Lipschitz continuous in $y$ on\linebreak
$\left\{\left(y,t\right)\in\left[\tau_k,y_2\right]
\times\left[\tau_k,\tau_{k+1}\right]: y\geq t\right\}$ for
$\omega\in\Omega'$ of equation~(\ref{omeac}), and uniqueness of
$X^{(n)}$ on $\Omega'$ follows from Corollary~\ref{difbc3}.
%Define a continuous extension
%$\tilde{D}^{(n)}$ of $D^{(n)}$ the same way as $\tilde{D}$ in
%Section~\ref{tDexu} for $D$. It exists and is continuous.

% because of
%the extension-part of Assumption~\ref{tcont1}a and
%Assumption~\ref{tsup}b. However, $\tilde{D}^{(n)}$ need not be equal to
%$D^{(n)}$ itself for $y=\tau_k$. However, this is not a
%real problem, since with probability one $\tilde{D}^{(n)}$ is not used
%for $y=\tau_k$. Let's look at the proof of Theorem~\ref{detth}. Then
%we see that the solution $X^{(n)}(t)$ stays strictly above $\tau_k$ on
%$\left[\tau_k,\tau_{k+1}\right)$ if $Y>\tau_k$, which can be seen
%from (\ref{tx}): all $F_{Y^{(t)}|\overline{Z}_{\tau_k}}$ ($t\geq \tau_k$) are
%strictly increasing on their support $\left[\tau_k,y_2\right]$, as we saw in
%Section~\ref{tSdisc}, and thus a q-q transform
%$F^{-1}_{Y^{(t)}|\overline{Z}_{\tau_k}}\circ
%F_{Y^{\left(\tau_{k+1}\right)}|\overline{Z}_{\tau_k}}$ ($t\geq \tau_k$) can
%never map an $X\left(\tau_{k+1}\right)>\tau_k$ to an
%$X\left(t\right)\leq \tau_k$. So only if

\subsection{Bounding the difference between $X$ and $X^{(n)}$ in
terms of $D$ and $D^{(n)}$}

$X^{(n)}$ and $X$ satisfy the differential equations with the
continuous extensions of $D^{(n)}$ and $D$, $\tilde{D}^{(n)}$ and
$\tilde{D}$, respectively, on the closed intervals
$\left[t_1,t_2\right]$ as in Section~\ref{XXnaf}, because if $Y\geq
t$, both $X(t)\geq t$ (Section~\ref{tDexu}) and $X^{(n)}(t)\geq t$
(Section~\ref{tSdisc}). $\tilde{D}$ is Lipschitz continuous in $y$ on
these intervals, as shown at the end of Section~\ref{tDexu}. Therefore
it follows in a similar way as in Section~\ref{XXnaf} but with
Corollary~\ref{difbc3} instead of Theorem~\ref{difbc2}
that almost surely
\begin{equation} \bigl|X^{(n)}(t) -X(t)\bigr|
\leq \int_t^{\tau}
e^{C\cdot\left(s-t\right)}
\bigl|D\bigl(X^{(n)}\left(s\right),s\bigr)-
D^{(n)}\bigl(X^{(n)}\left(s\right),s\bigr)\bigr| ds
\label{tXbd}\end{equation}
and
\begin{equation} \sup_{t\in\left[0,\tau\right]}\bigl|X^{(n)}(t) -X(t)\bigr|
=  \int_0^{\tau} e^{C\cdot s}
\bigl|D\bigl(X^{(n)}\left(s\right),s\bigr)-
D^{(n)}\bigl(X^{(n)}\left(s\right),s\bigr)\bigr| ds,\label{tXschat}
\end{equation}
with $C=L_2/\varepsilon+C_2L_1/\varepsilon^2$.
%[[maybe v/^ ipv min and max if that often makes a difference.]]

%We apply Corollary~\ref{difbc1} with $z=X$ instead of $X^{(n)}$ since we
%hope that $X(t)\geq t$ if $Y\geq t$ and there is no hope that this is
%true for $X^{(n)}$. If we assume that $\lim_{y\downarrow
%t}D\left(y,t;\overline{Z}_t\right)\leq 1$ then we get the same way as in
%Section~\ref{Xas} that
%\begin{equation} \left|X^{(n)}(t) -X(t)\right|
%\leq \int_t^\tau e^{C\cdot\left(s-t\right)}
%\left|D\left(X\left(s\right),s\right)-
%D^{(n)}\left(X\left(s\right),s\right)\right| ds.
%\label{tXbd}\end{equation}
%If we do not assume that $\lim_{y\downarrow
%y}D\left(y,t;\overline{Z}_t\right)\leq 1$ we get the same expression but with
%$\tilde{D}$ and $\tilde{D}^{(n)}$.

\subsection{Convergence of $D^{(n)}$ to $D$}

This section proves that for all $\left(y,t\right)$ fixed,
$D^{(n)}\bigl(y,t;\overline{Z}_t^{(n)}\bigr)$ converges almost surely
to $D\left(y,t;\overline{Z}_t\right)$. First, consider
$\overline{Z}_t$ indicating the person is dead at time $t$. Then,
$D\left(y,t;\overline{Z}_t\right)=0$. Note that the probability the
person died at exactly time $t$ is $0$. Therefore, almost surely, for
$n$ large enough $\overline{Z}_t^{(n)}$ indicates the person is dead
at the last $\tau_k$ at or before time $t$. Thus, also
$D^{(n)}\bigl(y,t;\overline{Z}_t^{(n)}\bigr)=0$, so $D^{(n)}$
converges to $D$.  For $y<t$, $D\left(y,t;\overline{Z}_t\right)$ and
all $D^{(n)}\bigl(y,t;\overline{Z}_t^{(n)}\bigr)$ are $0$. Therefore
it suffices to consider $y\geq t$ and $\omega$ and $t$ for which
$\overline{Z}_t$ indicates the person is alive at time $t$. We start by
proving that for $y>t$, $D^{(n)}\bigl(y,t;\overline{Z}_t^{(n)}\bigr)$
converges almost surely to $D\left(y,t;\overline{Z}_t\right)$.
Equation~(\ref{tDeq}) implies that if $\overline{Z}_t$ indicates the
person is alive at time $t$, for $y>t$
\begin{equation*} D\left(y,t;\overline{Z}_t\right)=-\frac{\left.\frac{\partial}{\partial h}\right|_{h=0}
F_{Y^{(t+h)}|\overline{Z}_t}\left(y\right)}{\frac{\partial}{\partial
      y}F_{Y^{(t)}|\overline{Z}_t}\left(y\right)},\end{equation*} and
equation~(\ref{tDneq}) implies that, since the person is alive at the
last $\tau_k$ at or before time $t$, for $y>t$
\begin{equation*} D^{(n)}\left(y,t;\overline{Z}_t^{(n)}\right)=-\frac{E\left[1_{\left\{{\rm alive}\;{\rm at}\;t\right\}}\left.\frac{\partial}{\partial h}\right|_{h=0}
F_{Y^{(t+h)}|\overline{Z}_t}\left(y\right)\big|\overline{Z}_t^{(n)}\right]}
{E\left[1_{\left\{{\rm alive}\;{\rm at}\;t\right\}}\frac{\partial}{\partial y}
F_{Y^{(t)}|\overline{Z}_t}\left(y\right)\big|\overline{Z}_t^{(n)}\right]}.\end{equation*}
We will apply L\'evy's Upward Theorem (see e.g.~\citet{Wil} page 134),
which is allowed since
%~\ref{Lev}
$\left.\frac{\partial}{\partial h}\right|_{h=0}
F_{Y^{(t+h)}|\overline{Z}_t}\left(y\right)$ and
$\frac{\partial}{\partial y} F_{Y^{(t)}|\overline{Z}_t}\left(y\right)$
are bounded because of Assumption~\ref{bdd+}. L\'evy's Upward Theorem
leads to
\begin{eqnarray} \lefteqn{E\left[1_{\left\{{\rm alive}\;{\rm at}\;t\right\}}\left.\frac{\partial}{\partial h}\right|_{h=0}
F_{Y^{(t+h)}|\overline{Z}_t}\left(y\right)\Big|\overline{Z}_t^{(n)}\right]}\nonumber\\
&&\;\hspace{0.5cm}\rightarrow E\left[1_{\left\{{\rm alive}\;{\rm at}\;t\right\}}\left.\frac{\partial}{\partial h}\right|_{h=0}
F_{Y^{(t+h)}|\overline{Z}_t}\left(y\right)\bigg|
\sigma\left(\cup_{n=1}^\infty\overline{Z}_t^{(n)}\right)
\right] \;\mbox{\rm a.s.}\label{th1}
\end{eqnarray}
and
\begin{equation} E\left[1_{\left\{{\rm alive}\;{\rm at}\;t\right\}}\frac{\partial}{\partial y}
F_{Y^{(t)}|\overline{Z}_t}\left(y\right)\Big|\overline{Z}_t^{(n)}\right] \rightarrow
E\left[1_{\left\{{\rm alive}\;{\rm at}\;t\right\}}\frac{\partial}{\partial
y}F_{Y^{(t)}|\overline{Z}_t}\left(y\right)\bigg|\sigma\left(\cup_{n=1}^\infty
\overline{Z}_t^{(n)}\right) \right] \;\mbox{\rm a.s.}. \label{th2}
\end{equation}
Replacing the conditioning on $\sigma\bigl(\cup_{n=1}^\infty
\overline{Z}_t^{(n)}\bigr)$ by conditioning on $\overline{Z}_t$ in
(\ref{th1}) and (\ref{th2}) is allowed because of
Lemma~\ref{Zt-}. Since for these $\overline{Z}_t$,
\begin{equation*}E\left[1_{\left\{{\rm alive}\;{\rm at}\;t\right\}}\left.\frac{\partial}{\partial h}\right|_{h=0}
F_{Y^{(t+h)}|\overline{Z}_t}\left(y\right)\Big|\overline{Z}_t\right]=1_{\left\{{\rm alive}\;{\rm at}\;t\right\}}
E\left[\left.\frac{\partial}{\partial h}\right|_{h=0}
F_{Y^{(t+h)}|\overline{Z}_t}\left(y\right)\Big|\overline{Z}_t\right]
\end{equation*}
and
\begin{equation*}E\left[1_{\left\{{\rm alive}\;{\rm at}\;t\right\}}\frac{\partial}{\partial
y}F_{Y^{(t)}|\overline{Z}_t}\left(y\right)\Big|\overline{Z}_t\right] =1_{\left\{{\rm alive}\;{\rm at}\;t\right\}}
E\left[\frac{\partial}{\partial
y}F_{Y^{(t)}|\overline{Z}_t}\left(y\right)\Big|\overline{Z}_t\right],\end{equation*}
this implies that for $y>t$ fixed
\begin{equation}
D^{(n)}\left(y,t;\overline{Z}_t^{(n)}\right)
\rightarrow D\left(y,t;\overline{Z}_t\right) \label{tDconv}\;\mbox{\rm a.s.}.
\end{equation}

To prove that also
\begin{equation}
D^{(n)}\left(t,t;\overline{Z}_t^{(n)}\right)
\rightarrow D\left(t,t;\overline{Z}_t\right) \label{tDconvt}\;\mbox{\rm a.s.},
\end{equation}
We will apply Lemma~\ref{plugl}. To do that it suffices to have
Lipschitz continuity of all $D^{(n)}$ and $D$ in $y$ with the same
Lipschitz constant, and if that is the case equation~(\ref{tDconvt})
follows.  That $D$ is Lipschitz continuous in $y$ with Lipschitz
constant $\frac{1}{\varepsilon}L_2+\frac{C_2}{\varepsilon^2}L_1$ was
shown in Section~\ref{tDexu}. That $D^{(n)}$ is Lipschitz continuous
in $y$ with Lipschitz constant
$\frac{1}{\varepsilon}L_2+\frac{C_2}{\varepsilon^2}L_1$ on $\Omega'$
was shown in Section~\ref{tSdisc}.

\subsection{$X^{(n)}(t)$ converges to $X(t)$ and $X(t)$ is measurable}

Equations~(\ref{tXschat}), (\ref{tDconv}) and (\ref{tDconvt}) are the
starting point here. If $Y\geq t$, both $X(t)\geq t$
(Section~\ref{tDexu}) and $X^{(n)}(t)\geq t$ (Section~\ref{tSdisc}), so
that the rest of the proof can be copied from Section~\ref{XnX}.

%This section can be copied from Section~\ref{Xas}, since we saw in
%Section~\ref{tDexu} that $X(t)\geq t$ if $Y\geq t$. Remark that the
%fact that the derivative of $X^{(n)}$ may be discontinuous in $t=Y$
%does not cause any trouble since $\overline{Z}_t$ jumps at $t=Y$ and thus was already
%treated separately in Section~\ref{Xas}.

\subsection{Conclusion}

This section can be copied from Section~\ref{Concl}.
%Ja

\subsection{Mimicking counterfactual outcomes: discrete-continuous time}

This section can be copied from Section~\ref{Disc2}.

\section{Web-Appendix: A simulation study}\label{simapp}

\subsection{Introduction}

This appendix provides further details on the design of the simulation study of Section~9. In the simulation study, we calibrated the distributions of the variables and the parameter values to HIV/AIDS data, perhaps the most salient example of application of structural nested models in the empirical literature. For details, see Section~\ref{HIV}.

Web-Appendix~\ref{simapp} is organized as follows. Section~\ref{add} presents the counterfactual outcomes and the treatment initiation process. The outcomes and treatment initiation are designed so that treatment predicts intermediate covariates which in turn predict future treatment: the type of setting structural nested models were designed for. Continuous-time structural nested models can be adopted when the treatment is initiated in continuous time. Consequently, in the simulations we impose that treatment decisions are adopted continuously in time. Section~\ref{norp} shows that our setting does not impose (local) rank preservation. In Section~\ref{HIV}, the distributions and parameters are calibrated to real data on HIV/AIDS. Section~\ref{Dder} derives the parametric form of the infinitesimal shift-function $D$. Section~\ref{Xpsisec} derives the solution to the differential equation (8): the ``mimicking'' variable $X_\psi(t)$ of Section~4. Section~\ref{est} describes the estimators. Section~\ref{res} describes the results of the simulation study.

\subsection{An additive model for treatment effect, and the treatment initiation process}\label{add}

In this simulation study no one is treated at time zero, and once treatment is initiated, it is never stopped. $Y^{(t)}$ is the counterfactual outcome had treatment been as given in reality until time $t$, and continued or initiated after that. For example, if treatment was initiated by time $t$ for a particular patient, $Y^{(t)}$ is the observed outcome for that patient, since he or she was already treated at time $t$ and treatment is never stopped. On the other hand, if treatment was not initiated by time $t$, $Y^{(t)}$ is the outcome had treatment been initiated at time $t$. Thus, in the definition of $Y^{(t)}$ in Section~\ref{sets}, the switch at time $t$ to ``some kind of baseline treatment regime $\overline{0}$'' is, in this case, ``treat continuously'' from time $t$ onwards. In the simulations, we study a setting with $t\in[0,2]$. The subscript $_{t}$ indicates the treatment initiation time, so for example $L_{1,t}$ indicates (counterfactual) covariates at time $1$ under ``treatment started at time $t$''. Similarly, the subscript $_{\infty}$ indicates (counterfactual) variables under no treatment. For example, $L_{2,\infty}$ indicates (counterfactual) covariates at time $2$ under no treatment. In the simulation design, the counterfactual covariates $L$ are as follows:
\begin{eqnarray*}
L_0&=&\tilde{L}_0+e_0,\\
L_{1,\infty}&=&\tilde{L}_0-\beta_0+e_{1,\infty},\\
L_{2,\infty}&=&\tilde{L}_0-2\beta_0+e_{2,\infty}\\
L_{1,t}&=&\tilde{L}_0-\beta_0+\theta (1-t)+e_{1,t} \;{\rm for}\;t\in\left[0,1\right], \;{\rm and}\; L_{1,\infty}\;{\rm otherwise}\\
L_{2,t}&=&\tilde{L}_0-2\beta_0+\psi (2-t)+e_{2,t},
\end{eqnarray*}
 where $\tilde{L}_0$ and the $e_{j,t}$ are random variables with values in $\mathbb{R}$. Notice that $(1-t)$ and $(2-t)$ are simply the durations of treatment until the respective covariate measurements. We will assume that the $e_{j,t}$ ($j=0,1,2$) are independent of $\tilde{L}_0$, and that the $e_{2,t}$ have a distribution function which does not depend on $t$. We will also assume that the $e_{2,t}$ are independent of all previous variables (and of the treatment initiation time, $T$, described below). In the simulations, $\psi\geq 0$ (a similar study could have been done for $\psi<0$). We define $Y_{t}=L_{2,t}$, the counterfactual outcome with treatment initiated at time $t$, which could potentially be observed at time $2$.

$T$ will be the treatment initiation time, with $T=\infty$ if treatment was not initiated in the time interval $[0,2]$. The treatment initiation time $T$ determines which of the above variables is observed. $Y=L_2=Y_{T}$ no matter when treatment is started. $L_1=L_{1,T}$ if $T\leq 1$ and $L_1=L_{1,\infty}$ if $T>1$. $T$ also determines what are the $Y^{(t)}$, with $Y^{(t)}=Y_t$ if $T>t$ and $Y^{(t)}=Y_T$ if $T\leq t$.

Suppose that the hazard of the treatment initiation time $T$, given the covariate history at time $t$ and given that treatment was not initiated before time $t$, is piecewise constant as follows:
\begin{equation*}
\lambda_T(t)=\left\{\begin{array}{ll}
\lambda^{(0)}_0 & {\rm if}\; L_0> c_0 \;{\rm and}\;t\in[0,1]\\
\lambda^{(0)}_1 & {\rm if}\; L_0\leq c_0 \;{\rm and}\;t\in[0,1]\\
\lambda^{(1)}_0 & {\rm if}\; L_{1,\infty} > c_1 \;{\rm and}\;t\in(1,2]\\
\lambda^{(1)}_1 & {\rm if}\; L_{1,\infty} \leq c_1 \;{\rm and}\;t\in(1,2],
\end{array}\right.
\end{equation*}
for constants $c_0$ and $c_1$ in $\mathbb{R}$. Notice that $T$ depends on $\tilde{L}_0$, $e_{0,\infty}$, and $e_{1,\infty}$, if $\lambda^{(0)}_0\neq \lambda^{(0)}_1$ or $\lambda^{(1)}_0\neq \lambda^{(1)}_1$.

In the simulation study, treatment can be initiated in continuous time, but the covariates are only measured at times $0$, $1$, and $2$, so that the treatment and covariate history up to time $t$, $\overline{Z}_t$, consists of the treatment information up to time $t$ and $L_0$, $(L_0,L_1)$, or $(L_0,L_1,L_2)$, depending on whether $t\in[0,1)$, $t\in[1,2)$, or $t=2$.

In the simulations, treatment affects later outcomes, and time-dependent covariates ($L_1$) which depend on previous treatment also predict future treatment and the outcome of interest. This is the type of setting structural nested models were developed for.

\subsection{No rank preservation in the simulations}\label{norp}

The outcomes described in Section~\ref{add} are not rank preserving. %The reason is, that the $e_{2,t}$'s can be different for all $t$, although most likely they will be correlated. In addition,
For two patients with the same treatment history, the complete observed data are the same if the sum of $\tilde{L}_0$ and the $e_{j,T}$ (for $j=0,1,2$) are all three the same. However, under an alternative treatment, the outcomes for two patients with the same observed data can be different. In fact, in this simulation study, they are different with probability one. This is easily seen because with probability one, the value of $\tilde{L}_0$ is not the same for these two patients. Under rank preservation, two patients with the same observed data ($(L_0,L_1,L_2,T,Y)$ the same for both patients) also would have had the same outcomes had they both followed the same alternative treatment ($Y_t$ the same for both patients). Thus, rank preservation does not hold in this simulation study.
%In case of rank preservation, $e_{2,t}$ would equal $e_{2,\infty}$ for all $t$. This was off. The above is much better!

\subsection{Choice of parameter values in the simulation study}\label{HIV}

In the simulation study, we calibrate the distributions of the variables and the parameter values to HIV/AIDS data. We focus on the first two years since HIV diagnosis. Time zero is the time of HIV diagnosis. The outcome variable is the CD4 count, a commonly used marker of the state of the immune system of HIV-positive patients. The usual treatment for HIV-positive patients is ART, antiretroviral treatment. ART is not always initiated immediately after diagnosis. ART initiation time often depends on the last measured CD4 count. When the CD4 count is at or below $350$ copies/ml, HIV-positive patients are much more likely to initiate ART than when the CD4 count is above $350$ copies/ml. Thus, in the simulation study we choose $c_1=c_0=350$. Intermediate CD4 counts are affected by previous treatment and predict both future treatment and the final outcome $Y$, the CD4 count at year two. %Estimation of the effect of ART is perhaps the most salient example of application of structural nested models in the empirical literature.

Based on a histogram of the first measured CD4 count in the AIEDRP data (Acute Infection and Early Disease Research Program, see \cite{Hecht}), and based on the median and IQR estimates of the first measured CD4 count in \cite{Althoff2}, we choose to simulate $L_0$ so that the square root of $L_0$ is approximately normal. According to \cite{Althoff2}, ``the median CD4 count at presentation increased from $256$ cells/mm3 (interquartile range, $96–-455$ cells/mm3) to $317$ cells/mm3 interquartile range (IQR), $135-–517$ cells/mm3) from 1997 to 2007.''
For our first scenario, we choose $\sqrt{\tilde{L}_0}\sim {\cal N}(17,8^2)$. We simulate $e_k\sim {\cal N}(0,20^2)$, so it has a relatively small standard deviation. Based on a preliminary simulation with one million observations, the median $L_0$ in this simulation scenario is $294$, IQR $135-501$; these values are close to the empirical values.

In our simulations, the probability of treatment initiation in the first year is $0.70$ for patients with a baseline CD4 count below $350$ and $0.30$ for patients with a baseline CD4 count above $350$. While in clinical practice, patients with a higher CD4 count are less likely to be treated, with $350$ often used as a cut-off, the $0.70$ and $0.30$ values are not chosen based on data, because treatment guidelines have been changing considerably over time in the past few years and differ by country (\cite{panel4,WHO3}). The simulation values of the treatment initiation parameters ensure that all patients have a considerable probability of being untreated and also a considerable probability of being treated. If patients with specific covariates are either always treated or always untreated, we cannot estimate the effect of treatment for these patients (because it is impossible to distinguish the treatment effect from the reason why the treatment was given). Thus, in the simulations we choose $\lambda^{(0)}_0=-log(1-0.3)$ and $\lambda^{(0)}_1=-log(1-0.7)$. For treatment initiation during the second year, we choose the same parameter values: $\lambda_0^{(1)}=-log(1-0.3)$ and $\lambda_1^{(1)}=-log(1-0.7)$. This implies that for any untreated covariate history, the probability of ever initiating treatment is $.3+(.7\times.3)$ to $.7+(.3\times .7)$, or $0.51$ to $0.91$. Based on the estimates in \cite{LokAIDS} and \cite{Victor}, the median CD4 count could increase by about 200 between ART initiation and one year later, and the effect of one year of ART is about $12\times 25=300$. Therefore, we choose: $\beta_0=100$, and $\theta_0=\psi_0=300$. Table~1, setting 1 describes the results of this simulation scenario.

In a second simulation scenario, setting 2 in Table~1, we introduce more variation around the signals. In the second scenario, we increase the variance of $e_k$ to $200^2$, with $\sqrt{\tilde{L}_0}\sim {\cal N}(17,6^2)$. Based on a preliminary simulation with one million observations, the median $L_0$ in scenario~2 is $313$, IQR $125-501$. The other parameters are as in scenario~1. In the third simulation study, setting 3 in Table~1, we introduce even more variation around the signals. In the third scenario, we increase the variance of $e_k$ to $300^2$, with $\sqrt{\tilde{L}_0}\sim {\cal N}(18,3^2)$. Based on a preliminary simulation with one million observations, the median $L_0$ in scenario~3 is $333$, IQR $119-545$. The other parameters are as in scenario~1.
%I also tried: $300^2$ and ${\cal N}(17.5,3^2)$: median 315, IQR $102-527$ based on 500.000 patients. This is mimsim5. I choose 18, which is mimsim4.

% $e_k\sim N(0,300^2)$ $with $\sqrt{\tilde{L}_0}\sim {\cal N}(17,4^2)$: based on a simulation with 1 thousand patients, this lead to a median $L_0$ of $309$, IQR $91-525$. The bottom number is a little low.

\subsection{Calculating the infinitesimal shift function $D$}%Not read Aug 21. But worked through it Sep 18.
\label{Dder}

This section calculates $D$ for the simulation study. First note that if treatment had already started by time $t$, there is no difference between $Y^{(t)}$ and $Y^{(t+h)}$, so that $D(y,t; \overline{Z}_t)=0$ for all $y$. Therefore, we focus on calculating $D$ for $\overline{Z}_t$ such that $T>t$. Let $t$ be given. We only need to derive $F_{Y^{(t+h)}|\overline{Z}_t}(y)$ for $h>0$ small. Therefore, we restrict calculations to $h\in[t,t+h_0]$ such that $[t+h_0]$ lies within either $[0,1)$ or $[1,2)$, depending on which of these two intervals contains $t$. Let $e_2$ be any random variable which is independent of $\overline{Z}_t$ and the treatment process and which has the same distribution as the $e_{2,t}$. Denote the actual duration of treatment until time $t+h$ by the random variable $cum(t+h)$. We derive:
\begin{eqnarray}\label{FYth}
F_{Y^{(t+h)}|\overline{Z}_t}(y)
&=&P\left(Y^{(t+h)}\leq y|\overline{Z}_t\right)\nonumber\\
&=&P\left(\tilde{L}_0-2\beta_0+\psi\Bigl(cum(t+h)+2-(t+h)\Bigr)+e_{2,min(T,t+h)}\leq y|\overline{Z}_t\right)\nonumber\\
&=&P\left(\tilde{L}_0-2\beta_0+\psi\Bigl(cum(t+h)+2-(t+h)\Bigr)+e_{2}\leq y|\overline{Z}_t\right)\nonumber\\
&=&E\left[P(\tilde{L}_0-2\beta_0+\psi\Bigl(cum(t+h)+2-(t+h)\Bigr)+e_{2}\leq y|\overline{Z}_t,e_2,\tilde{L}_0)|\overline{Z}_t\right]\nonumber\\
&=&E\left[P(cum(t+h)\leq \frac{1}{\psi}(y-e_2-\tilde{L}_0+2\beta_0)-(2-t-h)|\overline{Z}_t,e_2,\tilde{L}_0)|\overline{Z}_t\right].
\end{eqnarray}

Next, since we restrict to $\overline{Z}_t$ such that $T>t$, we have that $0\leq cum(t+h)\leq h$, and because of the way the treatment initiation process was simulated, for $x\in[0,h]$,
\begin{eqnarray}\label{cum}
P\left(cum(t+h)\leq x|\overline{Z}_t,e_2,\tilde{L}_0\right)
&=&P\left(cum(t+h)\leq x|\overline{Z}_t\right)\nonumber\\
&=&P\left((t+h)-T\leq x|\overline{Z}_t\right)\nonumber\\
%&=&1-P\left(cum(t+h)> x|\overline{Z}_t\right)\nonumber\\
%&=&1-P\left((t+h)-T> x|\overline{Z}_t\right)\nonumber\\
%&=&1-P\left(T-t<+h-x|\overline{Z}_t\right)\nonumber\\
&=&P\left(T-t\geq h-x|\overline{Z}_t\right)\nonumber\\
&=&\left\{\begin{array}{ll} e^{-\lambda^{(\lfloor t\rfloor )}_0(h-x)} &{\rm if}\;L_{\lfloor t\rfloor}> c_{\lfloor t\rfloor }\\
                            e^{-\lambda^{(\lfloor t\rfloor )}_1(h-x)} &{\rm if}\;L_{\lfloor t\rfloor }\leq c_{\lfloor t\rfloor },
\end{array}\right.
\end{eqnarray}
where $\lfloor t\rfloor$ is the floor of $t$, the largest integer less than or equal to $t$. In the first line of (\ref{cum}) we use that, in our simulation design, the rate of treatment initiation does not depend on $(e_2,\tilde{L}_0)$, given $\overline{Z}_t$, which includes $L_{\lfloor t\rfloor }$.

For $T>t$, $\overline{L}_{\lfloor t\rfloor }\leq c_{\lfloor t\rfloor }$, and $\lambda_1=\lambda^{(\lfloor t\rfloor )}_1$, since $cum(t+h)\in[0,h]$ for patients with $T>t$, it follows from equations~(\ref{FYth}) and~(\ref{cum}) that
\begin{eqnarray}
F_{Y^{(t+h)}|\overline{Z}_t}(y)
&=&0\cdot P\left(\frac{1}{\psi}(y-e_2-\tilde{L}_0+2\beta_0)-(2-t-h)<0|\overline{Z}_t\right)\nonumber\\
&&+1\cdot P\left(\frac{1}{\psi}(y-e_2-\tilde{L}_0+2\beta_0)-(2-t-h)>h|\overline{Z}_t\right)\nonumber\\
&&+E\left[e^{-\lambda_1\left(h-\left(\frac{1}{\psi}(y-e_2-\tilde{L}_0+2\beta_0)-(2-t-h)\right)\right)}|\right.\nonumber\\
&&\left.|\frac{1}{\psi}(y-e_2-\tilde{L}_0+2\beta_0)-(2-t-h)\in (0,h),\overline{Z}_t\right]\nonumber\\
&&\cdot P\left(\frac{1}{\psi}(y-e_2-\tilde{L}_0+2\beta_0)-(2-t-h)\in (0,h)|\overline{Z}_t\right)\nonumber\\
&=& P\left((y-e_2-\tilde{L}_0+2\beta_0)-\psi(2-t)>0|\overline{Z}_t\right)\nonumber\\
&&+E\left[e^{-\lambda_1\left(-\frac{1}{\psi}(y-e_2-\tilde{L}_0+2\beta_0)+(2-t)\right)}\right.\nonumber\\
&&\left.|\frac{1}{\psi}(y-e_2-\tilde{L}_0+2\beta_0)-(2-t-h)\in (0,h),\overline{Z}_t\right]\nonumber\\
&&\cdot P\left((y-e_2-\tilde{L}_0+2\beta_0)-\psi(2-t)\in (-\psi h,0)|\overline{Z}_t\right)\nonumber\\
&=& P\left(Y^{(t)}< y|\overline{Z}_t\right)\nonumber\\
&&+E\left[e^{-\lambda_1\left(-\frac{1}{\psi}(y-e_2-\tilde{L}_0+2\beta_0)+(2-t)\right)}\right.\nonumber\\
&&\left.|\frac{1}{\psi}(y-e_2-\tilde{L}_0+2\beta_0)-(2-t-h)\in (0,h),\overline{Z}_t\right]\nonumber\\
&&\cdot P\left((y-Y^{(t)})\in (-\psi h,0)|\overline{Z}_t\right)\nonumber\\
&=& P\left(Y^{(t)}\leq y|\overline{Z}_t\right)\nonumber\\
&&+E\left[e^{-\lambda_1\left(-\frac{1}{\psi}(y-e_2-\tilde{L}_0+2\beta_0)+(2-t)\right)}\right.\nonumber\\
&&\left.|\frac{1}{\psi}(y-e_2-\tilde{L}_0+2\beta_0)-(2-t-h)\in (0,h),\overline{Z}_t\right]\nonumber\\
&&\cdot P\left(Y^{(t)}\in (y,y+\psi h)|\overline{Z}_t\right)\nonumber\\
&=& F_{Y^{(t)}|\overline{Z}_t}(y)\nonumber\\
&&+E\left[e^{-\lambda_1\left(-\frac{1}{\psi}(y-e_2-\tilde{L}_0+2\beta_0)+(2-t)\right)}\right.\nonumber\\
&&\left.|\frac{1}{\psi}(y-e_2-\tilde{L}_0+2\beta_0)-(2-t)\in (-h,0),\overline{Z}_t\right]\nonumber\\
&&\cdot \left(F_{Y^{(t)}|\overline{Z}_{t}}(y+\psi h)-F_{Y^{(t)}|\overline{Z}_t}(y)\right).\label{mess}
\end{eqnarray}
In our simulation study, $(\tilde{L}_0,e_2)$ has a continuous conditional distribution $f_{(\tilde{L}_0,e_2)|\overline{Z}_t}$. Therefore, conditional on a value of $\overline{Z}_t$ such that $T>t$, we have that
\begin{eqnarray*}
\lefteqn{E\left[e^{-\lambda_1\left(-\frac{1}{\psi}(y-e_2-\tilde{L}_0+2\beta_0)+(2-t)\right)}|\frac{1}{\psi}(y-e_2-\tilde{L}_0+2\beta_0)-(2-t)\in (-h,0),\overline{Z}_t\right]}\\
&=&\int_{-\infty}^\infty d\tilde{l}_0\int_{-(2-t)\psi+y-\tilde{l}_0+2\beta_0}^{(-(2-t)+h)\psi+y-\tilde{l}_0+2\beta_0}de_2 f_{(\tilde{L}_0,e_2)|\overline{Z}_t}(\tilde{l}_0,e_2)e^{-\lambda_1\left(-\frac{1}{\psi}(y-e_2-\tilde{L}_0+2\beta_0)+(2-t)\right)},
\end{eqnarray*}
which is continuously differentiable in $h\geq 0$ for $(y,t,\overline{Z}_t)$ fixed, with some derivative, $g(y,t,h,\overline{Z}_t)$.
%where the first part can be proven by conditioning first on $(\tilde{L}_0,e_2,\overline{Z}_t)$, and if we assume for example that $\tilde{L}_0+e_2$ %has a continuous distribution given $\overline{Z}_t$ if $\overline{Z}_t$ indicates no treatment at time $t$.
Therefore, equation~(\ref{mess}) implies that
\begin{eqnarray*}
\frac{\partial}{\partial h} F_{Y^{(t+h)}|\overline{Z}_t}(y)
&=&
g(y,t,h,\overline{Z}_t)\left(F_{Y^{(t)}|\overline{Z}_{t}}(y+\psi h)-F_{Y^{(t)}|\overline{Z}_t}(y)\right)\\
&&+E\left[e^{-\lambda_1\left(-\frac{1}{\psi}(y-e_2-\tilde{L}_0+2\beta_0)+(2-t)\right)}\right.\\
&&\left.|\frac{1}{\psi}(y-e_2-\tilde{L}_0+2\beta_0)-(2-t)\in (-h,0),\overline{Z}_t\right]\\
&&\cdot \psi f_{Y^{(t)}|\overline{Z}_{t}}(y+\psi h),
\end{eqnarray*}
where $f_{Y^{(t)}|\overline{Z}_{t}}$ is the density of $Y^{(t)}$ given $\overline{Z}_t$. Notice that since $T>t$, $Y^{(t)}=Y_t$ given $\overline{Z}_t$, and the density $f_{Y^{(t)}|\overline{Z}_{t}}=f_{Y_t|\overline{Z}_{t}}$ exists and is continuous.
Letting $h\downarrow 0$, it follows that
\begin{equation*}
\left.\frac{\partial}{\partial h}\right|_{h=0} F_{Y^{(t+h)}|\overline{Z}_t}(y)
=\psi f_{Y^{(t)}|\overline{Z}_t}(y).
\end{equation*}
Clearly,
\begin{equation*}
\frac{\partial}{\partial y} F_{Y^{(t)}|\overline{Z}_t}(y)=f_{Y^{(t)}|\overline{Z}_t}(y).
\end{equation*}
Hence, because of equation~(12), %\ref{Deq}
\begin{equation}\label{Dpar}
D(y,t;\overline{Z}_t)=-\psi 1_{{\rm untreated}\;{\rm at}\;t}.
\end{equation}
The same derivation can be used for $L_{\lfloor t\rfloor }>c_{\lfloor t\rfloor }$.

In fact, it can be shown that if the counterfactual covariates and the counterfactual outcomes are as in this simulation study, $D=-\psi 1_{{\rm untreated}\;{\rm at}\;t}$ if $T$ has a continuous conditional density $f_{T|\overline{Z}_t}(y)$ given $\overline{Z}_t$ for $y\in[t,t+h_0]$ for some $h_0>0$. This is beyond the scope of the current article. In future work, we also plan to address multidimensional $\psi$, as well as parameterizing $D$ for survival outcomes.

\subsection{Calculating $X_{\psi}(t)$}\label{Xpsisec}

As a consequence of Section~4 and equation~(\ref{Dpar}), it follows that $X_{\psi}(t)$ is the solution to
\begin{equation*}
dX_{\psi}(t)/dt=-\psi 1_{{\rm untreated}\;{\rm at}\;t}
\end{equation*}
with end condition $X_{\psi}(2)=Y$ (recall that time $2$ is the time the outcome is measured). Therefore,
\begin{equation*}
X_{\psi}(t)=Y+\psi (min(T,2)-t)1_{T>t},
\end{equation*}
where $(min(T,2)-t)1_{T>t}$ is the duration of the patient not being on treatment between time $t$ and time $2$.%Checked it is indeed the solution to the differential equation Sep 15.

\subsection{Estimating equations when treatment initiation follows a piecewise exponential model}\label{est}

Suppose that we know that treatment initiation follows a piecewise exponential model, with parameter depending on a discretized covariate, measured at time $0$ and time $1$. In the simulation, we assume
\begin{equation*}
\lambda_T\left(t|\overline{Z}_t\right)
=\left\{\begin{array}{ll}
\lambda^{(0)}_0 & {\rm if}\; L_0> c_0 \;{\rm and}\;t\in[0,1]\\
\lambda^{(0)}_1 & {\rm if}\; L_0\leq c_0 \;{\rm and}\;t\in[0,1]\\
\lambda^{(1)}_0 & {\rm if}\; L_1> c_1 \;{\rm and}\;t\in(1,2]\\
\lambda^{(1)}_1 & {\rm if}\; L_1\leq c_1 \;{\rm and}\;t\in(1,2]
\end{array}\right.
\end{equation*}
for known constants $c_0$ and $c_1$ in $\mathbb{R}$, and for $\lambda^{(0)}_0$, $\lambda^{(0)}_1$, $\lambda^{(1)}_0$, and $\lambda^{(1)}_1$ unknown values in $\mathbb{R}$. To select from the many possible estimating equations for $\psi$ provided in Theorem~5.2, %\ref{thec},
we follow the approach of \cite{Enc}, which was proved to lead to consistent estimation in \cite{plss} provided the main result of the current article holds true. Below we explain why this approach works in the context of this simulation study. \cite{Enc} proposed to add $\alpha$ times a function of $X_\psi(t)$ and $\overline{Z}_{t-}$ to the model for treatment initiation $\lambda_T$, and find the parameter $\psi$ such that adding this function has no effect on the estimated hazard (that is, the $\psi$ that leads to $\hat{\alpha}=0$); that particular $\psi$ will be the estimate $\hat{\psi}$. The underlying observation for this procedure is that given $\overline{Z}_{t-}$, for the true $\psi$, $X_\psi(t)$ should not help to predict treatment changes (Section~5). Adding $\alpha$ times a function of $X_\psi(t)$ and $\overline{Z}_{t-}$ to $\lambda_T$ can be done in many different ways. For simplicity of calculations, we choose to add $\alpha X_\psi(0)$, a function of $X_\psi(t)$ and $\overline{Z}_{t-}$, to the model for treatment initiation $\lambda_T$ in the time interval $[0,1]$, and $\alpha X_\psi(1)$ to the model for treatment initiation $\lambda_T$ in the time interval $[1,2]$, both in a way similar to a Cox proportional hazards component: as a factor $e^{\alpha X_{\psi}(0)}$ and $e^{\alpha X_{\psi}(1)}$, respectively. Let $\delta_i^{(0)}=1_{T_i\leq 1}$, $\delta_i^{(1)}=1_{1<T_i\leq 2}$, $Z_i(0)=1_{L_0\leq c_0}$, and $Z_i(1)=1_{L_1\leq c_1}$. The partial likelihood for the model extended with $X_\psi$ as described above is
\begin{eqnarray*}
L\left(\lambda^{(0)}_0,\lambda^{(0)}_1,\lambda^{(1)}_0,\lambda^{(1)}_1,\alpha\right)
&=&\prod_{i=1}^n \left(\lambda^{(0)}_{Z_i(0)}e^{\alpha X_{\psi,i}(0)}\right)^{\delta_i^{(0)}} e^{-\lambda^{(0)}_{Z_i(0)}e^{\alpha X_{\psi,i}(0)}min(T_i,1)} \\
&&\prod_{i=1}^n \left(\lambda^{(1)}_{Z_i(1)}e^{\alpha X_{\psi,i}(1)}\right)^{\delta_i^{(1)}}e^{-\lambda^{(1)}_{Z_i(1)}\left(1-\delta_i^{(0)}\right)e^{\alpha X_{\psi,i}(1)}(min(T_i,2)-1)}.
\end{eqnarray*}
The log likelihood is:
\begin{eqnarray*}
\lefteqn{log L\left(\lambda^{(0)}_0,\lambda^{(0)}_1,\lambda^{(1)}_0,\lambda^{(1)}_1,\alpha\right)}\\
&=&\sum_{i=1}^n \delta_i^{(0)}log\left(\lambda^{(0)}_{Z_i(0)}e^{\alpha X_{\psi,i}(0)}\right)-\lambda^{(0)}_{Z_i(0)}min(T_i,1)e^{\alpha X_{\psi,i}(0)}\\
&& + \delta_i^{(1)}log\left(\lambda^{(1)}_{Z_i(1)}e^{\alpha X_{\psi,i}(1)}\right)-\left(1-\delta_i^{(0)}\right) \lambda^{(1)}_{Z_i(1)}(min(T_i,2)-1)e^{\alpha X_{\psi,i}(1)}\\
&=&\sum_{i=1}^n \delta_i^{(0)}Z_i(0)log\lambda^{(0)}_{1}-Z_i(0)\lambda^{(0)}_{1}(min(T_i,1))e^{\alpha X_{\psi,i}(0)}\\
&& +\sum_{i=1}^n\delta_i^{(0)}(1-Z_i(0))log\lambda^{(0)}_{0}-(1-Z_i(0))\lambda^{(0)}_{0}min(T_i,1)e^{\alpha X_{\psi,i}(0)}\\
&&+ \sum_{i=1}^n\delta_i^{(1)}Z_i(1)log\lambda^{(1)}_{1}- \left(1-\delta_i^{(0)}\right)Z_i(1)\lambda^{(1)}_{1}(min(T_i,2)-1)e^{\alpha X_{\psi,i}(1)}\\
&&+ \sum_{i=1}^n \delta_i^{(1)}(1-Z_i(1))log\lambda^{(1)}_{0}- \left(1-\delta_i^{(0)}\right)(1-Z_i(1))\lambda^{(1)}_{0}(min(T_i,2)-1)e^{\alpha X_{\psi,i}(1)}\\
&&+\sum_{i=1}^n \delta_i^{(0)}\alpha X_{\psi,i}(0)+ \delta_i^{(1)}\alpha X_{\psi,i}(1).
\end{eqnarray*}
Following \cite{plss}, to calculate $\hat{\psi}$, we take the derivative of this expression with respect to\linebreak $(\lambda^{(0)}_0,\lambda^{(0)}_1,\lambda^{(1)}_0,\lambda^{(1)}_1,\alpha)$, then set $\alpha=0$, and solve for $(\hat{\lambda}^{(0)}_0,\hat{\lambda}^{(0)}_1,\hat{\lambda}^{(1)}_0,\hat{\lambda}^{(1)}_1,\hat{\psi})$; as indicated below, consistency of the estimator will follow from Theorem~5.2. %\ref{thec}).
We obtain the estimating equations:
\begin{eqnarray}\label{EE}
\lefteqn{0=\left.\frac{\partial}{\partial \lambda^{(0)}_0,\lambda^{(0)}_1,\lambda^{(1)}_1,\lambda^{(1)}_0,\alpha}\right|_{\alpha=0}log L\left(\lambda^{(0)}_0,\lambda^{(0)}_1,\lambda^{(1)}_0,\lambda^{(1)}_1,\alpha\right)}\nonumber\\
&=&\sum_{i=1}^n \left(\begin{array}{c}
\frac{\delta_i^{(0)}(1-Z_i(0))}{\lambda^{(0)}_{0}}-(1-Z_i(0))min(T_i,1)\\
\frac{\delta_i^{(0)}Z_i(0)}{\lambda^{(0)}_{1}}-Z_i(0)min(T_i,1)\\
\frac{\delta_i^{(1)}(1-Z_i(1))}{\lambda^{(1)}_{0}}- \left(1-\delta_i^{(0)}\right)(1-Z_i(1))(min(T_i,2)-1)\\
\frac{\delta_i^{(1)}Z_i(1)}{\lambda^{(1)}_{1}}- \left(1-\delta_i^{(0)}\right)Z_i(1)(min(T_i,2)-1)\\
G_i(\psi)
\end{array}\right),
\end{eqnarray}
with
\begin{eqnarray}
G_i(\psi)&=&-\left(Z_i(0)\lambda^{(0)}_{1}+(1-Z_i(0))\lambda^{(0)}_{0}\right)(min(T_i,1))X_{\psi,i}(0)\nonumber\\
&&-\left(Z_i(1)\lambda^{(1)}_{1}+(1-Z_i(1))\lambda^{(1)}_{0}\right)\left(1-\delta_i^{(0)}\right)(min(T_i,2)-1)X_{\psi,i}(1)\nonumber\\
&&+\delta_i^{(0)} X_{\psi,i}(0)+ \delta_i^{(1)} X_{\psi,i}(1).\label{psiEE}
\end{eqnarray}
According to Theorem~5.2, with $h_t(X(t),\overline{Z}_{t-})=X(0)$ for $t\in[0,1]$ and $h_t(X(t),\overline{Z}_{t-})=X(1)$ for $t\in(1,2]$, %\ref{thec}),
these are indeed unbiased estimating equations. Solving the estimating equations, in the first step, the $\lambda$'s are estimated by their maximum likelihood estimates without adding $X_\psi$ to the model:
\begin{equation}\label{hat}
\left(\begin{array}{c}
\hat{\lambda}^{(0)}_0\\
\hat{\lambda}^{(0)}_1\\
\hat{\lambda}^{(1)}_0\\
\hat{\lambda}^{(1)}_1
\end{array}\right)
=\left(\begin{array}{c}
\frac{\sum_{i=1}^n\delta_i^{(0)}(1-Z_i(0))}{\sum_{i=1}^n(1-Z_i(0))min(T_i,1)}\\
\frac{\sum_{i=1}^n\delta_i^{(0)}Z_i(0)}{\sum_{i=1}^nZ_i(0)min(T_i,1)}\\
\frac{\sum_{i=1}^n\delta_i^{(1)}(1-Z_i(1))}{\sum_{i=1}^n(1-Z_i(1))1_{T_i>1}(min(T_i,2)-1)}\\
\frac{\sum_{i=1}^n\delta_i^{(1)}Z_i(1)}{\sum_{i=1}^nZ_i(1)1_{T_i>1}(min(T_i,2)-1)}
\end{array}\right).
\end{equation}
General theory, see e.g.~\cite{Vaart}, says these estimates for the hazard are consistent and asymptotically normal.
In the second step, $\psi$ is then estimated by plugging these estimates for the hazard in $\sum_{i=1}^n G_i(\psi)$, with $G_i(\psi)$ as in equation~(\ref{psiEE}), and solving for $\psi$. In the simulations, $X_\psi$ is linear in $\psi$, see Section~\ref{Xpsisec}, and thus $\sum_{i=1}^n G_i(\psi)$ is also linear in $\psi$. Therefore, solving for $\psi$ requires solving a linear, in our case even one-dimensional, equation. If the coefficient before $\psi$ is non-zero, there is a unique solution $\hat{\psi}$. In addition, the expectation of $G_i(\psi)$ at the true $\lambda$ is linear in $\psi$, so $EG_i(\psi)=0$ has a unique solution if the coefficient on $\psi$ in the linear equation $EG_i(\psi)$ is non-zero, which guarantees consistency and asymptotic normality (\cite{Vaart} Chapter 5). We conclude from Section~\ref{Xpsisec} and equation~(\ref{psiEE}) that a consistent estimator of $\psi$ can be defined as $\hat{\psi}=- \sum_{i=1}^n A_{1i}/ \sum_{i=1}^n A_{2i}$, where
\begin{eqnarray*}
A_{1i}&=&-Y_i\left(Z_i(0)\hat{\lambda}^{(0)}_{1}+(1-Z_i(0))\hat{\lambda}^{(0)}_{0}\right)min(T_i,1)\nonumber\\
&&-Y_i\left(Z_i(1)\hat{\lambda}^{(1)}_{1}+(1-Z_i(1))\hat{\lambda}^{(1)}_{0}\right)\left(1-\delta_i^{(0)}\right)(min(T_i,2)-1)\nonumber\\
&&+Y_i\delta_i^{(0)}+ Y_i\delta_i^{(1)}
\end{eqnarray*}
and
\begin{eqnarray*}
A_{2i}&=&-\left(Z_i(0)\hat{\lambda}^{(0)}_{1}+(1-Z_i(0))\hat{\lambda}^{(0)}_{0}\right)min(T_i,1)min(T_i,2)\nonumber\\
&&-\left(Z_i(1)\hat{\lambda}^{(1)}_{1}+(1-Z_i(1))\hat{\lambda}^{(1)}_{0}\right)\left(1-\delta_i^{(0)}\right)(min(T_i,2)-1)^2\nonumber\\
&&+\delta_i^{(0)}min(T_i,2) + \delta_i^{(1)}(min(T_i,2)-1).
\end{eqnarray*}

We choose to add $X_{\psi}(0)$ to the prediction model for treatment changes in the time interval $[0,1]$ and to add $X_\psi(1)$ in the time interval $(1,2]$. Optimally choosing the function of $X_\psi(t)$ and $\overline{Z}_{t-}$ to add to the prediction model for treatment changes is an interesting topic for future research.

Algebra shows that in this simulation study, the bias and the MSE of $\hat{\psi}$ do not depend on the values of $\psi_0$ or $\theta_0$ (for $\theta_0$ this is easily seen by noticing that the estimators depend only on the $Y_i$, the $T_i$ (which depend only on pre-treatment variables), and pre-treatment variables). These algebraic calculations were confirmed by simulating scenario 2 with $\theta_0=\psi_0=300$ replaced by $\theta_0=\psi_0=100$, which lead to the same bias and MSE as scenario 2 itself. %This was mimsim3
 Therefore, we did not vary $\psi_0$ or $\theta_0$ in the simulation study.

\subsection{Results of the simulation study}\label{res}

We ran a simulation study with $n=500$, $1000$, $2000$, $5000$, and $10000$, and 5000 repetitions each. The results for the three settings described in Section~\ref{HIV} are presented in Section~\ref{sim}, Table~1.

In this simulation study, both for small and large samples, the bias of the estimators is small. In all three settings and for all sample sizes considered (including the small sample size $n=100$), the MSE of the estimators arises mostly from the variance, not from the bias. Also, if the true parameter $\psi$ equals $300$ as in this simulation study, for $n=500$, $\sqrt{MSE}/\psi=0.04$ in setting~1, and~$0.08$ in setting~3. Thus, the estimates are already precise in relatively small samples. Because, as we noted before, the MSE in this simulation study does not depend on the true parameter $\psi$, a larger sample size would be required to obtain precise estimators of small true parameter values $\psi$.

Table~1 also shows that the mean squared error decreases appropriately as the sample size increases. In all three settings, the MSE times the sample size is roughly constant. This indicates that the large-sample theory in Section~\ref{est} (which follows from the fact that the estimating equations are unbiased) provides a reasonable approximation for the rate of convergence of the estimator in finite samples.

We conclude that adding $X_{\psi}(0)$ to the prediction model for treatment changes in the time interval $[0,1]$ and adding $X_\psi(1)$ in the time interval $(1,2]$, as described in Section~\ref{est}, provides estimators with good finite-sample properties in this simulation study.

Finally, to confirm consistency of the estimators, we ran a simulation study with one dataset and $n$ equal to one million, in all three settings. This resulted in $\hat{\psi}=299.545$ for setting 1, $\hat{\psi}=299.863$ for setting 2, and $\hat{\psi}=299.673$ for setting 3, all very close to the true value of $\psi=300$ in the simulation study.
%I ran n=100, psi-hat=302. I may not have enough variability built in! Perhaps the variance of the $e_i$ needs to go up! Let's run $n=500$. Ah, perhaps it is OK: psi-hat=292.

We conclude that in this simulation study, continuous-time structural nested models perform extremely well.

\section{Web-Appendix: Some facts about conditioning}\label{Cond}

The following definition and two theorems on existence and uniqueness
of conditional distributions can be found in
\citet{Bau} Section 10.3, in a different formulation. The first is a
definition of conditional distributions. A conditional distribution of
a random variable $X$ with values in ${\mathbb R}$ is more than just a
set of
conditional probabilities $P\left(X\leq x|{\cal G}\right)$ for
$x\in{\mathbb R}$: it is also a probability measure on ${\mathbb
R}$. Conditional probabilities always exist; a conditional
distribution always exists e.g.\ if $X$ takes values in ${\mathbb R}$, but not
in general. Conditional probabilities are
almost surely unique; under conditions the same is true for
conditional distributions.

\begin{defi}
\label{cdist} Let $X:\left(\Omega, {\cal F}\right)\rightarrow \left({\cal
X},{\cal A}\right)$ be a random variable on a probability space
$\left(\Omega, {\cal F},P\right)$ with values in a measurable space
$\left({\cal X},{\cal A}\right)$. Let ${\cal G}\subset {\cal F}$ be a
sub-$\sigma$-algebra. Then $P_{X|{\cal G}}:\Omega\times {\cal
A}\rightarrow {\mathbb R}$ is a conditional distribution of $X$ given
${\cal G}$ if
\begin{enumerate}[a)]
\item $\forall A\in{\cal A}$: $\omega\rightarrow P_{X|{\cal
G}}\left(\omega,A\right)$
is a version of $P\left(X\in A|{\cal G}\right)$, i.e. it is ${\cal
G}$-measurable and $\forall G\in{\cal G}$:
\begin{equation*}
\int_{G}P_{X|{\cal G}}\left(\omega,A\right)dP\left(\omega\right)
=\int_{G}1_{A}\left(X\right)dP=P\left(G\cap
X^{-1}\left(A\right)\right).\end{equation*}
\item $P_{X|{\cal G}}\left(\omega,\cdot\right)$ is a probability
measure on $\left({\cal X},{\cal A}\right)$.
\end{enumerate}
\end{defi}
\noindent If $X$ takes values in ${\mathbb R}$ the distribution function
belonging to the probability measure $P_{X|{\cal G}}$ is often
denoted by $F_{X|{\cal G}}$.

\begin{thm} Let $X:\left(\Omega, {\cal F}\right)\rightarrow \left({\cal
X},{\cal A}\right)$ be a random variable on a probability space
$\left(\Omega, {\cal F},P\right)$ with values in a measurable space
$\left({\cal X},{\cal A}\right)$. Suppose that ${\cal A}$ is a
countably generated $\sigma$-algebra and ${\cal G}\subset {\cal F}$ is a
sub-$\sigma$-algebra.
If $P_{X|{\cal G}}$ and $P^*_{X|{\cal G}}$ are
two conditional distributions of $X$ given ${\cal G}$ then they are
almost surely the same in the sense that there exists a $P$-null set
$N\in {\cal F}$ such that for all
$\omega\in\Omega\setminus N$ and all $A\in {\cal A}$,
\begin{equation*}P_{X|{\cal G}}\left(\omega,A\right)=P^*_{X|{\cal
G}}\left(\omega,A\right).\end{equation*}
\end{thm}

A topological space $E$ is called Polish if it has a countable dense
subset and there exists a metric that generates the topology and for
which the space is complete.
An example of a Polish space is ${\mathbb R}^k$ with the usual
topology.
\begin{thm} \label{Pool} Let $X:\left(\Omega,{\cal F}\right)\rightarrow
\left(E,{\cal B}\left(E\right)\right)$ be a random variable on a
probability space
$\left(\Omega,{\cal F},P\right)$ with values in a Polish space
$E$ with its Borel-$\sigma$-algebra.
Then for every $\sigma$-algebra ${\cal G}\subset {\cal F}$ there
exists a conditional distribution $P_{X|{\cal G}}$.
\end{thm}
The next theorem is very useful in combination with
Theorem~\ref{Pool}. Suppose that $Z$ is a random variable on
$\left(\Omega,{\cal F}\right)$ with values in the space of cadlag
functions on $[a,b]$, $D[a,b]$, equipped with the $\sigma$-algebra
generated by the coordinate projections. Then Theorems~\ref{Pool}
and~\ref{Dpools} imply that for any $\sigma$-algebra ${\cal G}\subset
{\cal F}$, $Z$ has a conditional distribution given ${\cal G}$.
\begin{thm} \label{Dpools} Suppose that $a,b\in{\mathbb R}$ are finite. Then
$D\left[a,b\right]$ with the Skorohod topology is a Polish
space. Furthermore, the $\sigma$-algebra on $D\left[a,b\right]$
generated by the Skorohod topology is the same as the $\sigma$-algebra
on $D\left[a,b\right]$ generated by the coordinate projections.
\end{thm}
\noindent The first statement of this theorem can be found in \citet{Bill68},
Chapter~3, the second statement is Theorem~14.5 in the same book.

The next lemma is an easy consequence of the existence of
conditional distributions:

\begin{lem} \label{Qok} Let $X$ and $Y$ be random variables on a probability
space
$\left(\Omega, {\cal F},P\right)$ with values in
$\left({\mathbb R},{\cal B}\right)$, with ${\cal B}$ the Borel-$\sigma$-algebra
on ${\mathbb R}$. Suppose that ${\cal G}\subset {\cal F}$ is a
sub-$\sigma$-algebra. Then there exist conditional distributions
$P_{X|{\cal
G}}$ and $P_{Y|{\cal G}}$. If moreover for every $x\in {\mathbb Q}$,
$P\left(X\leq x|{\cal G}\right)=P\left(Y\leq x|{\cal G}\right)$ a.s.,
then $P_{X|{\cal G}}=P_{Y|{\cal G}}$ a.s.\ in the sense that there
exists a $P$-null set $N\in{\cal F}$ such that for all $\omega
\in\Omega\setminus N$ and all $B\in{\cal B}$,
\begin{equation*}P_{X|{\cal G}}\left(\omega,B\right)=P_{Y|{\cal
G}}\left(\omega,B\right).\end{equation*}
\end{lem}

\noindent{\bf Proof.}  Existence of conditional distributions follows
from Theorem \ref{Pool} since $\left({\mathbb R},{\cal B}\right)$ is a
Polish space. Furthermore a probability measure on $\left({\mathbb
R},{\cal B}\right)$ is completely determined by its values on
$\left(-\infty, x\right]$ for $x\in{\mathbb Q}$. So it is enough to
prove that there exists a $P$-null set $N\in{\cal F}$ such that
\begin{equation} \label{nul}\forall \omega\in\Omega\setminus N \;\;\;\;\;\;\;\;\;
\forall x\in{\mathbb
Q}\;\;\;\;\;\;\;\;\; P_{X|{\cal
G}}\left(\omega,\left(-\infty,x\right]\right) = P_{Y|{\cal
G}}\left(\omega,\left(-\infty,x\right]\right).\end{equation}
But for every $x\in{\mathbb Q}$,
\begin{eqnarray*} P_{X|{\cal G}}\left(\omega,\left(-\infty,x\right]\right)
&=&P\left(X\leq x|{\cal G}\right)\;\mbox{\rm a.s.}\\
&=&P\left(Y\leq x|{\cal G}\right)\;\mbox{\rm a.s.}\\
&=&P_{Y|{\cal G}}\left(\omega,\left(-\infty,x\right]\right)\;\mbox{\rm a.s.}
\end{eqnarray*}
Define
\begin{equation*}N=\cup_{x\in{\mathbb Q}}\left\{\omega:P_{X|{\cal
G}}\left(\omega,\left(-\infty,x\right]\right) \neq P_{Y|{\cal
G}}\left(\omega,\left(-\infty,x\right]\right)\right\}.\end{equation*}
This is a countable union of null sets, so a null set, and it
satisfies (\ref{nul}).
\hfill $\Box$

$ $

For the proof of Lemma~\ref{un} and Lemma~\ref{gZ}, the
following two lemma's are used.
The first is well-known.
\begin{lem}\label{exp} If $P_{X|{\cal G}}$ is a conditional distribution of
$X$ given
${\cal G}$ then
\begin{equation*}E\left[f\left(X\right)|{\cal G}\right]=\int
f\left(x\right)dP_{X|{\cal G}}\left(x\right)\hspace{0.5cm}\;\mbox{\rm a.s.}.\end{equation*}
\end{lem}
\begin{lem}\label{concz} Suppose that $Z$ and $Y$ are random variables on a
probability
space $\left(\Omega,{\cal F},P\right)$ with values in Polish spaces
$\left({\mathcal Y},{\cal A}_1\right)$
and $\left({\cal Z},{\cal A}_2\right)$, respectively. Then
\begin{equation*}\left(\omega,A\right)\longmapsto
\int_{A}\delta_{z',z\left(\omega\right)}
P_{Y|Z=z\left(\omega\right)}\left(dy\right)dz'\end{equation*}
$:\Omega\times\sigma\left({\cal A}_1\times {\cal
A}_2\right)\rightarrow {\mathbb
R}$ is a version of $P_{\left(Y,Z\right)|Z}$, i.e.\ it is a
conditional distribution function of $\left(Y,Z\right)$ given $Z$.
\end{lem}

\noindent{\bf Proof.}
Define $\tilde{P}\left(\omega,A\right)=\int_{A}\delta_{z',z\left(\omega\right)}
P_{Y|Z=z\left(\omega\right)}\left(dy\right)dz'$.
Condition a and b of Definition~\ref{cdist} have to be checked for $\tilde{P}$.
Condition b:
for $\omega$ fixed it is indeed a probability measure on
$\left(\left({\mathcal Y}\times{\cal Z}\right),\sigma\left({\cal A}_1\times{\cal A}_2\right)\right)$ (concentrated on
$z=z\left(\omega\right)$).

Condition a: first it is shown that for any $A$ of the form $A_1\times
A_2$ with $A_1\in{\cal A}_1$ and
$A_2\in{\cal A}_2$, $\omega\rightarrow
\tilde{P}\left(\omega,A_1\times A_2\right)$ is a version of
$P\left(\left(Y,Z\right)\in\left(A_1\times A_2\right)|Z\right)$. Equivalently, for $A$ of the form $A_1\times
A_2$ and $G\in\sigma\left(Z\right)$, so $G$ of the form
$Z^{-1}\left(B\right)$ with $B\in{\cal A}_2$,
\begin{equation*}\int_{G}\tilde{P}\left(\omega,A\right)dP\left(\omega\right)
=\int_{G}1_{A}\left(\left(Y,Z\right)\right)dP=P\bigl(G\cap
\left(Y,Z\right)^{-1}\left(A\right)\bigr).\end{equation*}
This can be shown as follows:
\begin{eqnarray*} \int_{Z^{-1}\left(B\right)}
\tilde{P}\left(\omega,A_1\times A_2\right)dP\left(\omega\right)
&=&\int_{B}
\left(\int_{A_1\times A_2}\delta_{z',z}
P_{Y|Z=z}\left(dy'\right)dz'\right)
dP_{Z}\left(z\right)\\
&=&\int_{B\cap A_2}
\left(\int_{A_1}
P_{Y|Z=z}\left(dy'\right)\right)
dP_{Z}\left(z\right)\\
&=&\int_{B\cap A_2}
P\left(Y\in A_1|Z=z\right)
dP_{Z}\left(z\right)\\
&=&P\left(Y^{-1}\left(A_1\right)\cap Z^{-1}\left(B\cap
A_2\right)\right)\\
&=&P\bigl(Z^{-1}\left(B\right)\cap\left(Y,Z\right)^{-1}\left(A_1\times A_2\right)\bigr).
\end{eqnarray*}

Next it is shown that this is sufficient. Notice first that since
$\left(\left({\cal Y}\times {\cal Z}\right),\sigma\left({\cal
A}_1\times{\cal A}_2\right)\right)$ is a Polish space, there exists a
conditional distribution $P_{\left(Y,Z\right)|Z}$. We show that
$\tilde{P}$ and $P_{\left(Y,Z\right)|Z}$ are almost surely equal,
using the Uniqueness Theorem on page 27 of \citet{Bau}. Remark that
both ${\cal A}_1$ and ${\cal A}_2$ are countably generated, say by
${\cal A}_1^0$ and ${\cal A}_2^0$, so that $\sigma\left({\cal
A}_1\times{\cal A}_2\right)$ is countably generated by ${\cal
A}_1^0\times {\cal A}_2^0$: every $A_1\times A_2$ is an element of
$\sigma\bigl({\cal A}_1^0\times{\cal A}_2^0\bigr)$, since $A_1\times
A_2=A_1\times {\cal Z}\cap {\cal Y}\times A_2$. To apply the
Uniqueness Theorem we need a generator which is intersection-stable
(i.e., finite intersections of elements in ${\cal A}_1^0\times{\cal
A}_2^0$ are still in ${\cal A}_1^0\times{\cal A}_2^0$). ${\cal
A}_1^0\times{\cal A}_2^0$ need not be intersection-stable, but as in
the proof of Theorem 10.3.4 in \citet{Bau}: when finite intersections
of elements in ${\cal A}_1^0\times{\cal A}_2^0$ are added to the
generator it stays countable. Notice that these finite intersections are
still of the form $A_1\times A_2$ since $\left(A_1\times
A_2\right)\cap \left(B_1\times B_2\right)=\left(A_1\cap
B_1\right)\times\left(A_2\times B_2\right)$, and notice moreover that
this leads to a countable intersection-stable generator ${\cal
A}_1^1\times{\cal A}_2^1$.  Because of the former paragraph, for all
$A_1\times A_2$ with $A_1\in{\cal A}_1$ and $A_2\in{\cal A}_2$,
$\tilde{P}\left(\omega,A_1^1\times A_2^1\right)$ is a version of
$P\left(\left(Y,Z\right)\in A_1^1\times A_2^1|Z\right)$, so
$\tilde{P}\left(\omega,A_1^1\times
A_2^1\right)=P_{\left(Y,Z\right)|Z}\left(\omega,A_1^1\times
A_2^1\right)$ a.s. Hence because of the countability
\begin{equation*}
\cup_{A_1^1\times A_2^1: A_1^1\in{\cal A}_1^1, A_2^1\in{\cal A}_2^1}
\;\;\bigl\{\omega:\tilde{P}\left(\omega,A_1^1\times A_2^1\right)\neq
P_{\left(Y,Z\right)|Z}\left(\omega,A_1^1\times A_2^1\right)\bigr\}
\end{equation*}
is a null set.
% remark that since
%$\left(\left({\mathcal Y}\times{\cal Z}\right),\sigma\left({\cal
%A}_1\times{\cal A}_2\right)\right)$ is a Polish space there exists a
%conditional distribution $P_{\left(Y,Z\right)|Z}$.
Thus the Uniqueness Theorem on page 27 of \citet{Bau} implies that
indeed $\tilde{P}$ and $P_{\left(Y,Z\right)|Z}$ are equal except for on
this null set.  \hfill $\Box$

\begin{lem} \label{un}
Suppose that $Y$ has a continuous conditional distribution function
$F_{Y|Z}$ given $Z$. Then $F_{Y|Z}\left(Y\right)$ is uniformly
distributed on $\left[0,1\right]$ and independent of $Z$.
%[[Richard likes above much better. given $Z$.]]
\end{lem}

\noindent{\bf Proof.}
Because of Lemma~\ref{Qok} it suffices to prove that for all
$x\in\left[0,1\right]$, \linebreak
$P\left(F_{Y|Z}\left(Y\right)\leq x|Z\right)=x$ a.s.
This can be done as follows. Define
\begin{equation*}
F_{Y|Z}^{-1}\left(x+\right)=\sup\left\{y: F_{Y|Z}(y)\leq x\right\}.
\end{equation*}
Then
\begin{equation*}
P\left(F_{Y|Z}\left(Y\right)\leq x|Z=z\right)
=P\bigl(Y\leq F_{Y|Z}^{-1}\left(x+\right)|Z=z\bigr),
\end{equation*}
since $F_{Y|Z}(Y)\leq x$ implies that $Y\leq \sup\left\{y:
F_{Y|Z}(y)\leq x\right\}=F_{Y|Z}^{-1}\left(x+\right)$ and since $Y\leq
F_{Y|Z}^{-1}\left(x+\right)=\sup\left\{y: F_{Y|Z}(y)\leq x\right\}$
implies that $F_{Y|Z}(Y)\leq x$ by continuity of $F_{Y|Z}$. Hence
\begin{eqnarray*}
P\left(F_{Y|Z}\left(Y\right)\leq x|Z=z\right)
&=&E\Bigl[1_{\left\{Y\leq F_{Y|Z}^{-1}\left(x+\right)\right\}}
\,\big|Z=z\Bigr]\\
&=&\int_{\left(y,z'\right):y\leq F_{Y|Z=z'}^{-1}\left(x+\right)}
P_{\left(Y,Z\right)|Z=z}\left(dy,dz'\right)\\
&=&\int_{\left(y,z'\right):y\leq F_{Y|Z=z'}^{-1}\left(x+\right)}
\delta_{z,z'}F_{Y|Z=z}\left(dy\right) dz'\\
&=&\int_{y:y\leq F_{Y|Z=z}^{-1}\left(x+\right)}
F_{Y|Z=z}\left(dy\right)\\
&=& F_{Y|Z}\bigl(F_{Y|Z}^{-1}\left(x+\right)\bigr)
=x \hspace{0.5cm} \;\mbox{\rm a.s.},
\end{eqnarray*}
where Lemma~\ref{exp} is used in the second line, Lemma~\ref{concz} in
the third line, and continuity of $F_{Y|Z}$ in the last line.  \hfill
$\Box$

\begin{lem} \label{gZ} Suppose that $X$ is uniformly distributed on
$\left[0,1\right]$ and independent of $Z$ and that $F_{Y|Z}$ is a
conditional distribution function of $Y$ given $Z$. Then
$F_{Y|Z}^{-1}\left(X\right)$ has conditional distribution function
$F_{Y|Z}$ given $Z$.
%[[Richard likes above much better than given $Z$.]]
\end{lem}

\noindent{\bf Proof.}
Because of Lemma~\ref{Qok} it suffices to prove that for
all $s$, \linebreak
$P\bigl(F_{Y|Z}^{-1}\left(X\right)\leq
s|Z\bigr)=F_{Y|Z}\left(s\right)$ a.s. This can be done
as follows:
\begin{eqnarray*} P\bigl(F_{Y|Z}^{-1}\left(X\right)\leq s|Z=z\bigr)
&=&P\bigl(F^{-1}_{Y|Z}(X)\leq F^{-1}_{Y|Z}\circ F_{Y|Z}(s)|Z=z\bigr)\\
&=&P\bigl(X\leq F_{Y|Z}\left(s\right)|Z=z\bigr)\\
&=&E\Bigl[1_{\left\{X\leq
F_{Y|Z}\left(s\right)\right\}}\,\big|Z=z\Bigr]\\
&=&\int_{\left(x',z'\right):x'\leq F_{Y|Z=z'}\left(s\right)}
P_{\left(X,Z\right)|Z=z}\left(dx',dz'\right)\\
&=&\int_{\left(x',z'\right):x'\leq F_{Y|Z=z'}\left(s\right)}
\delta_{z,z'} F_{X|Z=z}\left(dx'\right) dz'\\
&=&\int_{x':x'\leq F_{Y|Z=z}\left(s\right)} F_{X|Z=z}\left(dx'\right)\\
&=&F_{Y|Z}\left(s\right) \;\;\;\;\;\;\;\;\; \;\mbox{\rm a.s.}.
\end{eqnarray*}
In the second line it is used that if $X\leq F_{Y|Z}(s)$ then also
$F^{-1}_{Y|Z}(X)\leq F^{-1}_{Y|Z}\circ F_{Y|Z}(s)$, and moreover that
if $X>F_{Y|Z}(s)$ then also, since $F_{Y|Z}(s)$ is in the range of
$F_{Y|Z}$ and conditional distribution functions are right continuous,
$F^{-1}_{Y|Z}(X)> F^{-1}_{Y|Z}\circ F_{Y|Z}(s)$. In the fourth line I
use Lemma~\ref{exp}, in the fifth line Lemma~\ref{concz} is used, and in
the last line it is used that $X$ is uniformly distributed on
$\left[0,1\right]$ given $Z$.  \hfill $\Box$

\begin{lem} \label{anpm}
If $X$ is a random variable taking values in ${\mathbb R}$ and for every
bounded Lipschitz continuous function $f:{\mathbb R}\rightarrow {\mathbb R}$
\begin{equation*}
E\left[f\left(X\right)|Z\right]=E\left[f\left(Y\right)|Z\right] \;\mbox{\rm a.s.}
\end{equation*}
then $X$ has the same conditional distribution as $Y$ given $Z$.
\end{lem}

\noindent{\bf Proof.}
Because of Lemma~\ref{Qok} it suffices to show that for every
$x\in {\mathbb R}$, \linebreak
$P\left(X\leq x|Z\right)=P\left(Y\leq x|Z\right)$ a.s.

Analogously to a proof of the Portmanteau Lemma, define
\begin{equation*}
f_m\left(y\right)=m \;{\rm
d}\left(y,\left(-\infty,x\right]\right)\wedge 1\end{equation*}
for $m=1,2,\ldots$.
Then $0\leq f_m\uparrow 1_{\left(x,\infty\right)}$ as $m\rightarrow
\infty$ and $f_m$ is bounded
and Lipschitz, so that
$E\left[f_m\left(X\right)|Z\right]=E\left[f_m\left(Y\right)|Z\right]$
a.s.
The remaining part is straightforward:
\begin{eqnarray*}
P\left(X\leq x|Z\right)
&=&E\left[1_{\left(-\infty,x\right]}\left(X\right)|Z\right]\\
&=&E\left[1-1_{\left(x,\infty\right)}\left(X\right)|Z\right]\\
&=&1-E\left[1_{\left(x,\infty\right)}\left(X\right)|Z\right] \;\mbox{\rm a.s.}\\
&=&1-\lim_{m\rightarrow \infty}E\left[f_m\left(X\right)|Z\right] \;\mbox{\rm a.s.}\\
&=&1-\lim_{m\rightarrow \infty}E\left[f_m\left(Y\right)|Z\right] \;\mbox{\rm a.s.}\\
&=&P\left(Y\leq x|Z\right) \;\mbox{\rm a.s.},
\end{eqnarray*}
where in the fourth line the conditional Monotone Convergence Theorem
(see e.g.\ \cite{FaZw}) is used.
\hfill $\Box$

$ $

\noindent{\bf Proof of Lemma~\ref{Zt-}.}
Remark that $\sigma\bigl(\overline{Z}_t^{(n)}\bigr)$ is increasing in $n$,
and that for $t$ on the infinite grid
$\sigma\bigl(\cup_{n=1}^{\infty}\overline{Z}_t^{(n)}\bigr)=\sigma\left(\overline{Z}_t\right)$ and for
$t$ not on the infinite grid
$\sigma\bigl(\cup_{n=1}^{\infty}\overline{Z}_t^{(n)}\bigr)=
\sigma\left(\overline{Z}_{t-}\right)$, where
$\overline{Z}_{t-}=\left(Z(s):s<t\right)$.  But the probability
that $Z$ jumps at time $t$ is equal to zero. Therefore
$E\left[X|\overline{Z}_t\right]=E\left[X|\overline{Z}_{t-}\right]$ a.s.: any version of
$E\left[X|\overline{Z}_{t-}\right]$ is a version of $E\left[X|\overline{Z}_t\right]$.
This can be seen as follows. $E\left[X|\overline{Z}_{t-}\right]$ is trivially
$\sigma\left(\overline{Z}_t\right)$-measurable. So it still has to be checked
that for any measurable $f:\overline{{\cal Z}}_t\rightarrow
\mathbb{R}$ for which $E\bigl(\bigl|X
f\left(\overline{Z}_t\right)\bigr|\bigr)<\infty$, $E \left(X
f\left(\overline{Z}_t\right)\right)= E\left( E\left[X|\overline{Z}_{t-}\right]
f\left(\overline{Z}_t\right)\right)$. So let $f$ with $E\bigl(\bigl|X
f\left(\overline{Z}_t\right)\bigr|\bigr)<\infty$ be given. Define
$g:\overline{{\cal Z}}_{t-}\rightarrow \overline{{\cal Z}}_{t}$ as the
``continuous'' extension:
\begin{equation*}g\left(\overline{z}_{t-}\right)(s)=\left\{\begin{array}{ll}
z(s) & {\rm if}\;s<t\\
\lim_{u\uparrow t} z(u)&{\rm if} \;s=t.\end{array}\right.
\end{equation*}
Then
\begin{eqnarray*} E \left(X f\left(\overline{Z}_t\right)\right)
&=& E\left( X f\left(g\left(\overline{Z}_{t-}\right)\right)\right)\\
&=&  E\left( E\left[X|\overline{Z}_{t-}\right]
f\left(g\left(\overline{Z}_{t-}\right)\right)\right)\\
&=&  E\left( E\left[X|\overline{Z}_{t-}\right] f\left(\overline{Z}_{t}\right)\right),
\end{eqnarray*}
where in the first and the last line it is used that the probability that
$Z$ jumps at time $t$ is equal to zero.  Therefore the conditional
expectation of $X$ given $\overline{Z}_t$ is almost surely equal to
the conditional expectation of $X$ given
$\sigma\bigl(\cup_{n=1}^{\infty}\overline{Z}_t^{(n)}\bigr)$.  \hfill
$\Box$

\section{Web-Appendix: A corollary of the Local Inverse Function Theorem}\label{LIFT}

{\bf Continuation of the proof of Lemma~\ref{quotc}.}  It is easy to
see that $\phi$ is differentiable at $(0,y_0)$ with non-singular
derivative. Therefore, the Local Inverse Function Theorem implies that
there exists an open neighbourhood $V_{h_0,y_0}$ of
$\left(h_0,y_0\right)$ such that $W=\phi\left(V_{h_0,y_0}\right)$ is
open and $\left.\phi\right|_{V_{h_0,y_0}}:V_{h_0,y_0}\rightarrow W$ is
a $C^1$-diffeomorphism. Hence $\phi^{-1}$ exists and is $C^1$.

Notice that $\phi^{-1}\left(h,x\right)$ must have the form
$\left(h,y\right)$ with $y$ satisfying $F_h\left(y\right)=x$. For
$\left(h,x\right)\in W$ such $y$ is unique, since all $F_h$ are
non-decreasing by assumption and $F_h'\left(y\right)$ is non-zero on
$V_{h_0,y_0}$. Thus $F_h^{-1}\left(x\right)$ is well-defined on $W$,
and it follows that
\begin{equation*}
\phi^{-1}\left(h,x\right)=\left(h,F^{-1}_h\left(x\right)\right).
\end{equation*}

Both $\phi$ and $\phi^{-1}$ are $C^1$, so the chain
rule can be applied to calculate
\begin{eqnarray*} \left(\begin{array}{ll} 1 & 0\\
                             0 & 1\end{array}\right)
&=&D\left(\phi \phi^{-1}\right)\left(h,x\right)\\
&=&\left(D\phi\right)\left(\phi^{-1}\left(h,x\right)\right)
\cdot\left(D\phi^{-1}\right)\left(h,x\right)\\
&=&\left(\begin{array}{ll} 1 & 0\\
   \frac{\partial}{\partial h} F_h\left(y\right)   &
F_h'\left(y\right)\end{array}\right) \cdot
\left(\begin{array}{ll} 1 & 0\\
   \frac{\partial}{\partial h} F^{-1}_h\left(x\right)   &
\left(F^{-1}_h\right)'\left(x\right)\end{array}\right),
\end{eqnarray*}
with $y=F_h^{-1}\left(x\right)$. Lemma~\ref{quotc} follows by comparing the
bottom left entries of the matrices on the left- and right hand side
of this equation.  \hfill $\Box$

\section{Web-Appendix: Lipschitz continuity and differentiability}
%Lemma's about
\label{Lips}

The following lemma can be useful for proving Lipschitz continuity of
quotients of functions.

\begin{lem} \label{Lip}
Suppose that $f$ and $g$ are functions from ${\mathbb R}$ to ${\mathbb R}$
which are Lipschitz continuous with Lipschitz
constants $L_f$ resp.\ $L_g$. Suppose furthermore that $g\geq \varepsilon$
for some $\varepsilon>0$ and $\left|f\right|\leq C$ for some $C>0$.
Then $f/g$ is
Lipschitz continuous with Lipschitz constant
e.g.\
$L_f/\varepsilon+C\, L_g/\varepsilon^2$.
\end{lem}

\noindent{\bf Proof.}
\begin{eqnarray*}
\left|\frac{f\left(x_1\right)}{g\left(x_1\right)}-\frac{f\left(x_2\right)}{g\left(x_2\right)}\right|
&\leq&\left|\frac{f\left(x_1\right)}{g\left(x_1\right)}-\frac{f\left(x_2\right)}{g\left(x_1\right)}\right|
+\left|\frac{f\left(x_2\right)}{g\left(x_1\right)}-\frac{f\left(x_2\right)}{g\left(x_2\right)}\right|\\
&=&\left|\frac{1}{g\left(x_1\right)}\right|\left|f\left(x_1\right)-f\left(x_2\right)\right|
+\left|\frac{f\left(x_2\right)}{g\left(x_1\right)g\left(x_2\right)}\right|\left|g\left(x_1\right)-g\left(x_2\right)\right|\\
&\leq& \frac{1}{\varepsilon} L_f\left|x_1-x_2\right|+
\frac{C}{\varepsilon^2}L_g\left|x_1-x_2\right|.\hspace*{5.8cm} \Box
\end{eqnarray*}

The next lemma deals with a continuous function $f$ on a
closed interval which is continuously differentiable on the interior
of that interval. If $f'$ can be continuously extended to the closed
interval, then $f$ is continuously differentiable on the closed
interval.

\begin{lem} \label{difext}
Suppose that $f$ is continuous on $\left[t_1,t_2\right]$ and $f$ is
continuously differentiable on $\left(t_1,t_2\right)$. Suppose
furthermore that $f'$ has a continuous extension to
$\left[t_1,t_2\right]$. Then $f$ is differentiable from the right at
$t_1$ with derivative $\lim_{t\downarrow t_1}f'\left(t\right)$ and
differentiable from the left at $t_2$ with derivative $\lim_{t\uparrow
t_2}f'\left(t\right)$.
\end{lem}

\noindent{\bf Proof.}
I just prove the statements for $t_1$; the proof for $t_2$ is
similar. Define $g\left(t\right)=f\left(t_1\right)+\int_{t_1}^t
f'\left(x\right)dx$. Then $g$ is continuous and continuously
differentiable on $\left[t_1,t_2\right)$ with derivative $f'(t)$ on
$\left(t_1,t_2\right)$ and $\lim_{t\downarrow t_1}f'\left(t_1\right)$
at $t_1$. It suffices to show that $f=g$ on $\left[t_1,t_2\right)$, since
then $f$ has the same properties as $g$ on $\left[t_1,t_2\right)$. To do
this, remark first that $f-g$ is constant on $\left(t_1,t_2\right)$
since it is differentiable there with derivative $0$. Because $f-g$ is
continuous on $\left[t_1,t_2\right)$, $f-g$ is also constant on
$\left[t_1,t_2\right)$. $\left(f-g\right)\left(t_1\right)=0$. Thus $f=g$ on
$\left[t_1,t_2\right)$.
\hfill $\Box$

\section{Web-Appendix: Some theory about differential equations}
\label{difeqapp}

\begin{thm} \label{isun}
Suppose that a function $D\left(y,t;\overline{Z}_t\right)$ satisfies
\begin{enumerate}[a)]
\item \emph{(continuity between the jump times of $Z$).} If
$Z$ does not jump in $\left(t_1,t_2\right)$ then
$D\left(y,t;\overline{Z}_t\right)$ is continuous in $\left(y,t\right)$ on
$\left[t_1,t_2\right)$ and can be continuously extended to
$\left[t_1,t_2\right]$.
\item \emph{(Lipschitz continuity).} For
each $\omega\in\Omega$ there exists a constant $L\left(\omega\right)$
such that
\begin{equation*}
\left|D\left(y,t;\overline{Z}_t\right)-D\left(z,t;\overline{Z}_t\right)\right|\leq
L\left(\omega\right)\left|y-z\right|
\end{equation*}
for all $t\in\left[0,\tau\right]$ and all $y,z$.
\end{enumerate}
Suppose furthermore that for each $\omega\in\Omega$ there are no more
than finitely many jump times of $Z$. Then, for each
$t_0\in\left[0,\tau\right]$ and $y_0\in\mathbb{R}$, there is a unique
continuous solution $x\left(t;t_0,y_0\right)$ to
\begin{equation*}
dx(t)/dt=D\left(x(t),t;\overline{Z}_t\right)
\end{equation*}
with boundary condition $x\left(t_0\right)=y_0$ and this solution is
defined for all $t\in\left[0,\tau\right]$.
\end{thm}
\noindent This theorem follows from well-known results about
differential equations, see e.g.\ \citet{anDe} Chapter~2.

For the next theorem we also refer to \citet{anDe} Chapter 2. It is a
consequence of Gronwall's lemma.
\begin{thm} \label{difb} Suppose that $I$ is an open or closed
interval in ${\mathbb R}$,
$f:I\times {\mathbb R}^n\rightarrow {\mathbb R}^n$ is continuous and
$C:I\rightarrow\left[0,\infty\right)$ is continuous, and suppose that
\begin{equation}
\parallel f\left(x,y\right)-f\left(x,z\right)\parallel\leq
C\left(x\right) \parallel y-z\parallel
\label{Lc2}\end{equation}
for all $x\in I$ and $y,z\in {\mathbb R}^n$.
Then, for every $x_0\in I$ and $y_0\in\mathbb{R}$,
there is a unique solution $y\left(x\right)$ of
$y'\left(x\right)=f\left(x,y(x)\right)$ with
$y\left(x_0\right)=y_0$, and this solution is defined
for all $x\in I$. If
$g:I\times {\mathbb R}^n\rightarrow {\mathbb R}^n$ is continuous and
$z:I\rightarrow {\mathbb R}^n$ is a
solution of
$z'\left(x\right)=g\left(x,z(x)\right)$ then
\begin{eqnarray*} \lefteqn{\parallel y\left(x\right)
-z\left(x\right)\parallel}\\
&\leq&
e^{\int_{x_0}^x C\left(\xi\right)d\xi}
\parallel y\left(x_0\right)-z\left(x_0\right)\parallel
+\int_{x_0}^x e^{\int_{\xi}^x C\left(\eta\right)d\eta}
\parallel f\left(\xi,z\left(\xi\right)\right)-
g\left(\xi,z\left(\xi\right)\right)\parallel d\xi
\end{eqnarray*}
for all $x,x_0\in I$ with $x_0\leq x$.
\end{thm}
In \citet{anDe} the interval is always an open interval, but as is
generally known this can be overcome by extending both $f$ and $g$
outside the closed interval $I$ by taking the values at the boundary
of $I$. This preserves the Lipschitz- and continuity conditions.
Existence and uniqueness on all of finitely many intervals implies
global existence and uniqueness; this is the way one often applies
this theorem.

This article is about a differential equation with end condition at
$\tau$, so interested lies in $x,x_0$ with $x\leq x_0$.
The following corollary can be used.
\begin{cor}
\label{difbc1}
Suppose that the conditions of Theorem~\ref{difb} are satisfied.
Then, for every $x_0\in I$ and $y_0\in\mathbb{R}^n$,
there is a unique solution $y\left(x\right)$ of
$y'\left(x\right)=f\left(x,y(x)\right)$ with
$y\left(x_0\right)=y_0$, and this solution is defined
for all $x\in I$. If
$g:I\times {\mathbb R}^n\rightarrow {\mathbb R}^n$ is continuous and
$z:I\rightarrow {\mathbb R}^n$ is a
solution of
$z'\left(x\right)=g\left(x,z(x)\right)$ then
\begin{eqnarray*} \lefteqn{\parallel y\left(x\right)
-z\left(x\right)\parallel}\\
&\leq&
e^{\int_{x}^{x_0} C\left(s\right)ds}
\parallel y\left(x_0\right)-z\left(x_0\right)\parallel
+\int_{x}^{x_0} e^{\int_{x}^s C\left(\eta\right)d\eta}
\parallel f\left(s,z\left(s\right)\right)-
g\left(s,z\left(s\right)\right)\parallel ds
\end{eqnarray*}
for all $x,x_0$ with $x\leq x_0$.
\end{cor}

\noindent{\bf Proof.} Define
$\tilde{y}\left(t\right)=y\left(x_0-t\right)$. Then
\begin{eqnarray*}
\tilde{y}'\left(t\right)&=&\frac{\partial}{\partial
t}y\left(x_0-t\right)\\
&=&-y'\left(x_0-t\right)\\
&=&-f\left(x_0-t,y\left(x_0-t\right)\right)\\
&=&-f\left(x_0-t,\tilde{y}\left(t\right)\right)\\
&=&\tilde{f}\left(t,\tilde{y}\left(t\right)\right)
\end{eqnarray*}
where $\tilde{f}\left(t,y\right)=-f\left(x_0-t,y\right)$. So
$\tilde{y}\left(t\right)=y\left(x_0-t\right)$ is a solution of the
differential equation
$\tilde{y}'\left(t\right)=\tilde{f}\left(t,y\left(t\right)\right)$
with boundary condition
$\tilde{y}\left(0\right)=y\left(x_0\right)=y_0$.
Applying Theorem~\ref{difb} on $\tilde{y}$ concludes the proof, as follows.
\begin{eqnarray*}
\parallel y(x)-z(x)\parallel
&=& \parallel y(x-x_0+x_0)-z(x-x_0+x_0)\parallel\\
&=& \parallel y(x_0-(x_0-x))-z(x_0-(x_0-x))\parallel\\
&=&\parallel \tilde{y}(x_0-x) - \tilde{z}(x_0-x)\parallel\\
&=&\parallel \tilde{y}(t)-\tilde{z}(t)\parallel
\end{eqnarray*}
with $t=x_0-x\geq 0.$ Notice that since because of
equation~(\ref{Lc2}),
\begin{equation*}
\parallel \tilde{f}\left(t,y\right)-\tilde{f}\left(t,z\right)\parallel\leq
C\left(x_0-t\right) \parallel y-z\parallel=:\tilde{C}(t)
\parallel y-z \parallel,
\end{equation*}
with $\tilde{C}(t)=C(x_0-t)$.
Hence Theorem~\ref{difb} implies that
\begin{eqnarray*}
\parallel y(x)-z(x)\parallel
&\leq& e^{\int_0^t\tilde{C}(\xi)d\xi}
\parallel\tilde{y}(0)-\tilde{z}(0)\parallel
+\int_0^t e^{\int_0^t \tilde{C}\left(\eta\right)d\eta}
\parallel\tilde{f}\left(\xi,\tilde{z}\left(\xi\right)\right)
-\tilde{g}\left(\xi,\tilde{z}(\xi)\right)\parallel d\xi\\
&=&e^{\int_0^t C(x_0-\xi)d\xi}\parallel y(x_0-0)-z(x_0-0)\parallel\\
&&+\int_0^t e^{\int_0^t C\left(x_0-\eta\right)d\eta}
\parallel\tilde{f}\left(\xi,\tilde{z}\left(\xi\right)\right)
-\tilde{g}\left(\xi,\tilde{z}(\xi)\right)\parallel d\xi.
\end{eqnarray*}
For the first term a change of variables is done; $\xi$ from $0$ to
$t$, define $s=x_0-\xi$; $d\xi=-ds$. $0\leq\xi\leq t$; $s$ from $x_0-0$
to $x_0-t=x_0-(x_0-x)=x$. Therefore, the first term is equal to
\begin{equation*}
e^{-\int_{x_0}^x C\left(s\right)ds}\parallel y(x_0)-z(x_0)\parallel
= e^{\int_{x}^{x_0} C\left(s\right)ds}\parallel y(x_0)-z(x_0)\parallel.
\end{equation*}
For the second term similar changes of variables can be done,
resulting in Corollary~\ref{difbc1}.
\hfill $\Box$\\

\noindent{\bf Proof of Theorem~\ref{difbc2}.} Write $I=\left[x_1,x_2\right]$.
In order to apply Theorem~\ref{difb}
define an extension $\tilde{f}:{\mathbb R}\times {\mathbb
R}\rightarrow {\mathbb R}$ of $f$ as follows:
\begin{equation*} \tilde{f}\left(x,y\right)=\left\{
\begin{array}{ll}
f\left(x,y\right)&{\rm if} \;
\left(x,y\right)\in I\times \left[y_1,y_2\right]\\
f\left(x_1,y\right)&{\rm if} \;
\left(x,y\right)\in\left(-\infty,x_1\right]\times\left[y_1,y_2\right]\\
f\left(x_2,y\right)&{\rm if} \;
\left(x,y\right)\in \left(x_2,\infty\right)\times \left[y_1,y_2\right]\\
\tilde{f}\left(x,y_1\right)&{\rm if}\; y<y_1\\
\tilde{f}\left(x,y_2\right)&{\rm if}\; y>y_2.
%0&{\rm if}\; y\notin\left[y_1,y_2\right].
\end{array}\right.
\end{equation*}
If there exists a unique solution of the differential
equation with $\tilde{f}$ and this solution stays in
$\left[y_1,y_2\right]$, then this
solution is also the unique continuous solution of the differential
equation with $f$.

On the differential equation with $\tilde{f}$, Theorem~\ref{difb} will
be applied. $\tilde{f}$ is continuous on ${\mathbb R}\times{\mathbb
R}$ because $f$ is continuous on $I\times\left[y_1,y_2\right]$ and
$f\left(x,y_1\right)=0=f\left(x,y_2\right)$ for every $x\in I$.  Also
there exists a continuous $\tilde{C}$ satisfying equation~(\ref{Lc2}):
define $\tilde{C}$ as an extension of $C$ as follows:
\begin{equation*}\tilde{C}\left(x\right)=\left\{
\begin{array}{ll}
C\left(x\right)&{\rm if} \; x\in I\\
C\left(x_1\right)&{\rm if} \; x\in\left(-\infty,x_1\right)\\
C\left(x_2\right)&{\rm if} \; x\in \left(x_2,\infty\right).
\end{array}\right.
\end{equation*}
That this $\tilde{C}$ satisfies equation~(\ref{Lc2}) can easily be
checked by first considering $x\in\left[x_1,x_2\right]$ and reducing
different $x$ to $x_1$ and $x_2$.

Thus Theorem~\ref{difb} implies that there is a unique solution of the
differential equation with $\tilde{f}$.  That the solution stays in
$\left[y_1,y_2\right]$ is clear from the fact that $\tilde{f}=0$ for
$y\in\{y_1,y_2\}$ and the fact that the solution is unique.

Since $g$ can be extended the same way as $f$ and $z$ stays in
$\left[y_1,y_2\right]$ by assumption, the bound for $\left|
y\left(x\right)-z\left(x\right)\right|$ given by Theorem~\ref{difb}
also holds here.  The bound of Theorem~\ref{difbc2} follows with the same
reasoning from Corollary~\ref{difbc1}.  This finishes
the proof.  \hfill $\Box$

\begin{cor}
\label{difbc3} Suppose that $I=\left[x_1,x_2\right]\subset
\left[0,y_2\right]$ is a closed
interval in ${\mathbb R}$,
$f:\left\{\left(x,y\right)\in I\times \left[0,y_2\right]:y\geq
x\right\}\rightarrow {\mathbb R}$ is continuous
with for all $x\in I$, $f\left(x,y_2\right)=0$ and
$f\left(x,x\right)\leq 1$ and
$C:I\rightarrow\left[0,\infty\right)$ is continuous, and suppose that
\begin{equation*}
\left| f\left(x,y\right)-f\left(x,z\right)\right|\leq
C\left(x\right) \left| y-z\right|
\end{equation*}
for all $x\in I$ and $y,z\in \left[x,y_2\right]$.
Then for every $y_0\in\left[x_2,y_2\right]$ there is a unique solution
$y\left(x\right)$ of
$y'\left(x\right)=f\left(x,y(x)\right)$ with
final condition $y\left(x_2\right)=y_0$, and this solution is defined
for all $x\in I$. Furthermore $y\left(x\right)\in\left[x,y_2\right]$
for all $x\in I$. Suppose that
$g:\left\{\left(x,y\right)\in I\times \left[0,y_2\right]:y\geq
x\right\}\rightarrow {\mathbb R}$ is continuous and
$z:I\rightarrow \left[0,y_2\right]$ is a
solution of
$z'\left(x\right)=g\left(x,z(x)\right)$
with $z\left(x\right)\geq x$ then
\begin{eqnarray*} \lefteqn{\left| y\left(x\right)
-z\left(x\right)\right|}\\
&&\hspace{0.5cm}\leq
e^{\int_{x}^{x_2} C\left(s\right)ds}
\left| y\left(x_2\right)-z\left(x_2\right)\right|
+\int_{x}^{x_2} e^{\int_{x}^s C\left(\eta\right)d\eta}
\left| f\left(s,z\left(s\right)\right)-
g\left(s,z\left(s\right)\right)\right| ds
\end{eqnarray*}
for all $x\in I$.
\end{cor}

\noindent{\bf Proof.}
This can be proved the same way as Corollary~\ref{difbc2}, if one
defines
\begin{equation*} \tilde{f}\left(x,y\right)=\left\{
\begin{array}{ll}
f\left(x,y\right)&{\rm if} \;
\left(x,y\right)\in I\times \left[0,y_2\right]:y\geq x\\
f\left(x,y_2\right)&{\rm if} \;
x\in I\;{\rm and}\;y>y_2\\
f\left(x,x\right)&{\rm if} \;
x\in I\;{\rm and}\;y<x\\
\tilde{f}\left(x_1,y\right)&{\rm if}\; x<x_1\\
\tilde{f}\left(x_2,y\right)&{\rm if} \; x>x_2.
\end{array}\right.
\end{equation*}
Remark that the solution $y\left(x\right)$ stays in
$\left(x,y\right)\in I\times \left[0,y_2\right]:y\geq x$ for $x\in I$
since $f\left(x,x\right)\leq 1$ and $f\left(x,y_2\right)=0$.
\hfill $\Box$

$ $

\noindent Remark that if it is not known whether $f\left(x,x\right)\leq 1$ but
it is known that $f$ is continuous in $\left(x,y\right)$ and Lipschitz
continuous in $y$ on the set mentioned in Corollary~\ref{difbc3}, then
the proof above shows that if a solution $y\left(x\right)$ exists for
which $y\left(x\right)\geq x$ for $x\leq x_2$ then this solution is
unique.

\section{Web-Appendix: Convergence Theorems}
\label{Conv}

A lemma with a corollary:
\begin{lem} \label{plugl} Suppose that the random functions
$f^\omega:\left[y_1,y_2\right]\rightarrow {\mathbb R}$ and
$f_n^\omega:\left[y_1,y_2\right]\rightarrow {\mathbb R}$ ($n=1,2,\ldots$)
are `asymptotically uniformly equicontinuous with probability one',
i.e.\ there exists $\Omega'\subset
\Omega$ with $P\left(\Omega'\right)=1$ such that for all
$\omega\in\Omega'$: $\forall \varepsilon>0$ $\exists \delta>0$ $\exists
N$: $\forall n\geq N$:
\begin{equation*} \left|y-z\right|<\delta \; \Rightarrow
\left\{ \begin{array}{l} \left|f_n^\omega\left(y\right)-
f_n^\omega\left(z\right)\right| <\varepsilon\\
\left|f^\omega\left(y\right)-
f^\omega\left(z\right)\right| <\varepsilon. \end{array} \right.
\end{equation*}
Suppose furthermore that for all $y\in\left(y_1,y_2\right)\cap {\mathbb
Q}$,
$f_n^\omega\left(y\right)\rightarrow f^\omega\left(y\right)$ a.s.
Then
\begin{equation*}
\sup_{y\in\left[y_1,y_2\right]}
\left|f_n^\omega\left(y\right)-f^\omega\left(y\right)\right|
\rightarrow 0 \;\;\;\;\;\; \;\mbox{\rm a.s.}\end{equation*}
\end{lem}
Remark: the regularity condition for Lemma~\ref{plugl} is e.g.\
satisfied if there is
a Lipschitz constant $L$ such that all $f^\omega$ and $f^\omega_n$ are
Lipschitz
continuous in $y$ with Lipschitz constant $L$ (define
$\delta=\varepsilon/L$).

$ $

\noindent{\bf Proof.}
%of Lemma~\ref{plugl}.\\
Define $\Omega''=\left\{\omega:f^\omega_n\left(y\right)\rightarrow
f^\omega\left(y\right) \; \forall y\in{\mathbb Q}\cap
\left(y_1,y_2\right)\right\}$. Then $\Omega''$ has probability one
($\Omega$ minus countably many null sets). Define $\Omega_0=\Omega'\cap
\Omega''$. Then also $\Omega_0$ has probability one. We show that
for all $\omega\in\Omega_0$:
$\sup_{y\in\left[y_1,y_2\right]}\left|f_n^\omega\left(y\right)-f^\omega\left(y\right)\right|\rightarrow
0$.

Let $\omega\in\Omega_0$ and $\varepsilon>0$ be given. To show: there
exists an $N$ such that $\forall n\geq N$:\linebreak
$\sup_{y\in\left[y_1,y_2\right]}\left|f_n^\omega\left(y\right)-f^\omega\left(y\right)\right|<\varepsilon$.
Choose $N_1$ and $\delta>0$ such that for all $n\geq N_1$:
\begin{equation*} \left|y-z\right|<\delta \; \Rightarrow
\left\{ \begin{array}{l} \left|f_n^\omega\left(y\right)-
f_n^\omega\left(z\right)\right| <\varepsilon/3\\
\left|f^\omega\left(y\right)-
f^\omega\left(z\right)\right| <\varepsilon/3. \end{array} \right.
\end{equation*}
This is possible because $\omega\in\Omega'$. Next choose
$y^{\left(1\right)},\ldots,y^{\left(N_2\right)}\in {\mathbb Q}\cap
\left(y_1,y_2\right)$ such that for all
$y\in\left[y_1,y_2\right]$ there is a $y^{\left(i\right)}$ with
$|y-y^{\left(i\right)}|<\delta$. After this choose $N_3$
such that for all $n\geq N_3$:
\begin{equation*}\max_{1\leq i\leq
N_2}\bigl|f_n^\omega\bigl(y^{\left(i\right)}\bigr)-f^\omega\bigl(y^{\left(i\right)}\bigr)\bigr|
<\varepsilon/3.\end{equation*} This is possible because
$\omega\in\Omega''$ and the number of $y^{(i)}$'s is finite. Then for
$n\geq N=\max\{N_3,N_1\}$:
\begin{eqnarray*} \left|f_n^\omega\left(y\right)-f^\omega\left(y\right)\right|
&\leq& \min_{1\leq i \leq
N_2}\Bigl(\bigl|f_n^\omega\left(y\right)
-f_n^\omega\bigl(y^{\left(i\right)}\bigr)\bigr|+
\bigl|f^\omega\bigl(y^{\left(i\right)}\bigr)
-f^\omega\left(y\right)\bigr|\Bigr)\\
&&\hspace{1cm}+\max_{1\leq
i \leq N_2}
\bigl(\bigl|f_n^\omega\bigl(y^{\left(i\right)}\bigr)
-f^\omega\bigl(y^{\left(i\right)}\bigr)\bigr|\bigr)\\
&<&\varepsilon/3+\varepsilon/3+\varepsilon/3=\varepsilon.\hspace*{8cm}\Box
\end{eqnarray*}

\begin{cor} \label{plug} Under the conditions of Lemma~\ref{plugl}, if $X_n$
is a series of random variables with values in $\left[y_1,y_2\right]$, then
\begin{equation*}
\left|f_n^\omega\left(X_n\right)-f^\omega\left(X_n\right)\right|\rightarrow
0 \;\;\;\;\;\; \;\mbox{\rm a.s.}\end{equation*}
\end{cor}

%I include a well-known corollary of Lebesgue's Dominated Convergence
%Theorem, to have the precise conditions at hand.
%\begin{thm} \label{difint}
%Suppose that $\mu$ is a measure on a measurable space $\left({\cal
%X},{\cal A}\right)$.  Suppose that for all $t\in\left[a,b\right]$,
%$f_t$ is a measurable function on ${\cal X}$. Suppose furthermore that
%for every $t$, $\frac{\partial}{\partial t} f_t\left(\cdot\right)$
%exists. Suppose that there exists an integrable $g$ such that for all
%$t$: $\left|\frac{\partial}{\partial t}f_t\right|\leq g$ and that
%$f_{t_0}$ is integrable for some $t_0\in\left[a,b\right]$.  In that
%case $f_t$ is integrable for all $t$ and
%\begin{equation*} \frac{\partial}{\partial t} \int f_t d\mu = \int
%\frac{\partial}{\partial t} f_t d\mu
%.\end{equation*}
%\end{thm}

\addcontentsline{toc}{chapter}{Bibliography}
\bibliographystyle{chicago}
\bibliography{ref}

\end{document}